# Hard ball systems are completely hyperbolic

By Nándor Simányi and Domokos Szász*

**Abstract**

We consider the system of $N$ ($\geq 2$) elastically colliding hard balls with masses $m_1, \ldots, m_N$, radius $r$, moving uniformly in the flat torus $\mathbb{T}^\nu_L = \mathbb{R}^\nu / L \cdot \mathbb{Z}^\nu$, $\nu \geq 2$. It is proved here that the relevant Lyapunov exponents of the flow do not vanish for almost every $(N+1)$-tuple $(m_1, \ldots, m_N; L)$ of the outer geometric parameters.

## 1. Introduction

The proper mathematical formulation of Ludwig Boltzmann's ergodic hypothesis, which has incited so much interest and discussion in the last hundred years, is still not clear. For systems of elastic hard balls on a torus, however, Yakov Sinai, in 1963, [Sin(1963)] gave a stronger, and at the same time mathematically rigorous, version of Boltzmann's hypothesis: *The system of an arbitrarily fixed number $N$ of identical elastic hard balls moving in the $\nu$-torus $\mathbb{T}^\nu = \mathbb{R}^\nu / \mathbb{Z}^\nu$ ($\nu \geq 2$) is ergodic — of course, on the submanifold of the phase space specified by the trivial conservation laws.* Boltzmann used his ergodic hypothesis when laying down the foundations of statistical physics, and its various forms are still intensively used in modern statistical physics. The importance of Sinai's hypothesis for the theory of dynamical systems is stressed by the fact that the interaction of elastic hard balls defines the only physical system of an arbitrary number of particles in arbitrary dimension whose dynamical behaviour has been so far at least guessed — except for the completely integrable system of harmonic oscillators. (As to the history of Boltzmann's hypothesis, see the recent work [Sz(1996)].)

Sinai's hypothesis was partially based on the physical arguments of Krylov's 1942 thesis (cf. [K(1979)] and its afterword written by Ya. G. Sinai, [Sin(1979)]), where Krylov discovered that hard ball collisions provide effects

*Research supported by the Hungarian National Foundation for Scientific Research, grants OTKA-7275, OTKA-16425, and T-026176.

analogous to the hyperbolic behaviour of geodesic flows on compact manifolds of constant negative curvature, exploited so beautifully in the works of Hedlund, [He(1939)] and Hopf, [Ho(1939)].

The aim of the present paper is to establish that *hard ball systems are, indeed, fully hyperbolic*, i.e. all relevant Lyapunov exponents of these systems are nonzero almost everywhere. Our claim holds for typical $(N+1)$-tuples $(m_1, \ldots, m_N; L) \in \mathbb{R}_+^{N+1}$ of the outer geometric data of the system; that is — in contrast to earlier results — we do not require that particles have identical masses, though, on the other hand, we have to exclude a countable union of proper submanifolds of $(N+1)$-tuples $(m_1, \ldots, m_N; L)$ (which set is, in fact, very likely to be empty).

Full hyperbolicity combined with Katok-Strelcyn theory (see [K-S(1986)]) immediately provides that *the ergodic components of these systems are of positive measure*. Consequently, there are at most countably many of them, and, moreover, on each of them the system is K-mixing. Our methods so far do not give the expected global ergodicity of the systems considered (an additional hypothesis to provide that is formulated in Section 6).

The equality of the *radii* of the balls is not essential, but, for simplicity, it will be assumed throughout. For certain values of the radii — just think of the case when they are large — the phase space of our system decomposes into a finite union of different connected components, and these connected parts certainly belong to different ergodic components. Now, according to the wisdom of the ergodic hypothesis, these connected components are expected to be just *the* ergodic components of the system, and on each of them the system should also possess the Kolmogorov-mixing property.

Let us first specify the model and formulate our result.

Assume that, in general, a system of $N (\geq 2)$ balls, identified as $1, 2, \ldots, N$, of masses $m_1, \ldots, m_N$ and of radius $r > 0$ are given in $\mathbb{T}_L^\nu = \mathbb{R}^\nu / L \cdot \mathbb{Z}^\nu$, the $\nu$-dimensional cubic torus with sides $L$ ($\nu \geq 2$). Denote the phase point of the $i$th ball by $(q_i, v_i) \in \mathbb{T}_L^\nu \times \mathbb{R}^\nu$. A priori, the configuration space $\tilde{\mathbf{Q}}$ of the $N$ balls is a subset of $\mathbb{T}_L^{N \cdot \nu}$: from $\mathbb{T}_L^{N \cdot \nu}$ we cut out $\binom{N}{2}$ cylindric scatterers:

$$(1.1) \qquad \tilde{C}_{i,j} = \left\{ Q = (q_1, \ldots, q_N) \in \mathbb{T}_L^{N \cdot \nu} : \|q_i - q_j\| < 2r \right\},$$

$1 \leq i < j \leq N$, or in other words $\widetilde{\mathbf{Q}} := \mathbb{T}_L^{N\nu} \setminus \bigcup_{1 \leq i < j \leq N} \tilde{C}_{i,j}$. The energy $H = \frac{1}{2} \sum_1^N m_i v_i^2$ and the total momentum $P = \sum_1^N m_i v_i$ are first integrals of the motion. Thus, without loss of generality, we can assume that $H = \frac{1}{2}$, $P = 0$. (If $P \neq 0$, then the system has an additional conditionally periodic or periodic motion.) Now, for these values of $H$ and $P$, we define our dynamical system.

*Remark.* It is clear that actual values of $r$ and $L$ are not used, but merely their ratio $L/r$, is relevant for our model. Therefore, throughout this paper we fix the value of $r > 0$, and will only consider (later in Sections 3–4) the size $L$ of the torus as a variable.

In earlier works (cf. [S-Sz(1995)] and the references therein), where the masses were identical, one could and — to obtain ergodicity — had to fix a center of mass. For different masses — as observed in [S-W(1989)] — this is generically not possible and we shall follow a different approach.

The equivalence relation $\Psi$ over $\widetilde{\mathbf{Q}}$, defined by $Q \sim_\Psi Q^*$ *if and only if there exists an* $a \in \mathbb{T}_L^\nu$ *such that for every* $i \in [1, N]$, $q_i^* = q_i + a$ allows us to introduce $\mathbf{Q} := \widetilde{\mathbf{Q}}/\Psi$ (the equivalence relation means, in other words, that for the internal coordinates $q_i - q_j = q_i^* - q_j^*$ holds for every $i, j \in [1, N]$). The set $\mathbf{Q}$, a compact, flat Riemannian manifold with boundary will actually be the configuration space of our system, whereas its phase space will be the ellipsoid bundle $\mathbf{M} := \mathbf{Q} \times \mathcal{E}$, where $\mathcal{E}$ denotes the ellipsoid $\sum_{i=1}^N m_i (dq_i)^2 = 1$, $\sum_{i=1}^N m_i dq_i = 0$. Clearly, $d := \dim \mathbf{Q} = N\nu - \nu$, and $\dim \mathcal{E} = d - 1$. The well-known Liouville measure $\mu$ is invariant with respect to the evolution $S^\mathbb{R} := \{S^t : t \in \mathbb{R}\}$ of our dynamical system defined by elastic collisions of the balls of masses $m_1, \ldots, m_N$ and their uniform free motion. Here we have two remarks:

(i) The collision laws for a pair of balls with different masses are well-known from mechanics and will also be reproduced in the equation (3.9);

(ii) The dynamics can, indeed, be defined for $\mu$ — a.e. phase point; see the corresponding references in Section 2.

The dynamical system $(\mathbf{M}, S^\mathbb{R}, \mu)_{\vec{m}, L}$ is called *the standard billiard ball system with the outer geometric parameters* $(\vec{m}; L) \in \mathbb{R}_+^{N+1}$.

Denote by $\tilde{R}_0 := \tilde{R}_0(N, \nu, L)$ the interval of those values of $r > 0$, for which the interior $\operatorname{Int} \mathbf{M} = \operatorname{Int} \mathbf{Q} \times \mathcal{E}$ of the phase space of the standard billiard ball flow is connected. From $\tilde{R}_0$ we should also exclude an at most countable number of values of $r$ where the nondegeneracy condition of [B-F-K(1998)] (see the Definition on page 697) fails. The resulting set will be denoted by $R_0$.

The basic result of our paper is the following:

MAIN THEOREM. *For $N \geq 2, \nu \geq 2$ and $r \in R_0$, none of the relevant Lyapunov exponents of the standard billiard ball system $(\mathbf{M}, \{S^\mathbb{R}\}, \mu)_{\vec{m}, L}$ vanishes — apart from a countable union of proper analytic submanifolds of the outer geometric parameters* $(\vec{m}; L) \in \mathbb{R}_+^{N+1}$.

An interesting consequence of this theorem is the following

COROLLARY. *For the good set of typical outer geometric parameters $(\vec{m}; L)$ the ergodic components of the system have positive measure and on each of them the standard billiard ball flow has the* K-*property.*

*Note* 1. In general, if $r \notin R_0$, then Int **M** decomposes into a finite number of connected components. Our results extend to these cases, too.

*Note* 2. As it will be seen in Section 2 (see Lemma 2.1) our system is isomorphic to a *semi-dispersing billiard*, and, what is more, is a semi-dispersing billiard in a weakly generalized sense, and, consequently, we can use the theory of semi-dispersing billiards. As has been proved in recent manuscripts by N. I. Chernov and C. Haskell, [C-H(1996)] on one hand, and by D. Ornstein and B. Weiss, [O-W(1998)] on the other hand, the K-mixing property of a semi-dispersing billiard flow on a positive ergodic component actually implies its Bernoulli property, as well.

*Note* 3. As to the basic results (and history) concerning the Boltzmann and Sinai hypotheses we refer to the recent survey [Sz(1996)].

The basic notion in the theory of semi-dispersing billiards is that of the *sufficiency* of a phase point or, equivalently, of its orbit. The conceptual importance of sufficiency can be explained as follows (for a technical introduction and our prerequisites, see Section 2): In a suitably small neighbourhood of a (typical) phase point of a dispersing billiard the system is *hyperbolic*; i.e. its relevant Lyapunov exponents are not zero. For a semi-dispersing billiard the same property is guaranteed for sufficient points only! Physically speaking, a phase point is sufficient if its trajectory encounters in its history all possible degrees of freedom of the system.

Our proof of the full hyperbolicity is the first major step in the basic strategy for establishing global ergodicity of semi-dispersing billiards as was suggested in our series of works with A. Krámli (see the references), and, in fact, some elements of that approach are also used in this paper. In the sense of our strategy for proving global ergodicity, initiated in [K-S-Sz(1991)] and explained in the introductions of [K-S-Sz(1992)] and of [Sim(1992)-I] , there are two fundamental parts in the demonstration of hyperbolicity, once a combinatorial property, called *richness* of the symbolic collision sequence of a trajectory, had been suitably defined:

(1) The "geometric-algebraic considerations" on the codimension of the manifolds describing the nonsufficient trajectory segments with a combinatorially rich symbolic collision structure;
(2) Proof of the fact that the set of phase points with a combinatorially non-rich collision sequence has measure zero.

In contrast to earlier proofs of ergodic properties of hard ball systems our method for obtaining full hyperbolicity will not be inductive; more precisely, we do not use the hyperbolicity of smaller systems. This has become possible since our Theorem 5.1 (settling step (2) here) which is a variant of the so-called *weak ball-avoiding theorems*, can now be proven without any inductive assumption. Nevertheless, the proof of our crucial Key Lemma 4.1, which copes with part (1), is inductive and, indeed, one of the main reasons to use *varying masses* is just that by choosing the mass of one particle to be equal to zero permits us to use an inductive assumption on the smaller system (but about its algebraic behaviour, only). As just mentioned, to prove our theorem we have introduced varying masses, and can only claim hyperbolicity for typical $(N+1)$-tuples $(m_1, \ldots, m_N; L) \in \mathbb{R}_+^{N+1}$ of outer geometric parameters. Furthermore, we will, in Section 3, *complexify the dynamics*. The complexification is required, on one hand, by the fact that our arguments in Section 4 use some algebraic tools that assume the ground field to be algebraically closed. Another advantage of the complexification is that, in the inductive derivation of Key Lemma 4.1, one does not have to worry about the sheer existence of orbit segments with a prescribed symbolic collision sequence: they do exist, thanks to the algebraic closedness of the complex field.

The paper is organized as follows: Section 2 is devoted to prerequisites. Section 3 then describes the complexified dynamics, while in Section 4 we establish Key Lemma 4.1 — both for the complexified and the real dynamics — settling, in particular, part (1) of the strategy. Section 5 provides the demonstration of our Main Theorem through the aforementioned Theorem 5.1, and, finally, Section 6 contains some comments and remarks.

## 2. Prerequisites
## Semi-dispersing billiards and hard ball systems

Our approach is based on a simple observation:

LEMMA 2.1. *The standard billiard ball flow with mass vector $\vec{m} \in \mathbb{R}_+^N$ is isomorphic to a semi-dispersing billiard.*

Because of the importance of the statement we sketch the proof.

*Proof.* Introduce new coordinates in the phase space as follows: For $(Q, V) \in \mathbf{M}$ let $\hat{q}_i = \sqrt{m_i} q_i$, $\hat{v}_i = \sqrt{m_i} v_i$, where thus $\hat{q}_i \in \sqrt{m_i}\, \mathbb{T}_L^\nu$ and $\sum_{i=1}^N \hat{v}_i^2 = 1$. Moreover, the vector $\hat{V} = (\hat{v}_1, \ldots, \hat{v}_N)$ necessarily belongs to the hyperplane $\sum_{i=1}^N \sqrt{m_i}\, \hat{v}_i = 0$. Denote $\mathbb{T}_{\vec{m}, L} = \prod_{i=1}^N (\sqrt{m_i}\, \mathbb{T}_L^\nu)$ and

$$(2.2) \quad \widehat{C}_{i,j} = \left\{ \hat{Q} = (\hat{q}_1, \ldots, \hat{q}_N) \in \mathbb{T}_{\vec{m}, L} : \left\| \frac{\hat{q}_i}{\sqrt{m_i}} - \frac{\hat{q}_j}{\sqrt{m_j}} \right\| < 2r \right\}$$

and $\widetilde{\hat{\mathbf{Q}}} = \mathbb{T}_{\vec{m},L} \setminus \cup_{1 \le i < j \le N} \hat{C}_{i,j}$. The equivalence relation $\hat{\Psi}_{\vec{m}}$ over $\mathbb{T}_{\vec{m},L}$ is as follows: $\hat{Q} \sim_{\hat{\Psi}_{\vec{m}}} \hat{Q}^*$ if and only if there exists an $a \in \mathbb{T}_L^\nu$ such that for every $i \in [1, N]$, $\hat{q}_i^* = \hat{q}_i + \sqrt{m_i}a$ allows us to define $\hat{\mathbf{Q}} = \widetilde{\hat{\mathbf{Q}}}/\hat{\Psi}_{\vec{m}}$. Let $\hat{\mathbf{M}}_{\vec{m},L} = \hat{\mathbf{Q}} \times \mathbb{S}^{d-1}$ be the unit tangent bundle of $\hat{\mathbf{Q}}$, and denote by $d\hat{\mu}$ the probability measure $const \cdot d\hat{Q} \cdot d\hat{V}$, where $d\hat{V}$ is the surface measure on the $(d-1)$−sphere $\mathbb{S}^{d-1}$, and $d\hat{Q}$ is the Lebesgue-measure on $\mathbb{T}_{\vec{m},L}$. Now the standard billiard ball system $(\mathbf{M}, S^\mathbb{R}, \mu)_{\vec{m},L}$ with mass vector $\vec{m}$ and the billiard system $(\hat{\mathbf{M}}_{\vec{m},L}, \hat{S}^\mathbb{R}, \hat{\mu})$ are isomorphic.

Indeed, we can reduce the question to the case of one-dimensional particles because for both models the velocity components perpendicular to the normal of impact remain unchanged. The claimed isomorphy for one-dimensional particles, however, is well-known, and for its simple proof we can refer, for instance, to Section 4 of Chapter 5 in [C-F-S(1981)]. Thus, the point is that in the isomorphic flow $\left(\hat{\mathbf{M}}_{\vec{m},L}, \hat{S}^\mathbb{R}, \hat{\mu}\right)$, the velocity transformations at collisions become orthogonal reflections across the tangent hyperplane of the boundary of $\hat{\mathbf{Q}}$; see also (3.9) for the mentioned velocity transformation.

The fact that the billiard system $(\hat{\mathbf{M}}_{\vec{m},L}, \hat{S}^\mathbb{R}, \hat{\mu})$ is semi-dispersing, is obvious, because the scattering bodies in $\hat{\mathbf{Q}}$ are cylinders built on ellipsoid bases (therefore this system is a cylindric billiard as introduced in [Sz(1993)]). □

For convenience and brevity, we will throughout use the concepts and notation, related to semi-dispersing billiards and hard ball systems, of the papers [K-S-Sz(1990)] and [Sim(1992)-I-II], respectively, and will only point out where and how different masses play a role.

*Remark* 2.3. By slightly generalizing the notion of semi-dispersing billiards with allowing Riemannian metrics different from the usual ones (see (2.2) of [K-S-Sz(1990)]), we could have immediately identified the standard billiard ball system $(\mathbf{M}, S^\mathbb{R}, \mu)_{\vec{m},L}$ with mass vector $\vec{m}$ with a (generalized) semi-dispersing billiard. (In the coordinates of the proof of Lemma 2.1, the Riemannian metric is

$$(d\rho)^2 = \sum_{i=1}^N \left(\|d\hat{q}_i\|^2 + \|ds\|^2\right)$$

where $\|ds\|^2$ is the square of the natural Riemannian metric on the unit sphere $\mathbb{S}^{d-1}$.) Then the standard billiard ball system $(\mathbf{M}, S^\mathbb{R}, \mu)_{\vec{m},L}$ itself will be a (generalized) semi-dispersing billiard. The advantage is that — as is easy to see — the results of [Ch-S(1987)] and [K-S-Sz(1990)] remain valid for this class, too. For simplifying our exposition therefore, we will omit the change of coordinates of Lemma 2.1 and will be using the notion of semi-dispersing billiards in this slightly more general sense in which the model of the introduction is a semi-dispersing billiard.

An often used abbreviation is the shorthand $S^{[a,b]}x$ for the trajectory segment $\{S^t x : a \leq t \leq b\}$. The natural projections from $\mathbf{M}$ onto its factor spaces are denoted, as usual, by $\pi : \mathbf{M} \to \mathbf{Q}$ and $p : \mathbf{M} \to \mathcal{S}^{N \cdot \nu - \nu - 1}$ or, sometimes, we simply write $\pi(x) = Q(x) = Q$ and $p(x) = V(x) = V$ for $x = (Q,V) \in \mathbf{M}$. Any $t \in [a,b]$ with $S^t x \in \partial M$ is called *a collision moment or collision time*.

As pointed out in previous works on billiards, the dynamics can only be defined for trajectories where the moments of collisions do not accumulate in any finite time interval (cf. Condition 2.1 of [K-S-Sz(1990)]). An important consequence of Theorem 5.3 of [V(1979)] is that — for semi-dispersing billiards under the nondegeneracy condition mentioned before our Main Theorem — there are *no trajectories at all with a finite accumulation point of collision moments* (see also [G(1981)] and [B-F-K(1998)]).

As a result, for an arbitrary nonsingular orbit segment $S^{[a,b]}x$ of the standard billiard ball flow, there is a uniquely defined maximal sequence $a \leq t_1 < t_2 < \cdots < t_n \leq b : n \geq 0$ of collision times and a uniquely defined sequence $\sigma_1 < \sigma_2 < \cdots < \sigma_n$ of "colliding pairs"; i.e. $\sigma_k = \{i_k, j_k\}$ whenever $Q(t_k) = \pi(S^{t_k}x) \in \partial \tilde{C}_{i_k, j_k}$. The sequence $\Sigma := \Sigma(S^{[a,b]}x) := (\sigma_1, \sigma_2, \ldots, \sigma_n)$ is called the *symbolic collision sequence* of the trajectory segment $S^{[a,b]}x$.

*Definition* 2.4. We say that the *symbolic collision sequence* $\Sigma = (\sigma_1, \ldots, \sigma_n)$ is *connected* if the collision graph of this sequence:

$$\mathcal{G}_\Sigma := (\mathcal{V} = \{1, 2, \ldots, N\}, \mathcal{E}_\Sigma := \{\{i_k, j_k\} : \text{ where } \sigma_k = \{i_k, j_k\}, \ 1 \leq k \leq n\})$$

is connected.

*Definition* 2.5. We say that the *symbolic collision sequence* $\Sigma = (\sigma_1, \ldots, \sigma_n)$ is *C-rich*, $C$ being a natural number, if it can be decomposed into at least $C$ consecutive, disjoint collision sequences in such a way that each of them is connected.

### Neutral subspaces, advance and sufficiency

Consider a *nonsingular* trajectory segment $S^{[a,b]}x$. Suppose that $a$ and $b$ are *not moments of collision*. Before defining the neutral linear space of this trajectory segment, we note that the tangent space of the configuration space $\mathbf{Q}$ at interior points can be identified with the common linear space

$$(2.6) \qquad \mathcal{Z} = \left\{ (w_1, w_2, \ldots, w_N) \in (\mathbb{R}^\nu)^N : \sum_{i=1}^N m_i w_i = 0 \right\}.$$

*Definition* 2.7. The *neutral space* $\mathcal{N}_0(S^{[a,b]}x)$ *of the trajectory segment* $S^{[a,b]}x$ *at time zero* $(a < 0 < b)$ is defined by the following formula:

$$\mathcal{N}_0(S^{[a,b]}x) = \{W \in \mathcal{Z} : \exists (\delta > 0) \text{ such that } \forall \alpha \in (-\delta, \delta),$$
$$p\left(S^a\left(Q(x) + \alpha W, V(x)\right)\right) = p(S^a x) \text{ and}$$
$$p\left(S^b\left(Q(x) + \alpha W, V(x)\right)\right) = p(S^b x)\}.$$

It is known (see (3) in Section 3 of [S-Ch (1987)]) that $\mathcal{N}_0(S^{[a,b]}x)$ is a linear subspace of $\mathcal{Z}$, and $V(x) \in \mathcal{N}_0(S^{[a,b]}x)$. The neutral space $\mathcal{N}_t(S^{[a,b]}x)$ of the segment $S^{[a,b]}x$ at time $t \in [a,b]$ is defined as follows:

$$(2.8) \qquad \mathcal{N}_t(S^{[a,b]}x) = \mathcal{N}_0\left(S^{[a-t,b-t]}(S^t x)\right).$$

It is clear that the neutral space $\mathcal{N}_t(S^{[a,b]}x)$ can be canonically identified with $\mathcal{N}_0(S^{[a,b]}x)$ by the usual identification of the tangent spaces of $\mathbf{Q}$ along the trajectory $S^{(-\infty,\infty)}x$ (see, for instance, Section 2 of [K-S-Sz(1990)] ).

Our next definition is that of the *advance*. Consider a nonsingular orbit segment $S^{[a,b]}x$ with symbolic collision sequence $\Sigma = (\sigma_1, \ldots, \sigma_n)$ $(n \geq 1)$ as at the beginning of the present section. For $x = (Q, V) \in \mathbf{M}$ and $W \in \mathcal{Z}$, $\|W\|$ sufficiently small, denote $T_W(Q, V) := (Q + W, V)$.

*Definition* 2.9. For any $1 \leq k \leq n$ and $t \in [a,b]$, the advance

$$\alpha(\sigma_k) : \; \mathcal{N}_t(S^{[a,b]}x) \to \mathbb{R}$$

is the *unique linear extension* of the linear functional defined in a sufficiently small neighbourhood of the origin of $\mathcal{N}_t(S^{[a,b]}x)$ in the following way:

$$\alpha(\sigma_k)(W) := t_k(x) - t_k(S^{-t} T_W S^t x).$$

It is now time to bring up the basic notion of *sufficiency* of a trajectory (segment). This is the most important necessary condition for the proof of the fundamental theorem for semi-dispersing billiards; see Condition (ii) of Theorem 3.6 and Definition 2.12 in [K-S-Sz(1990)] .

*Definition* 2.10.
(1) The *nonsingular trajectory segment* $S^{[a,b]}x$ (a and b are not supposed to be moments of collision) is said to be *sufficient* if and only if the dimension of $\mathcal{N}_t(S^{[a,b]}x)$ ($t \in [a,b]$) is minimal, i.e. $\dim \mathcal{N}_t(S^{[a,b]}x) = 1$.
(2) The *trajectory segment* $S^{[a,b]}x$ *containing exactly one singularity* is said to be *sufficient* if and only if both branches of this trajectory segment are sufficient.

For the notion of trajectory branches see, for example, the end of Section 2 in [Sim(1992)-I].

*Definition* 2.11. The *phase point* $x \in \mathbf{M}$ with *at most one singularity* is said to be sufficient if and only if its whole trajectory $S^{(-\infty,\infty)}x$ is sufficient, which means, by definition, that some of its bounded segments $S^{[a,b]}x$ are sufficient.

In the case of an orbit $S^{(-\infty,\infty)}x$ with exactly one singularity, sufficiency requires that both branches of $S^{(-\infty,\infty)}x$ be sufficient.

### Connecting Path Formula for particles with different masses

The Connecting Path Formula, abbreviated as CPF, was discovered for particles with identical masses in [Sim(1992)-II] . Its goal was to give an explicit description (by introducing a useful system of linear coordinates) of the neutral linear space $\mathcal{N}_0(S^{[-T,0]}x_0)$ in the language of the "advances" of the occurring collisions by using, as coefficients, linear expressions of the (pre-collision and post-collision) velocity differences of the colliding particles. Since it relied upon the conservation of the momentum, it has been natural to expect that the CPF can be generalized for particles with different masses as well. The case is, indeed, this, and next we give this generalization for particles with different masses. Since its structure is the same as that of the CPF for identical masses, our exposition follows closely the structure of [Sim(1992)-II] .

Consider a phase point $x_0 \in \mathbf{M}$ whose trajectory segment $S^{[-T,0]}x_0$ is not singular, $T > 0$. In the forthcoming discussion the phase point $x_0$ and the positive number $T$ will be fixed. All the velocities, $v_i(t) \in \mathbb{R}^\nu$, $i \in \{1, 2, \ldots, N\}$, $-T \leq t \leq 0$, appearing in the considerations are velocities of certain balls at specified moments $t$ and always with the starting phase point $x_0$ ($v_i(t)$ is the velocity of the $i^{\text{th}}$ ball at time $t$). We suppose that the moments $0$ and $-T$ are not moments of collision. We label the balls by the natural numbers $1, 2, \ldots, N$ (so the set $\{1, 2, \ldots, N\}$ is always the vertex set of the collision graph) and we denote by $e_1, e_2, \ldots, e_n$ the collisions of the trajectory segment $S^{[-T,0]}x_0$ (i.e. the edges of the collision graph) so that the time order of these collisions is just the opposite of the order given by the indices. More definitions and notation:

1. $t_i = t(e_i)$ denotes the time of the collision $e_i$, so that $0 > t_1 > t_2 > \cdots > t_n > -T$.

2. If $t \in \mathbb{R}$ is not a moment of collision ($-T \leq t \leq 0$), then
$$\Delta q_i(t) : \mathcal{N}_0(S^{[-T,0]}x_0) \to \mathbb{R}^\nu$$
is a linear mapping assigning to every element $W \in \mathcal{N}_0(S^{[-T,0]}x_0)$ the displacement of the $i^{\text{th}}$ ball at time $t$, provided that the configuration displacement at time zero is given by $W$. Originally, this linear mapping is only defined for vectors $W \in \mathcal{N}_0(S^{[-T,0]}x_0)$ close enough to the origin, but it can be uniquely extended to the whole space $\mathcal{N}_0(S^{[-T,0]}x_0)$ by preserving linearity.

3. $\alpha(e_i)$ denotes the advance of the collision $e_i$; thus

$$\alpha(e_i): \mathcal{N}_0(S^{[-T,0]}x_0) \to \mathbb{R}$$

is a linear mapping ($i = 1, 2, \ldots, n$).

4. The integers $1 = k(1) < k(2) < \cdots < k(l_0) \leq n$ are defined by the requirement that for every $j$ ($1 \leq j \leq l_0$) the graph $\{e_1, e_2, \ldots, e_{k(j)}\}$ consists of $N - j$ connected components (on the vertex set $\{1, 2, \ldots, N\}$, as always) while the graph $\{e_1, e_2, \ldots, e_{k(j)-1}\}$ consists of $N - j + 1$ connected components and, moreover, we require that the number of connected components of the whole graph $\{e_1, e_2, \ldots, e_n\}$ be equal to $N - l_0$. It is clear from this definition that the graph

$$\mathcal{T} = \{e_{k(1)}, e_{k(2)}, \ldots, e_{k(l_0)}\}$$

does not contain any loop, especially $l_0 \leq N - 1$.

Here we make two remarks commenting on the above notions.

*Remark* 2.12. We often do not indicate the variable $W \in \mathcal{N}_0(S^{[-T,0]}x_0)$ of the linear mappings $\Delta q_i(t)$ and $\alpha(e_i)$, for we will not be dealing with specific neutral tangent vectors $W$ but, instead, we think of $W$ as a typical (running) element of $\mathcal{N}_0(S^{[-T,0]}x_0)$ and $\Delta q_i(t)$, $\alpha(e_i)$ as linear mappings defined on $\mathcal{N}_0(S^{[-T,0]}x_0)$ in order to obtain an appropriate description of the neutral space $\mathcal{N}_0(S^{[-T,0]}x_0)$.

*Remark* 2.13. If $W \in \mathcal{N}_0(S^{[-T,0]}x_0)$ has the property $\Delta q_i(0)[W] = \lambda v_i(0)$ for some $\lambda \in \mathbb{R}$ and for all $i \in \{1, 2, \ldots, N\}$ (here $v_i(0)$ is the velocity of the $i^{\text{th}}$ ball at time zero), then $\alpha(e_k)[W] = \lambda$ for all $k = 1, 2, \ldots, n$. This particular $W$ corresponds to the direction of the flow. In the sequel we shall often refer to this remark.

Let us fix two *distinct* balls $\alpha, \omega \in \{1, 2, \ldots, N\}$ that *are in the same connected component* of the collision graph $\mathcal{G}_n = \{e_1, e_2, \ldots, e_n\}$. The CPF expresses the relative displacement $\Delta q_\alpha(0) - \Delta q_\omega(0)$ in terms of the advances $\alpha(e_i)$ and the relative velocities occurring at these collisions $e_i$. In order to be able to formulate the CPF we need to define some graph-theoretic notions concerning the pair of vertices $(\alpha, \omega)$.

*Definition* 2.14. Since the graph $\mathcal{T} = \{e_{k(1)}, e_{k(2)}, \ldots, e_{k(l_0)}\}$ contains no loop and the vertices $\alpha, \omega$ belong to the same connected component of $\mathcal{T}$, there is a unique path $\Pi(\alpha, \omega) = \{f_1, f_2, \ldots, f_h\}$ in the graph $\mathcal{T}$ connecting the vertices $\alpha$ and $\omega$. The edges $f_i \in \mathcal{T}$ ($i = 1, 2, \ldots, h$) are listed successively along this path $\Pi(\alpha, \omega)$ starting from $\alpha$ and ending at $\omega$. The vertices of the path $\Pi(\alpha, \omega)$ are denoted by $\alpha = B_0, B_1, B_2, \ldots, B_h = \omega$ indexed along this path going from $\alpha$ to $\omega$, so the edge $f_i$ connects the vertices $B_{i-1}$ and $B_i$ ($i = 1, 2, \ldots, h$).

When trying to compute $\Delta q_\alpha(0) - \Delta q_\omega(0)$ by using the advances $\alpha(e_i)$ and the relative velocities at these collisions, it turns out that not only the collisions $f_i$ ($i = 1, 2, \ldots, h$) make an impact on $\Delta q_\alpha(0) - \Delta q_\omega(0)$, but some other adjacent edges too. This motivates the following definition:

*Definition* 2.15. Let $i \in \{1, 2, \ldots, h-1\}$ be an integer. We define the *set of $\mathcal{A}_i$ adjacent edges* at the vertex $B_i$ as follows:
$$\mathcal{A}_i = \{e_j : j \in \{1, 2, \ldots, n\} \text{ and } (t(e_j) - t(f_i)) \cdot (t(e_j) - t(f_{i+1})) < 0 \text{ and}$$
$$B_i \text{ is a vertex of } e_j\}.$$

We adopt a similar definition of the sets $\mathcal{A}_0, \mathcal{A}_h$ of adjacent edges at the vertices $B_0$ and $B_h$, respectively:

*Definition* 2.16.
$$\mathcal{A}_0 = \{e_j : 1 \leq j \leq n \text{ and } t(e_j) > t(f_1) \text{ and } B_0 \text{ is a vertex of } e_j\};$$
$$\mathcal{A}_h = \{e_j : 1 \leq j \leq n \text{ and } t(e_j) > t(f_h) \text{ and } B_h \text{ is a vertex of } e_j\}.$$

We note that the sets $\mathcal{A}_0, \mathcal{A}_1, \ldots, \mathcal{A}_h$ are not necessarily mutually disjoint.

Finally, we need to define the "contribution" of the collision $e_j$ to $\Delta q_\alpha(0) - \Delta q_\omega(0)$ which is composed from the relative velocities just before and after the moment $t(e_j)$ of the collision $e_j$.

*Definition* 2.17. For $i \in \{1, 2, \ldots, h\}$ the *contribution* $\Gamma(f_i)$ of the edge $f_i \in \Pi(\alpha, \omega)$ is given by the formula

$$\Gamma(f_i) = \begin{cases} v^-_{B_{i-1}}(t(f_i)) - v^-_{B_i}(t(f_i)), & \text{if } t(f_{i-1}) < t(f_i) \text{ and } t(f_{i+1}) < t(f_i); \\[1ex] v^+_{B_{i-1}}(t(f_i)) - v^+_{B_i}(t(f_i)), & \text{if } t(f_{i-1}) > t(f_i) \text{ and } t(f_{i+1}) > t(f_i); \\[1ex] \dfrac{1}{m_{B_{i-1}} + m_{B_i}} \Big[m_{B_{i-1}}\big(v^-_{B_{i-1}}(t(f_i)) - v^-_{B_i}(t(f_i))\big) + m_{B_i}\big(v^+_{B_{i-1}}(t(f_i)) \\ \qquad - v^+_{B_i}(t(f_i))\big)\Big], & \text{if } t(f_{i+1}) < t(f_i) < t(f_{i-1}) \\[1ex] \dfrac{1}{m_{B_{i-1}} + m_{B_i}} \Big[m_{B_{i-1}}\big(v^+_{B_{i-1}}(t(f_i)) - v^+_{B_i}(t(f_i))\big) + m_{B_i}\big(v^-_{B_{i-1}}(t(f_i)) \\ \qquad - v^-_{B_i}(t(f_i))\big)\Big], & \text{if } t(f_{i-1}) < t(f_i) < t(f_{i+1}). \end{cases}$$

Here $v^-_{B_i}(t(f_i))$ denotes the velocity of the $B_i^{\text{th}}$ particle just before the collision $f_i$ (occurring at time $t(f_i)$) and, similarly, $v^+_{B_i}(t(f_i))$ is the velocity of the same particle just after the mentioned collision. We also note that, by convention, $t(f_0) = 0 > t(f_1)$ and $t(f_{h+1}) = 0 > t(f_h)$. Apparently, the time order plays an important role in this definition.

*Definition* 2.18. For $i \in \{0, 1, 2, \ldots, h\}$ the *contribution* $\Gamma_i(e_j)$ of an edge $e_j \in \mathcal{A}_i$ is defined as follows:

$$\Gamma_i(e_j) = \text{sign}\big(t(f_i) - t(f_{i+1})\big) \frac{m_C}{m_{B_i} + m_C}$$
$$\cdot \big[\big(v^+_{B_i}(t(e_j)) - v^+_C(t(e_j))\big) - \big(v^-_{B_i}(t(e_j)) - v^-_C(t(e_j))\big)\big]$$

where $C$ is the vertex of $e_j$ different from $B_i$.

Here again we adopt the convention $t(f_0) = 0 > t(e_j)$ ($e_j \in \mathcal{A}_0$) and $t(f_{h+1}) = 0 > t(e_j)$ ($e_j \in \mathcal{A}_h$). We note that, by the definition of the set $\mathcal{A}_i$, exactly one of the two possibilities $t(f_{i+1}) < t(e_j) < t(f_i)$ and $t(f_i) < t(e_j) < t(f_{i+1})$ occurs. The subscript $i$ of $\Gamma$ is only needed because an edge $e_j \in \mathcal{A}_{i_1} \cap \mathcal{A}_{i_2}$ ($i_1 < i_2$) has two contributions at the vertices $B_{i_1}$ and $B_{i_2}$ which are just the endpoints of $e_j$.

We are now in the position of formulating the Connecting Path Formula:

PROPOSITION 2.19. *With the definitions and notation above, the following sum is an expression for $\Delta q_\alpha(0) - \Delta q_\omega(0)$ in terms of the advances and relative velocities of collisions*:

$$\Delta q_\alpha(0) - \Delta q_\omega(0) = \sum_{i=1}^{h} \alpha(f_i)\Gamma(f_i) + \sum_{i=0}^{h} \sum_{e_j \in \mathcal{A}_i} \alpha(e_j)\Gamma_i(e_j).$$

The proof of the proposition follows the proof of Simányi's CPF (Lemma 2.9 [Sim(1992)-II] ) with the only difference that Lemma 2.8 of [Sim(1992)-II] is replaced here by the following:

LEMMA 2.20. *If $e$ is a collision at time $t$ between the particles $B$ and $C$, then*

$$v^+_B(t) - v^+_C(t) = \frac{m_C}{m_B + m_C}\big(v^+_B(t) - v^+_C(t)\big) + \frac{m_B}{m_B + m_C}\big(v^-_B(t) - v^-_C(t)\big)$$

*and*

$$v^+_B(t) - v^-_B(t) = \frac{m_C}{m_B + m_C}\big[\big(v^+_B(t) - v^+_C(t)\big) - \big(v^-_B(t) - v^-_C(t)\big)\big].$$

*Proof.* The lemma is an easy consequence of the conservation of momentum for the collision $e$. □

*Remark* 2.21. In Section 3, we will complexify the system, and also allow complex masses, in particular. It is easy to see that the CPF of Proposition 2.19 and the whole discussion of this subsection will still hold for complex masses if we assume, in addition, that for every $e = \{B, C\}$, occurring in the symbolic collision sequence of $S^{[-T,0]}x_0$, $m_B + m_C \neq 0$.

## 3. The complexified billiard map
## Partial linearity of the dynamics

The aim of this section is to understand properly the algebraic relationship between the kinetic data of the billiard flow measured at different moments of time.

From now on we shall investigate orbit segments $S^{[0,T]}x_0$ ($T > 0$) of the standard billiard ball flow $\left(\tilde{\mathbf{M}}, \{S^t\}, \mu, \vec{m}, L\right)$. We note here that $\tilde{\mathbf{M}} = \tilde{\mathbf{Q}} \times \mathcal{E}$, where the configuration space $\tilde{\mathbf{Q}}$ is as defined right after (1.1) and $\mathcal{E}$ is the velocity sphere

$$\mathcal{E} = \left\{(v_1, \ldots, v_N) \in \mathbb{R}^{\nu N} \,\bigg|\, \sum_{i=1}^{N} m_i v_i = 0 \text{ and } \sum_{i=1}^{N} m_i \|v_i\|^2 = 1 \right\}$$

introduced in Section 1. Also note that in the geometric-algebraic considerations of the upcoming sections we do not use the equivalence relation $\Psi$ of the introduction. This is why we work in $\tilde{M}$ rather than in $M$. Later on even the conditions $\sum_{i=1}^{N} m_i v_i = 0$ and $\sum_{i=1}^{N} m_i \|v_i\|^2 = 1$ will be dropped (cf. Remark 3.14).

The symbolic collision sequence of $S^{[0,T]}x_0$ is denoted by

$$\Sigma\left(S^{[0,T]}x_0\right) = (\sigma_1, \ldots, \sigma_n).$$

The symbol $v_i^k = v_i^+(t_k) \in \mathbb{R}^\nu$ denotes the velocity $\dot{q}_i$ of the $i^{\text{th}}$ ball right *after* the $k^{\text{th}}$ collision $\sigma_k = \{i_k, j_k\}$ ($1 \leq i_k < j_k \leq N$) occurring at time $t_k$, $k = 1, \ldots, n$. Of course, there is no need to deal with the velocities right *before* collisions, since $v_i^-(t_k) = v_i^+(t_{k-1}) = v_i^{k-1}$. As usual, $q_i^k \in \mathbb{T}_L^\nu = \mathbb{R}^\nu / L \cdot \mathbb{Z}^\nu$ denotes the position of the center of the $i^{\text{th}}$ ball at the moment $t_k$ of the $k^{\text{th}}$ collision.

Let us now fix the symbolic collision sequence $\Sigma = (\sigma_1, \ldots, \sigma_n)$ ($n \geq 1$), and explore the algebraic relationship between the data

$$\left\{q_i^{k-1}, v_i^{k-1} \,|\, i = 1, \ldots, N\right\}$$

and

$$\left\{q_i^k, v_i^k \,|\, i = 1, \ldots, N\right\},$$

$k = 1, \ldots, n$. By definition, the data $\{q_i^0, v_i^0 | i = 1, \ldots, N\}$ correspond to the initial (noncollision) phase point $x_0$. We also set $t_0 = 0$.

Since we would like to carry out arithmetic operations on these data, the periodic positions $q_i^k \in \mathbb{T}_L^\nu$ are not suitable for this purpose. Therefore, instead of studying the genuine orbit segments $S^{[0,T]}x_0$, we will deal with their Euclidean liftings.

PROPOSITION 3.1. *Let $S^{[0,T]}x_0 = \{q_i(t), v_i(t) | 0 \le t \le T\}$ be an orbit segment as above, and assume that certain pre-images (Euclidean liftings) $\tilde{q}_i(0) = \tilde{q}_i^0 \in \mathbb{R}^\nu$ of the positions $q_i(0) \in \mathbb{T}_L^\nu = \mathbb{R}^\nu / L \cdot \mathbb{Z}^\nu$ are given. Then there is a uniquely defined, continuous, Euclidean lifting $\{\tilde{q}_i(t) \in \mathbb{R}^\nu | 0 \le t \le T\}$ of the given orbit segment that is an extension of the initial lifting $\{\tilde{q}_i^0 | i = 1, \ldots, N\}$. Moreover, for every collision $\sigma_k$ there exists a uniquely defined integer vector $a_k \in \mathbb{Z}^\nu$ – named the adjustment vector of $\sigma_k$ – such that*

$$(3.2) \qquad \left\| \tilde{q}_{i_k}^k - \tilde{q}_{j_k}^k - L \cdot a_k \right\|^2 - 4r^2 = 0.$$

*The orbit segment $\tilde{\omega} = \{\tilde{q}_i^k, v_i^k | i = 1, \ldots, N; k = 0, \ldots, n\}$ is called a lifted orbit segment with the system of adjustment vectors $\mathcal{A} = (a_1, \ldots, a_n) \in \mathbb{Z}^{n\nu}$.*

The proof of this proposition is straightforward and we omit it.

The next result establishes the already-mentioned polynomial relationships between the kinetic data

$$\left\{ \tilde{q}_i^{k-1}, v_i^{k-1} | i = 1, \ldots, N \right\}$$

and

$$\{\tilde{q}_i^k, v_i^k | i = 1, \ldots, N\},$$

$k = 1, \ldots, n$.

PROPOSITION 3.3. *Using all notation and notions from above, one has the following polynomial relations between the kinetic data at $\sigma_{k-1}$ and $\sigma_k$, $k = 1, \ldots, n$. In order to simplify the notation, these equations are written as if $\sigma_k$ were $\{1, 2\}$:*

$$(3.4) \qquad v_i^k = v_i^{k-1}, \quad i \notin \{1, 2\},$$

$$(3.5) \qquad m_1 v_1^k + m_2 v_2^k = m_1 v_1^{k-1} + m_2 v_2^{k-1}$$

*(conservation of the momentum),*
(3.6)
$$v_1^k - v_2^k = v_1^{k-1} - v_2^{k-1} - \frac{1}{2r^2} \left\langle v_1^{k-1} - v_2^{k-1}; \tilde{q}_1^k - \tilde{q}_2^k - L \cdot a_k \right\rangle \cdot \left( \tilde{q}_1^k - \tilde{q}_2^k - L \cdot a_k \right)$$

(*reflection of the relative velocity determined by the elastic collision*),

$$(3.7) \qquad \tilde{q}_i^k = \tilde{q}_i^{k-1} + \tau_k v_i^{k-1}, \quad i = 1, \ldots, N,$$

where the time slot $\tau_k = t_k - t_{k-1}$ ($t_0 := 0$) is determined by the quadratic equation

$$(3.8) \qquad \left\| \tilde{q}_1^{k-1} - \tilde{q}_2^{k-1} + \tau_k \left( v_1^{k-1} - v_2^{k-1} \right) - L \cdot a_k \right\|^2 = 4r^2.$$

The proof of this proposition is also obvious.

*Remark.* Observe that the new velocities $v_1^k$ and $v_2^k$ can be computed from (3.5)–(3.6) as follows:

$$(3.9) \quad v_1^k = v_1^{k-1} - \frac{m_2}{2r^2(m_1 + m_2)} \left\langle v_1^{k-1} - v_2^{k-1}; \tilde{q}_1^k - \tilde{q}_2^k - L \cdot a_k \right\rangle$$
$$\cdot \left( \tilde{q}_1^k - \tilde{q}_2^k - L \cdot a_k \right),$$
$$v_2^k = v_2^{k-1} + \frac{m_1}{2r^2(m_1 + m_2)} \left\langle v_1^{k-1} - v_2^{k-1}; \tilde{q}_1^k - \tilde{q}_2^k - L \cdot a_k \right\rangle$$
$$\cdot \left( \tilde{q}_1^k - \tilde{q}_2^k - L \cdot a_k \right).$$

We also note that the equations above even extend analytically to the case when one of the two masses, say $m_1$, is equal to zero:

$$(3.10)$$
$$v_1^k = v_1^{k-1} - \frac{1}{2r^2} \left\langle v_1^{k-1} - v_2^{k-1}; \tilde{q}_1^k - \tilde{q}_2^k - L \cdot a_k \right\rangle \cdot \left( \tilde{q}_1^k - \tilde{q}_2^k - L \cdot a_k \right),$$
$$v_2^k = v_2^{k-1}.$$

We call the attention of the reader to the fact that in our understanding the symbol $\langle .; . \rangle$ denotes the Euclidean inner product of $\nu$-dimensional real vectors and $\| . \|$ is the corresponding norm.

### The complexification of the billiard map

Given the pair $(\Sigma, \mathcal{A}) = (\sigma_1, \ldots, \sigma_n; a_1, \ldots, a_n)$, the equations (3.4)–(3.8) make it possible to iteratively compute the kinetic data $\left\{ \tilde{q}_i^k, v_i^k \,|\, i = 1, \ldots, N \right\}$ by using the preceding data $\left\{ \tilde{q}_i^{k-1}, v_i^{k-1} \,|\, i = 1, \ldots, N \right\}$. Throughout these computations one only uses the field operations and square roots. (The latter one is merely used when computing $\tau_k$ as the root of the quadratic equation (3.8). Obviously, since – at the moment – we are dealing with genuine, real orbit segments, in this dynamically realistic situation the equations (3.8) have two distinct, positive real roots, and $\tau_k$ is the smaller one.) Therefore, the data $\left\{ \tilde{q}_i^k, v_i^k \,|\, i = 1, \ldots, N \right\}$ can be expressed as certain algebraic functions of the initial data $\{\tilde{q}_i^0,\, v_i^0|\, i = 1, \ldots, N\}$, and the arising algebraic functions merely contain field operations, square roots, the radius $r$ of the balls, the size $L > 0$ of the torus and, finally, the masses $m_i$ as constants.

Since these algebraic functions make full sense over the complex field $\mathbb{C}$ and, after all, our proof of the theorem requires the complexification, we are now going to complexify the whole system by considering the kinetic variables, the size $L$, and the masses as complex ones and by retaining the polynomial equations (3.4)–(3.8). However, due to the ambiguity of selecting a root of (3.8) out of the two, it proves to be important to explore first the algebraic frame of the relations (3.4)-(3.8) connecting the studied variables. This is what we do now.

### The field extension associated with the pair $(\Sigma, \mathcal{A})$

To avoid misunderstanding we immediately stress that the field extensions $\mathbb{K} = \mathbb{K}(\Sigma; \mathcal{A})$ to be defined below will also depend on a sequence $\vec{\tau} = (\tau_0, \ldots, \tau_{n-1})$ of field elements to be introduced successively in Definition 3.11. If we also want to emphasize the dependence of $\mathbb{K}$ on $\vec{\tau}$, then we will write $\mathbb{K} = \mathbb{K}(\Sigma; \mathcal{A}) = (\Sigma, \mathcal{A}, \vec{\tau})$.

We are going to define the commutative function field $\mathbb{K} = \mathbb{K}(\Sigma; \mathcal{A})$ generated by all functions

$$\left\{ (\tilde{q}_i^k)_j,\, (v_i^k)_j|\, i = 1, \ldots, N;\, k = 0, \ldots, n;\, j = 1, \ldots, \nu \right\}$$

of the lifted orbit segments corresponding to the given parameters

$$(\Sigma, \mathcal{A}) = (\sigma_1, \ldots, \sigma_n; a_1, \ldots, a_n)$$

in such a way that the field $\mathbb{K}(\Sigma; \mathcal{A})$ incorporates all algebraic relations among these variables that are consequences of equations (3.4)–(3.8). (Here the subscript $j$ denotes the $j^{\text{th}}$ component of a $\nu$-vector.) In our setup the ground field of allowed constants (coefficients) is, by definition, the complex field $\mathbb{C}$. The precise definition of $\mathbb{K}_n = \mathbb{K}(\Sigma; \mathcal{A})$ is:

*Definition* 3.11. For $n = 0$ the field $\mathbb{K}_0 = \mathbb{K}(\emptyset; \emptyset)$ is the transcendental extension $\mathbb{C}(\mathcal{B})$ of the coefficient field $\mathbb{C}$ by the algebraically independent formal variables

$$\mathcal{B} = \left\{ (\tilde{q}_i^0)_j,\, (v_i^0)_j,\, m_i,\, L|\, i = 1, \ldots, N;\, j = 1, \ldots, \nu \right\}.$$

Suppose now that $n > 0$ and the commutative field $\mathbb{K}_{n-1} = \mathbb{K}(\Sigma'; \mathcal{A}')$ has already been defined, where $\Sigma' = (\sigma_1, \ldots, \sigma_{n-1})$, $\mathcal{A}' = (a_1, \ldots, a_{n-1})$. Then consider the quadratic equation $b_n \tau_n^2 + c_n \tau_n + d_n = 0$ in (3.8) with $k = n$ as a polynomial equation defining a new field element $\tau_n$ to be adjoined to the field $\mathbb{K}_{n-1} = \mathbb{K}(\Sigma'; \mathcal{A}')$. (Recall that the coefficients $b_n$, $c_n$, $d_n$ come from the field $\mathbb{K}_{n-1}$.)

There are two possibilities:

(i) The *quadratic polynomial* $b_n x^2 + c_n x + d_n$ *is reducible over the field* $\mathbb{K}_{n-1}$. Then $\mathbb{K}_n = \mathbb{K}_{n-1}$.

(ii) The *polynomial* $b_n x^2 + c_n x + d_n$ *is irreducible over the field* $\mathbb{K}_{n-1}$. Then we define $\mathbb{K}_n = \mathbb{K}(\Sigma; \mathcal{A})$ as the extension of $\mathbb{K}_{n-1} = \mathbb{K}(\Sigma'; \mathcal{A}')$ by the root $\tau_n$ of this irreducible polynomial.

*Remark* 3.12. The importance of the field $\mathbb{K}_n$ is underscored by the fact that this field encodes all algebraic relations among the kinetic data that follow from the polynomial equations (3.4)–(3.8). Furthermore, Proposition 3.3 gives an iterative computation rule for successively obtaining the kinetic variables

$$\left\{ \tilde{q}_i^k, \, v_i^k \,|\, i = 1, \ldots, N \right\},$$

$k = 0, \ldots, n$. The field $\mathbb{K}_n = \mathbb{K}(\Sigma, \mathcal{A})$ is, after all, the algebraic frame of such computations. We note that – everywhere in what follows – the symbols $\langle \,.\,;\,.\, \rangle$ and $\| \,.\, \|^2$ do not refer to an Hermitian inner product but, rather, $\langle x; y \rangle = \sum_{j=1}^{\nu} x_j y_j$ and $\|x\|^2 = \sum_{j=1}^{\nu} x_j^2$, so that these expressions retain their polynomial form.

*Remark* 3.13. Later on it will be necessary to express each kinetic variable $\left( \tilde{q}_i^k \right)_j$ and $\left( v_i^k \right)_j$ as an algebraic function of the initial variables of $\mathcal{B}$. This raises, however, an important question: Which one of the two roots $\tau_k$ of (3.8) should be considered during these computations? This is no problem if the polynomial $b_k x^2 + c_k x + d_k$, defining $\tau_k$, is irreducible over the field $\mathbb{K}_{k-1}$, because the two roots of this polynomial are algebraically equivalent over $\mathbb{K}_{k-1}$, and any of them can be chosen as $\tau_k$. However, if the polynomial $b_k x^2 + c_k x + d_k$ is reducible over the field $\mathbb{K}_{k-1}$, then it is necessary to make a decision and assign to one of the two roots (both in the field $\mathbb{K}_{k-1}$ now) the role of $\tau_k$. This is what we do. The result is the field $\mathbb{K}_n = \mathbb{K}(\Sigma; \mathcal{A})$ endowed with a distinguished n-tuple $(\tau_1, \ldots, \tau_n)$ of its elements. The field $\mathbb{K}_n$ with the distinguished $n$-tuple $\vec{\tau}$ is denoted by $\mathbb{K}_n = \mathbb{K}(\Sigma; \mathcal{A}; \vec{\tau})$. Then the algebraic object $\mathbb{K}(\Sigma; \mathcal{A}; \vec{\tau})$ completely controls the whole process of the iterative computation of the kinetic variables.

Note that in the genuine, real case the two roots of $b_k x^2 + c_k x + d_k$ are distinct, positive real numbers, and – by the nature of the billiard dynamics – the chosen value is always the smaller one.

Finally, we mention here a simple case when the polynomial $b_k x^2 + c_k x + d_k$ is reducible over the field $\mathbb{K}_{k-1}$. Namely, this takes place whenever $\sigma_{k-1} = \sigma_k$ and $a_{k-1} = a_k$. In that case the value $\tau_k \equiv 0$ is clearly a solution of the above polynomial $b_k x^2 + c_k x + d_k$, and, being so, that polynomial is reducible over the field $\mathbb{K}_{k-1}$; see also 3.31–3.33.

*Remark* 3.14. The reader may wonder why we did not postulate the algebraic dependencies $\sum_{i=1}^{N} m_i v_i^k = 0$, $\sum_{i=1}^{N} m_i \|v_i^k\|^2 = 1$. The answer is the following: The definition of the neutral linear space together with its characterization via the Connecting Path Formula (CPF, see Section 2 of this paper or Lemma 2.9 of [Sim(1992)-II]), i.e. the partial linearity of the billiard flow, are invariant under

(1) all uniform velocity translations (adding the same vector to all velocities) and

(2) time rescalings.

*Remark* 3.15. All fields occurring in this paper are only defined up to an isomorphism over the coefficient field $\mathbb{C}$. Therefore, the statements like "the field $\mathbb{K}_2$ is an extension of $\mathbb{K}_1$" should be understood as follows: $\mathbb{K}_1$ is a subfield of $\mathbb{K}_2$ after the natural identification of the generating variables bearing the same name. For the necessary notions, properties, and results from the theory of field extensions and Galois theory, the reader may look up, for instance, the book by I. Stewart, [St(1973)].

*Remark* 3.16. As said before (see Definition 3.11), the collection $\mathcal{B}$ of elements of the field $\mathbb{K}(\Sigma; \mathcal{A}; \vec{\tau})$ is a base of transcendence in that field over the subfield $\mathbb{C}$, and the degree of the extension $\mathbb{K}(\Sigma; \mathcal{A}; \vec{\tau}) : \mathbb{C}(\mathcal{B}) = \mathbb{K}_n : \mathbb{K}_0$ is a power of two. Moreover, each of the following sets is a generator for the field $\mathbb{K}_n$:

(a) $\mathbb{C} \cup \mathcal{B} \cup \{\tau_1, \ldots, \tau_n\}$;

(b) $\mathbb{C} \cup \{m_1, \ldots, m_N, L\} \cup \{(\tilde{q}_i^n)_j, (v_i^n)_j \mid i = 1, \ldots, N; j = 1, \ldots, \nu\} \cup \{\tau_1, \ldots, \tau_n\}$;

(c) $\mathbb{C} \cup \{m_1, \ldots, m_N, L\} \cup \{(\tilde{q}_i^k)_j, (v_i^k)_j \mid i = 1, \ldots, N; j = 1, \ldots, \nu; k = 0, \ldots, n\}$.

*Remark* 3.17. The above procedure of constructing field extensions is closely related to the classical theory of geometric constructions by a ruler and compass, and this is not surprising: the billiard trajectory is constructed by intersecting a straight line with a sphere and then reflecting it across the tangent hyperplane of the sphere. This is a sort of classical geometric construction in $\nu$ dimensions. For the details see, for instance, Sections 57 and 60 of [VDW (1955)], or Chapter 5 of [St(1973)].

Let us fix the pair $(\Sigma; \mathcal{A})$ and the $n$-tuple $\vec{\tau} = (\tau_1, \ldots, \tau_n)$ of elements of the field $\mathbb{K}(\Sigma, \mathcal{A})$. We are now defining a $(2\nu + 1)N + 1$-dimensional complex analytic manifold $\tilde{\Omega} = \tilde{\Omega}(\Sigma, \mathcal{A}, \vec{\tau})$ and certain holomorphic functions $\tilde{q}_i^k$, $v_i^k$ $: \tilde{\Omega} \to \mathbb{C}^\nu$, $m_i : \tilde{\Omega} \to \mathbb{C}$, $L : \tilde{\Omega} \to \mathbb{C}$, $\tau_k : \tilde{\Omega} \to \mathbb{C}$, $i = 1, \ldots, N$, $k = 0, \ldots, n$ ($\tau_0$ is not defined). But first we introduce:

*Definition* 3.18. Define the domain $D = D(\Sigma, \mathcal{A}, \vec{\tau}) \subset \mathbb{C}^{(2\nu+1)N+1}$ as the set of all complex $(2\nu+1)N+1$-tuples

$$\{(\tilde{q}_i^0)_j, (v_i^0)_j, m_i, L \,|\, i = 1, \ldots, N;\, j = 1, \ldots, \nu\}$$

for which

(a) The leading coefficient $\left\|v_{i_k}^{k-1} - v_{j_k}^{k-1}\right\|^2$ and the discriminant of the quadratic equation (3.8) is never zero, $k = 1, \ldots, n$. (The latter condition is equivalent to
$$\left\langle v_{i_k}^{k-1} - v_{j_k}^{k-1};\, \tilde{q}_{i_k}^k - \tilde{q}_{j_k}^k - L \cdot a_k \right\rangle \neq 0;$$
see also Remark 3.29.) And

(b) $m_{i_k} + m_{j_k} \neq 0$ for $k = 1, \ldots, n$, provided that the iterative computation of the kinetic variables is carried out,
$$\left\{(\tilde{q}_i^k)_j, (v_i^k)_j \,|\, i = 1, \ldots, N;\, j = 1, \ldots, \nu;\, k = 0, \ldots, n\right\}.$$
as governed by the polynomial equations (3.4)–(3.8) and by the fixed selection of elements $\vec{\tau}$.

*Remark* 1. When inverting the dynamics, we see that
$$\left\|v_{i_k}^{k-1} - v_{j_k}^{k-1}\right\|^2 = \left\|v_{i_k}^k - v_{j_k}^k\right\|^2,$$
and
$$\left\langle v_{i_k}^{k-1} - v_{j_k}^{k-1};\, \tilde{q}_{i_k}^k - \tilde{q}_{j_k}^k - L \cdot a_k \right\rangle = -\left\langle v_{i_k}^k - v_{j_k}^k;\, \tilde{q}_{i_k}^k - \tilde{q}_{j_k}^k - L \cdot a_k \right\rangle.$$

*Remark* 2. Note that we require the validity of (a) for any branch of the square root function, when computing the $\tau_k$'s with irreducible defining polynomials (3.8).

LEMMA 3.19. *The complement set $\mathbb{C}^{(2\nu+1)N+1} \setminus D$ is a proper, closed, algebraic subset of $\mathbb{C}^{(2\nu+1)N+1}$, especially a finite union of complex analytic submanifolds with codimension at least one, so that the open set $D \subset \mathbb{C}^{(2\nu+1)N+1}$ is connected and dense.*

*Proof.* We have the polynomial equations (3.4)–(3.8) making it possible to set up an algorithm for iteratively computing the kinetic variables measured at different times. The point is not just this iterative algorithm, but also its invertibility (time reversibility). The inverse process has similar algebraic properties, for it just means time reversal. One can easily prove by an induction on the length $n$ of $\Sigma$ that the iteratively defined dynamics determines a several-to-several mapping with maximum rank between a nonempty Zariski open set of the kinetic variables with superscript 1 and a nonempty Zariski open set of the kinetic variables measured right after the last reflection (with superscript $n$). (A Zariski open set is the complement of a closed algebraic set defined as the simultaneous zero set of finitely many polynomials.) Then one can go ahead with the induction from $n$ to $n+1$, because all possible obstructions to extending the process only occur on proper algebraic submanifolds of the Zariski open set of kinetic variables with superscript $n$. □

*Definition* 3.20. We define $\tilde{\Omega}(\Sigma, \mathcal{A}, \vec{\tau})$ as the set of all complex $(2\nu(n+1)N + N + n + 1)$-tuples

$$\omega = \left( \left(\tilde{q}_i^k\right)_j, (v_i^k)_j,\, m_i,\, L,\, \tau_l \,\middle|\, i=1,\ldots,N;\; k=0,\ldots,n;\; j=1,\ldots,\nu;\; l=1,\ldots,n \right)$$

for which these coordinates

(1) are interrelated by the equations (3.4)–(3.8),
(2) respect the choices of $\tau_k$ prescribed by $\vec{\tau}$ (cf. 3.11, and also 3.13),

and the vector

$$\vec{x}(\omega) = \left( \left(\tilde{q}_i^0\right)_j,\, (v_i^0)_j,\, m_i,\, L \,\middle|\, i=1,\ldots,N;\; j=1,\ldots,\nu \right)$$

of initial data belongs to the set $D(\Sigma, \mathcal{A}, \vec{\tau})$ defined above.

It is clear now that $\tilde{q}_i^k, v_i^k : \tilde{\Omega} \to \mathbb{C}^\nu$, $m_i : \tilde{\Omega} \to \mathbb{C}$, $L : \tilde{\Omega} \to \mathbb{C}$, and $\tau_k : \tilde{\Omega} \to \mathbb{C}$ ($i=1,\ldots,N$, $k=0,\ldots,n$) are holomorphic functions on $\tilde{\Omega} = \tilde{\Omega}(\Sigma, \mathcal{A}, \vec{\tau})$. The complex analytic manifold $\tilde{\Omega}$ endowed with the above holomorphic functions can justifiably be considered as the complexification of the $(\Sigma, \mathcal{A}, \vec{\tau})$-iterated billiard map.

*Remark.* The careful reader has certainly noted the following principal observation: the $\tau_k$'s figure in two different roles. They first denote successively chosen field elements of the extensions introduced in Definition 3.11, and secondly they also denote complex-valued functions on $\tilde{\Omega}$ (or multivalued complex functions on $D$; cf. Definitions 3.20 and 3.18).

*Remark.* It is very likely that the manifold $\tilde{\Omega}$ is connected. Nevertheless, in our proof we do not need this connectedness, and are not going to pursue the goal of proving it.

## The complex neutral space $\mathcal{N}(\omega)$

Fix a base point $\omega \in \tilde{\Omega}(\Sigma, \mathcal{A}, \vec{\tau})$. The tangent space $\mathcal{T}_\omega \tilde{\Omega}$ of $\tilde{\Omega}$ at $\omega$ consists of all complex vectors

$$x = \left( \delta\tilde{q}_i^0,\, \delta v_i^0,\, \delta m_i,\, \delta L \,\middle|\, i=1,\ldots,N \right) \in \mathbb{C}^{(2\nu+1)N+1};$$

thus $\mathcal{T}_\omega \tilde{\Omega}$ can be naturally identified with the complex vector space $\mathbb{C}^{(2\nu+1)N+1}$.
Set
(3.21)
$$\mathcal{N}(\omega) = \left\{ x = \left( \delta \tilde{q}_i^0, \, \delta v_i^0, \, \delta m_i, \, \delta L \,|\, i = 1, \ldots, N \right) \in \mathcal{T}_\omega \tilde{\Omega} \,\Big|\right.$$
$$\left. D_x(v_i^k) = D_x(m_i) = D_x(L) = 0 \text{ for } i = 1, \ldots, N;\, k = 0, \ldots, n \right\},$$

where $D_x(.)$ denotes the directional differentiation in the direction of $x$. Clearly, the complex linear subspace $\mathcal{N}(\omega)$ of $\mathbb{C}^{(2\nu+1)N+1}$ is the proper complex analogue of the neutral space of a genuine, real orbit segment.

*Remark.* Since $\delta v_i^0$, $\delta m_i$, and $\delta L$ must be equal to zero for a neutral vector $x \in \mathcal{N}(\omega)$, we shall simply write $x = (\delta \tilde{q}_1^0, \ldots, \delta \tilde{q}_N^0)$ instead of indicating the zero entries.

The proof of the following proposition is completely analogous to the one for the real case, and therefore we omit it.

PROPOSITION 3.22. *For every tangent vector $x \in \mathcal{N}(\omega)$ ($\omega \in \tilde{\Omega}(\Sigma, \mathcal{A}, \vec{\tau})$) and $1 \le k \le n$,*

(3.23) $$D_x \left( \tilde{q}_{i_k}^{k-1} - \tilde{q}_{j_k}^{k-1} \right) = (\alpha_k(x) - \alpha_{k-1}(x)) \cdot \left( v_{i_k}^{k-1} - v_{j_k}^{k-1} \right),$$

(3.24) $$D_x(t_k) = -\alpha_k(x),$$

*and*

(3.25) $$D_x(\tau_k) = D_x(t_k - t_{k-1}) = \alpha_{k-1}(x) - \alpha_k(x),$$

*where, by definition, $t_0 = 0 = \alpha_0(x)$. The functions $\alpha_k : \mathcal{N}(\omega) \to \mathbb{C}$ ($k = 1, \ldots, n$) are linear functionals. The name of the functional $\alpha_k$ is the advance of the $k^{\text{th}}$ collision $\sigma_k = \{i_k, j_k\}$; see also Section 2.*

This result markedly shows that, indeed, the vector space $\mathcal{N}(\omega)$ describes the linearity of the $(\Sigma, \mathcal{A}, \vec{\tau})$-iterated billiard map.

*Remark* 3.26. It is a matter of simple computation to convince ourselves that all assertions in Section 2 of the present paper and in [Sim(1992)-II] pertaining to the neutral space $\mathcal{N}(\omega)$ remain valid for the complexified dynamics. Here we only point out the two most important statements from those results:

(1) The vector space $\mathcal{N}(\omega)$ measures the ambiguity in determining the orbit segment $\omega$ purely by its velocity history $\{v_i^k(\omega) : i = 1, \ldots, N; k = 0, \ldots, n\}$ and the outer geometric data $(\vec{m}, L)$. In other words, this means that if, locally in the phase space $\tilde{\Omega}$ (in a small open set), two phase points $\omega_1$ and $\omega_2$ have the same velocities $(v_i^k(\omega_1) = v_i^k(\omega_2), i = 1, \ldots, N; k = 0, \ldots, n)$ and geometric parameters $(\vec{m}, L)$, then the initial data of $\omega_1$ and $\omega_2$ can only differ by a spatial translation by a vector from the space $\mathcal{N}(\omega_1) = \mathcal{N}(\omega_2)$. This statement is obviously reversible.

(2) The reflection laws for the neutral vectors are exactly the same as for velocities.

The Connecting Path Formula (CPF, Lemma 2.9 of [Sim(1992)-II], or Proposition 2.19 in this paper) is applicable, and it enables us to compute the neutral space $\mathcal{N}(\omega)$ via solving a homogeneous system of linear equations (now over the field $\mathbb{C}$)

$$(3.27) \qquad \sum_{k=1}^{n} \alpha_k \Gamma_{ik} = 0, \; i = 1, \ldots, n + P_\Sigma - N,$$

where each equation in (3.27) is a $\nu$-dimensional complex vector equation, the coefficients $\Gamma_{ik} = \Gamma_{ik}(\omega) \in \mathbb{C}^\nu$ are certain linear combinations of relative velocities $v_i^k(\omega) - v_j^k(\omega)$ (the coefficients of those linear combinations are just fractional linear expressions of the masses, see Definitions 2.17–2.18) and, finally, $P_\Sigma$ denotes the number of connected components of the collision graph of $\Sigma$. For a more detailed explanation of (3.27), see the beginning part of Remark 4.6 in the next section.

*Remark* 3.28. We note here that the ordering of the moments of collisions $t(f_i)$ (which plays a significant role in Definitions 2.17–2.18) is no longer meaningful over the unordered field $\mathbb{C}$. However, the use of the corresponding inequalities in 2.17–2.18 is purely technical/notational. Those inequalities only serve to introduce the combinatorial ordering of the collisions $f_i$ (which are $\sigma$'s) given by the indexes of the symbolic collision sequence $\Sigma = (\sigma_1, \sigma_2, \ldots, \sigma_n)$.

*Remark* 3.29. As a straightforward consequence of equation (3.6) in Proposition 3.3, we obtain that $v_i^{k-1} = v_i^k$ for $i = 1, \ldots, N$ (i.e. the $k^{\text{th}}$ reflection does not change the compound velocity) if and only if

$$\left\langle v_{i_k}^{k-1} - v_{j_k}^{k-1}; \tilde{q}_{i_k}^k - \tilde{q}_{j_k}^k - L \cdot a_k \right\rangle = 0,$$

and this is just the case of a tangential collision. An easy calculation shows that the mentioned tangentiality occurs if and only if the two roots $\tau_k$ of (3.8) coincide.

It is obvious that $\dim_{\mathbb{C}} \mathcal{N}(\omega)$ is at least $\nu + 1$ (as long as not all velocities are the same), because the flow direction and the uniform spatial translations are necessarily contained by $\mathcal{N}(\omega)$.

*Definition* 3.30 (*Sufficiency*). The *orbit segment* $\omega \in \tilde{\Omega}(\Sigma, \mathcal{A}, \vec{\tau})$ is said to be sufficient if and only if $\dim_{\mathbb{C}} \mathcal{N}(\omega) = \nu + 1$.

Finally, the main result of the next section (Key Lemma 4.1) will use the following:

*Definition* 3.31. The triple $(\Sigma, \mathcal{A}, \vec{\tau})$ has Property (A) if and only if the following assertion holds:

For every pair of indices $1 \leq k < l \leq n$ for which $\sigma_k = \sigma_l$, $a_k = a_l$, and $\sigma_j \cap \sigma_k = \emptyset$ for $j = k+1, \ldots, l-1$, we have that $\tau_l \not\equiv -\sum_{j=k+1}^{l-1} \tau_j$ in the algebraic selection $\vec{\tau}$ of the time slots. (Note that in this case the quadratic polynomial $b_l x^2 + c_l x + d_l$ is reducible over the field $\mathbb{K}_{l-1}$ because the time slot $\tau_l \equiv -\sum_{j=k+1}^{l-1} \tau_j \in \mathbb{K}_{l-1}$ is automatically a solution of (3.8).)

*Remark* 3.32. It is obvious that the schemes $(\Sigma, \mathcal{A}, \vec{\tau})$ corresponding to genuine, real billiard trajectory segments enjoy Property (A). Indeed, in those cases even the equations $\sigma_k = \sigma_l$, $a_k = a_l$, and $\sigma_j \cap \sigma_k = \emptyset$ $(j = k+1, \ldots, l-1)$ cannot hold simultaneously. Concerning the necessity for introducing Property (A), see Example I at the end of the next section.

*Remark* 3.33. For combinatorial schemes satisfying Property (A), the quadratic polynomial $b_l x^2 + c_l x + d_l$ is reducible over the field $\mathbb{K}_{l-1}$ (i.e. $\mathbb{K}_{l-1} = \mathbb{K}_l$) only in the case described above in 3.31. This will be a consequence of Main Lemma 4.21 (cf. Corollary 4.36).

## 4. The notion of richness. Richness implies sufficiency

The result of this section is

KEY LEMMA 4.1. *There exists a positive number $C(N)$ (depending merely on the number of balls $N \geq 2$) with the following property*: *If a symbolic sequence $\Sigma = (\sigma_1, \ldots, \sigma_n)$ is $C(N)$-rich and the triple $(\Sigma, \mathcal{A}, \vec{\tau})$ has Property* (A) *(see Definition* 3.31*), then* $\dim_{\mathbb{C}} \mathcal{N}(\omega) = \nu + 1$ *(i.e. the trajectory is sufficient) for almost every $\omega \in \tilde{\Omega}(\Sigma, \mathcal{A}, \vec{\tau})$. The real version of this result is also valid. (For the definition of a $C$-rich symbolic sequence, see* 2.5.*)*

*Remark* 4.1/a. Plainly, the exceptional zero measure subset of $\tilde{\Omega}(\Sigma, \mathcal{A}, \vec{\tau})$ mentioned in this lemma must actually be a countable union of proper, analytic submanifolds.

The aim of Sections 3 and 4, in particular that of Key Lemma 4.1, is to settle part (1) of the strategy detailed. The methods are mainly algebraic, and for an easier separation of our dynamical vs. algebraic arguments, we are now going to describe *what are the only properties of the dynamics (and nothing more) that we use.*

(1) The dynamical properties of an orbit segment $S^{[0,T]} x_0$ with a symbolic collision sequence $\Sigma$ are completely characterized by equations (3.4)–(3.8) if we, in addition, specify $\tau_k$ as the smaller of the two, necessarily positive, roots of equation (3.8). In the inductive steps of Lemmas 4.9, 4.39, and 4.40 we also use the remark to Proposition 3.3 saying that the collision laws (3.9), themselves consequences of (3.5)–(3.6), extend analytically to the case of zero mass; cf. (3.10).

(2) Since our goal is to show that orbits with a rich collision structure are typically sufficient or, equivalently, their neutral subspaces have typically minimal possible dimension, we also need an algebraic characterization of the neutral subspace. This is attained through the CPF (see Proposition 2.19 and the system of linear equations (3.27)), and our arguments also use some simple inferences about the CPF and the neutral subspaces from [Sim(1992)-II] .

(3) As we have seen in Section 3, the algebraic formalisms of points (1) and (2) can be complexified, and our forthcoming discussion is based upon this complexification.

*Remark* 4.1/b. In the course of proving the key lemma we will see that the thresholds $C(N)$ ($N \geq 2$) can actually be defined by the recursion $C(2) = 1$, $C(N) = \frac{N}{2} \cdot \max\{C(N-1), 3\}$ ($N \geq 3$); that is, $C(2) = 1$, $C(N) = \dfrac{3N!}{2^{N-1}}$ for $N \geq 3$. We note that, as proved in [K-S-Sz(1991)], for $N = 3$ the value $C(3) = 2$ already does the job.

*Proof of Key Lemma* 4.1. All the rest of this section is devoted to the inductive proof of the key lemma, and which be split into several lemmas.

First of all, we prove the original, complex version of the key lemma by an induction on the number $N$. Once this is completed, we will show that the validity of the key lemma carries over to the real case.

### The dichotomy. Algebraic characterization of sufficiency

LEMMA 4.2. *Consider the open, connected and dense subset $D$ of $\mathbb{C}^{(2\nu+1)N+1}$ formed by all possible complex $(2\nu+1)N+1$-tuples*

$$\vec{x}(\omega) = \left( \left(\tilde{q}_i^0(\omega)\right)_j, \left(v_i^0(\omega)\right)_j, m_i(\omega), L(\omega) \,\middle|\, i = 1, \ldots, N;\, j = 1, \ldots, \nu \right) \quad \left(\omega \in \tilde{\Omega}\right)$$

*of initial data of the complex $(\Sigma, \mathcal{A}, \vec{\tau})$-dynamics; see Definitions* 3.18 *and* 3.20. *We claim that there are finitely many complex polynomials $P_1(\vec{x}), \ldots, P_s(\vec{x})$ ($\vec{x} \in \mathbb{C}^{(2\nu+1)N+1}$), canonically determined by the discrete parameters $(\Sigma, \mathcal{A}, \vec{\tau})$, such that for every $\vec{x} \in D(\Sigma, \mathcal{A}, \vec{\tau})$ the following implication holds true*:

(4.3)
$$\left( \exists\, \omega \in \tilde{\Omega}(\Sigma, \mathcal{A}, \vec{\tau}) \text{ such that } \vec{x} = \vec{x}(\omega) \text{ and } \dim_{\mathbb{C}} \mathcal{N}(\omega) > \nu + 1 \right)$$
$$\implies P_1(\vec{x}) = P_2(\vec{x}) = \cdots = P_s(\vec{x}) = 0.$$

*Remark.* The left-hand side in the above equivalence precisely says that *some* branch of the multiple-valued $(\Sigma, \mathcal{A}, \vec{\tau})$-dynamics with initial data $\vec{x}$ is not sufficient, see also Remark 2 to Definition 3.18.

*Proof.* This proof is a typical application of the Connecting Path Formula (see Proposition 2.19 here or Lemma 2.9 and Proposition 3.4 of [Sim(1992)-II]). The latter asserts that

(4.4) $$\dim_{\mathbb{C}} \mathcal{N}(\omega) = \nu \cdot P_{\Sigma} + \dim_{\mathbb{C}} \{\alpha_1, \ldots, \alpha_n\},$$

where, as said before, $P_{\Sigma}$ denotes the number of connected components of the collision graph of $\Sigma$, and $\{\alpha_1, \ldots, \alpha_n\}$ is shorthand for the complex linear space of all possible $n$-tuples $(\alpha_1(x), \ldots, \alpha_n(x))$ of advances, $x \in \mathcal{N}(\omega)$.

*Remark* 4.5. It is worth noting here that our present formula (4.4) differs from its counterpart in Proposition 3.4 of [Sim(1992)-II] by an additional term $\nu$. This is, however, due to the fact that in the present approach we no longer have the reduction equation $\sum_{i=1}^{N} m_i \delta \tilde{q}_i = 0$. It follows easily from (4.4) that the sufficiency is equivalent to $\dim_{\mathbb{C}}\{\alpha_1, \ldots, \alpha_n\} = 1$; i.e., sufficiency means that all advances are equal to the same functional. Note that $\dim_{\mathbb{C}}\{\alpha_1, \ldots, \alpha_n\} = 1$ obviously implies $P_{\Sigma} = 1$.

A simple, but important consequence of the Connecting Path Formula is that the linear space $\{\alpha_1, \ldots, \alpha_n\}$ is the solution set of a homogeneous system of linear equations (3.27).

*Remark* 4.6. Note that, when applying the mentioned lemma, the left-hand side of the Connecting Path Formula has to be written as the relative velocity of the colliding balls multiplied by the advance of that collision. It follows immediately from the exposition of [Sim(1992)-II] that equations of this sort (arising from all CPF's) are the *only* constraints on the advance functionals. The reason why this is true is that the fulfillment of all CPF's precisely means that the relative displacement (variation of position) of every pair of particles right before their new collision is parallel to the relative incoming velocity of these particles (see also (3.23)), and this fact guarantees that the variation of the newly formed relative outgoing velocity will also be zero, just as required in (3.21).

Also note that the number $n + P_{\Sigma} - N$ of equations in (3.27) is equal to the number of all collisions $\sigma_k = \{i_k, j_k\}$ for which $i_k$ and $j_k$ are in the same connected component of the collision graph of $\{\sigma_1, \sigma_2, \ldots, \sigma_{k-1}\}$.

The determinant $D(M)$ of every minor $M$ of the coefficient matrix of (3.27) is a homogeneous velocity polynomial and, therefore, $D(M)$ is a holomorphic function on $\tilde{\Omega}(\Sigma, \mathcal{A}, \vec{\tau})$. Thus, there are finitely many velocity polynomials (actually, determinants) $R_1(\omega), \ldots, R_s(\omega)$ such that their simultaneous vanishing is equivalent to the nonsufficiency of the orbit segment $\omega$. By using the equations (3.4)–(3.8) (which define the complex dynamics recursively), each of these velocity polynomials $R_i(\omega)$ can be – essentially uniquely – written as an algebraic function $f_i(\vec{x}(\omega)) = f_i(\vec{x})$ of the initial data $\vec{x}$, and these algebraic functions $f_i(\vec{x})$ only contain (finitely many) field operations and square roots. By using the canonical method of successive elimination of the square roots from the equation $f_i(\vec{x}) = 0$, one obtains a complex polynomial $P_i(\vec{x})$ such that

$$P_i(\vec{x}) = 0 \iff f_i(\vec{x}) = 0.$$

We emphasize here that

(1) we understand the relation $f_i(\vec{x}) = 0$ in such a way that $f_i(\vec{x})$ becomes zero on *some* branch of the square root function when evaluating the multiple-valued algebraic function $f_i(\vec{x})$; cf. Remark 2 to Definition 3.18;
(2) the equivalence is claimed to only hold on $D$; cf. the formulation of Lemma 4.2.

The process of eliminating the square roots from an equation $f(\vec{x}) = 0$ (of the above type), however, requires a bit of clarification. The algebraic function $f$ has a natural representing element $\hat{f}$ in the field $\mathbb{K}_n$ (see Definition 3.11). Suppose that $\hat{f} \in \mathbb{K}_k$ but $\hat{f} \notin \mathbb{K}_{k-1}$. Then the field element $\hat{f} \in \mathbb{K}_k$ can be written uniquely in the form $\hat{f} = a + b\sqrt{\delta}$, where $a, b, \delta \in \mathbb{K}_{k-1}$, $b \neq 0$, $\sqrt{\delta} \notin \mathbb{K}_{k-1}$, and

$$\delta = \left\langle v_{i_k}^{k-1} - v_{j_k}^{k-1};\, \tilde{q}_{i_k}^{k-1} - \tilde{q}_{j_k}^{k-1} - L \cdot a_k \right\rangle^2 \\ - \left\| v_{i_k}^{k-1} - v_{j_k}^{k-1} \right\|^2 \cdot \left( \left\| \tilde{q}_{i_k}^{k-1} - \tilde{q}_{j_k}^{k-1} - L \cdot a_k \right\|^2 - 4r^2 \right)$$

is the discriminant of the quadratic equation (3.8) with the unknown $\tau_k$.

In the first step of the successive elimination of the square roots from the equation $f(\vec{x}) = 0$, we switch from the field element $\hat{f} = a + b\sqrt{\delta} \in \mathbb{K}_k$ (representing $f$) to the new field element

$$\hat{f}' = \left(a + b\sqrt{\delta}\right) \cdot \left(a - b\sqrt{\delta}\right) = a^2 - b^2\delta \in \mathbb{K}_{k-1},$$

which is a representing element of an algebraic function $f'(\vec{x})$ containing square roots than $f(\vec{x})$. Corresponding to the switching $\hat{f} \mapsto \hat{f}'$, we also replace the original equation $f(\vec{x}) = 0$ by the new equation $f'(\vec{x}) = 0$. Then we continue this process until we arrive at a field element

$$\hat{g}(\vec{x}) = \frac{P(\vec{x})}{Q(\vec{x})} \in \mathbb{K}_0 = \mathbb{C}(\vec{x}),$$

which is a complex rational function $\dfrac{P(\vec{x})}{Q(\vec{x})}$ of the initial variables $\vec{x}$. Of course, here we assume that the polynomials $P(\vec{x})$ and $Q(\vec{x})$ have no nontrivial common divisor.

Actually, as directly follows from Proposition 3.3 and from the method of eliminating the square roots (presented above), in each denominator $Q_l(\vec{x})$ ($1 \leq l \leq s$) associated to the algebraic function $f_l(\vec{x})$ there are only factors of the type $m_i + m_j$ or $\|v_i^k - v_j^k\|^2$, $1 \leq i < j \leq N$, $0 \leq k \leq n$. Therefore, none of the denominators $Q_l(\vec{x})$ vanishes on $D$. This explains property (2) formulated above.

As it also follows from the method of eliminating the square roots and from the actual form of the equations (3.4–3.8), the coefficients of the polynomials $P_i(\vec{x})$ are real.

By using some simple facts from Galois theory, we can easily detect an intimate relationship between the rational function $\hat{g}(\vec{x}) \in \mathbb{K}_0$ just constructed and the product $\alpha$ of all conjugate elements of $\hat{f} \in \mathbb{K}_n$. (The conjugates of $\hat{f}$ (over the base field $\mathbb{K}_0$) are sitting in the normal hull $\overline{\mathbb{K}_n} : \mathbb{K}_0$ of the field extension $\mathbb{K}_n : \mathbb{K}_0$. The element $\alpha \in \mathbb{K}_0$ is obviously the constant term in the normalized minimal polynomial of $\hat{f}$ over the field $\mathbb{K}_0$. For the elementary concepts and facts from Galois theory, the reader is referred to the book [St(1973)].) By establishing a natural matching between all conjugates of $a^2 - b^2\delta$ and the conjugates of the ordered pair $(a + b\sqrt{\delta}, a - b\sqrt{\delta})$, we can easily see that there exists an integer $l$ such that $\hat{g} = \alpha^{2^l}$. (In case $l < 0$ this should be understood as $\alpha = (\hat{g})^{2^{-l}}$.) Especially, $\alpha$ and $\hat{g}$ have the same irreducible factors both in the numerator and in the denominator. Another consequence is that the field element $\hat{f}$ is zero if and only if $P(\vec{x}) \equiv 0$.

So far we have seen that for each velocity determinant $R_i(\omega)$ ($i = 1, \ldots, s$) one can find a *canonically determined* polynomial $P_i(\vec{x})$ with the following property: For every vector $\vec{x} \in D(\Sigma, \mathcal{A}, \vec{\tau})$ the equation $P_i(\vec{x}) = 0$ is equivalent to the equation $f_i(\vec{x}) = 0$ *on some branch of the square root function when evaluating the algebraic function* $f_i(\vec{x})$. Hence Lemma 4.2 follows. □

COROLLARY 4.7 (dichotomy corollary). *As another consequence of the above identity $\hat{g} = \alpha^{2^l}$, we obtain that if $P \equiv 0$, then $f(\vec{x}) = 0$ for every $\vec{x} \in \mathbb{C}^{(2\nu+1)N+1}$ and on every branch of the square root function when evaluating the multiple-valued algebraic function $f(\vec{x})$. This also means that if all polynomials $P_1, \ldots, P_s$ in Lemma 4.2 happen to be identically zero, then no orbit segment $\omega \in \tilde{\Omega}(\Sigma, \mathcal{A}, \vec{\tau})$ is sufficient.*

*Therefore, thanks to the algebraic feature of the dynamics, there is a dichotomy: Either every $(\Sigma, \mathcal{A}, \vec{\tau})$-orbit segment is nonsufficient, or almost every such trajectory segment is sufficient.*

*Remark.* We note that the reason why we do not have an equivalence in (4.3) is that – for a typical $\vec{x} \in \mathbb{C}^{(2\nu+1)N+1}$ – the branches, on which $f_i(\vec{x}) = 0$, may be different for different values of $i$.

According to the above lemma, in order to prove Key Lemma 4.1, it is enough to show that at least one of the polynomials $P_i(\vec{x})$ is nonzero. This is exactly what we are going to do by an induction on the number $N \geq 2$.

LEMMA 4.8. *Pertaining to the polynomials $P_i(\vec{x})$ associated with $R_i(\omega) = f_i(\vec{x}(\omega))$ in Lemma 4.2, the following holds*:

*If the algebraic function $f_i(\vec{x})$ is identically zero on a nonempty, open subset $B$ of $\mathbb{R}^{(2\nu+1)N+1}$, then $f_i$ is algebraically trivial, i.e., $\hat{f}_i = 0$ or, equivalently, $P_i(\vec{x}) \equiv 0$.*

*Proof.* Denote the minimal polynomial of $\hat{f}_i$ over the base field $\mathbb{K}_0 = \mathbb{C}(\vec{x})$ by $M(x) = \sum_{j=0}^{k} c_j x^j$, $c_j \in \mathbb{K}_0$. Suppose, on the contrary, that $f_i$ vanishes on an open ball $B$ of $\mathbb{R}^{(2\nu+1)N+1}$. Since $f_i$ fulfills the identity $\sum_{j=0}^{k} c_j (f_i)^j \equiv 0$, the vanishing of $f_i$ on the ball $B$ implies that the element $c_0$ of the rational function field $\mathbb{K}_0$ vanishes on $B$, as well. Therefore, $c_0 = 0$ in $\mathbb{K}_0$, so that $M(x)$ is divisible by $x$. The irreducibility of $M$ yields $M(x) = x$, i.e., $\hat{f}_i = 0$. Hence the lemma follows. □

## The substitution $\{m_N = 0\}$. Derived schemes

The next lemma describes the $(\Sigma, \mathcal{A}, \vec{\tau})$-orbit segments in the case when the $N^{\text{th}}$ particle is "infinitely light" compared to the others, i.e., when $m_N = 0$ but $m_1 \cdot \ldots \cdot m_{N-1} \neq 0$.

LEMMA 4.9. *Suppose that $m_N = 0$ and $m_1 \cdot \ldots \cdot m_{N-1} \neq 0$. Assume further that there exist indices $1 \leq p < q \leq n$ such that $N \in \sigma_p$, $N \in \sigma_q$ and $N \notin \sigma_j$ for $j = p+1, \ldots, q-1$. Then an orbit segment $\omega \in \tilde{\Omega}(\Sigma, \mathcal{A}, \vec{\tau})$ is sufficient (i.e. $\dim_\mathbb{C} \mathcal{N}(\omega) = \nu + 1$) if*

(1) *the $\{1, 2, \ldots, N-1\}$-part*

$$\left\{\tilde{q}_i^k(\omega), v_i^k(\omega), m_i(\omega), L(\omega) \mid i = 1, \ldots, N-1; k \in I_N\right\}$$

*of the orbit segment is sufficient as an orbit segment of the particles $1, \ldots, N-1$ and*

(2) *the relative velocities $v_N^p(\omega) - v_{i_p}^p(\omega)$ and $v_N^{q-1}(\omega) - v_{i_q}^{q-1}(\omega)$ are not parallel.*

(Here, as usual, $i_p$ ($i_q$) is the index of the ball colliding with the $N^{\text{th}}$ particle at $\sigma_p$ ($\sigma_q$), and the index set $I_N \subset \{0, 1, \ldots, n\}$ contains 0 and those indices $i > 0$ for which $N \notin \sigma_i$.)

*Remark.* It is worth noting that, as follows easily from the conditions, $v_N^p(\omega) = v_N^{q-1}(\omega)$.

*Proof.* First of all, we observe that the motion of the $N^{\text{th}}$ ball with zero mass has absolutely no effect on the evolution of the $\{1, \ldots, N-1\}$-part of the trajectory segment, and this statement is also valid for the time evolution of the $\{1, \ldots, N-1\}$-part $\{\delta \tilde{q}_i^k \mid i = 1, \ldots, N-1\}$ of a neutral vector $(\delta \tilde{q}_1^0, \ldots, \delta \tilde{q}_N^0) \in \mathcal{N}(\omega)$; see also the equations (3.4), (3.7)–(3.8) and (3.10).

Consider now an arbitrary neutral vector $\delta Q = (\delta \tilde{q}_1^0, \ldots, \delta \tilde{q}_N^0) \in \mathcal{N}(\omega)$. According to (1) and the above principle, we can modify the neutral vector $\delta Q$ by a scalar multiple of the flow direction and by a uniform spatial displacement in such a way that $\delta \tilde{q}_i^0 = 0$ for $i = 1, \ldots, N-1$, and then $\delta \tilde{q}_i^k = 0$ remains true for the whole orbit segment, $i = 1, \ldots, N-1$, $k = 0, \ldots, n$. The neutrality of $\delta Q$ with respect to the entire orbit segment means, however, that

(4.10) $\quad \delta \tilde{q}_N^p = \alpha_p \cdot \left( v_N^p(\omega) - v_{i_p}^p(\omega) \right)$ and $\delta \tilde{q}_N^{q-1} = \alpha_q \cdot \left( v_N^{q-1}(\omega) - v_{i_q}^{q-1}(\omega) \right)$.

According to our hypothesis (2), the consequence of (4.10) and of the obvious equation $\delta \tilde{q}_N^p = \delta \tilde{q}_N^{q-1}$ is that $\alpha_p = \alpha_q = 0$, and thus $\delta Q = 0$. Hence the lemma follows. □

The crucial part of the inductive proof of Key Lemma 4.1, i.e. the substitution $m_N = 0$ (and $L = 0$, too), requires some preparatory thoughts and lemmas. The formulations given for the substitution $m_N = 0$ can be easily adapted to the substitution $L = 0$. Assume that a combinatorial-algebraic scheme $(\Sigma, \mathcal{A}, \vec{\tau})$ is given for the $N$-ball system $\{1, 2, \ldots, N\}$ so that Property (A) holds (see Definition 3.31). In the proof of Key Lemma 4.1 we want to use some combinatorial-algebraic $(N-1)$-schemes $(\Sigma', \mathcal{A}', \vec{\tau}')$ that

(i) govern the time evolution of the $\{1, 2, \ldots, N-1\}$-part of *some* $(\Sigma, \mathcal{A}, \vec{\tau})$-orbit segments with an infinitely light $N^{\text{th}}$ ball, i.e. $m_N = 0$ and

(ii) enjoy Property (A).

In the construction of such $(N-1)$-schemes $(\Sigma', \mathcal{A}', \vec{\tau}')$, the so-called *derived schemes*, we want to follow the guiding principles below, which also serve as the definition of the *derived schemes*:

*Definition* 4.11. We say that the $(N-1)$-scheme $(\Sigma', \mathcal{A}', \vec{\tau}')$ is a scheme derived from $(\Sigma, \mathcal{A}, \vec{\tau})$ by putting $m_N = 0$ if it is obtained as follows:

(1) First of all, we discard all symbols $\sigma_j$, $a_j$, and $\tau_j$ from $(\Sigma, \mathcal{A}, \vec{\tau})$ that correspond to a $\sigma_j$ containing the label $N$. This also means that we retain the other symbols $\sigma_j$, $a_j$ without change: they only get re-indexed, due to the dropping of the other symbols;

(2) As far as the selection of $\vec{\tau}'$ is concerned, we want to retain *all algebraic relations* among the variables that are encoded in the original scheme $(\Sigma, \mathcal{A}, \vec{\tau})$, i.e. all algebraic relations that follow from the scheme $(\Sigma, \mathcal{A}, \vec{\tau})$ and from $m_N = 0$;

(3) We want the derived scheme $(\Sigma', \mathcal{A}', \vec{\tau}')$ to enjoy Property (A); see 3.31.

It is straightforward that, following just the instructions in (1)-(2) above, one can easily construct such schemes $(\Sigma, \mathcal{A}, \vec{\tau})$ which, perhaps, do not have Property (A).

*Remark* 4.12. Part (2) of the above definition means the following: Initially, the algebra of $(\Sigma, \mathcal{A}, \vec{\tau})$ consists

(1) of the variables
$$\left\{ \left(\tilde{q}_i^k\right)_j, (v_i^k)_j, m_i, L, \tau_l \Big| i = 1, \ldots, N; \ k = 0, \ldots, n; \right.$$
$$\left. j = 1, \ldots, \nu; \ l = 1, \ldots, n \right\};$$

(2) of the equations (3.4–3.8), and finally
(3) of the choices of the signs of the square roots in the solutions of the quadratic equations (3.8) whenever the equations are reducible (cf. Remark 3.13).

To obtain the algebra of the derived scheme $(\Sigma', \mathcal{A}', \vec{\tau}')$, we put first $m_N = 0$ and then $\tau_l' = \sum_{j=k+1}^{l} \tau_j$ whenever for some $1 \leq k < l \leq n$ we have $N \notin \sigma_k \cup \sigma_l$, $N \in \bigcap_{j=k+1}^{l-1} \sigma_j$ (here the indexing of the time slots $\tau_l'$ of $(\Sigma', \mathcal{A}', \vec{\tau}')$ is lacunary). After trivial substitutions, the quadratic equation for $\tau_l$ leads to the natural quadratic equation for $\tau_l'$, and then, finally, we can cancel all the variables $\{(\tilde{q}_N^k)_j, (v_N^k)_j, m_N | \ 0 \leq k \leq n, 1 \leq j \leq \nu\}$, and all the equations they figure in as well. The collection of the remaining variables and equations with the corresponding vector $\vec{\tau}'$ make $(\Sigma', \mathcal{A}', \vec{\tau}')$.

*Remark* 4.13. Below we show the inheritance of the reducibility $\mathbb{K}_{l-1} = \mathbb{K}_l \Longrightarrow \mathbb{K}'_{l-1} = \mathbb{K}'_l$. We note, however, that this result *is not used* in the proof of Key Lemma 4.1. We merely include it – along with its nice consequence: part (3) of Lemma 4.37 – in order to give the reader a broader view of the subject. The subsequent proof of Lemma 4.40 only uses the weaker part (2) of 4.37 which, in turn, does not need the above mentioned inheritance of reducibility.

Another interpretation of part (2) of the above definition is the following one: Once we have a complex rational expression $g(z_1, \ldots, z_p)$ of arbitrary elements
$$z_1, \ldots, z_p \in \mathbb{K}'_n = \mathbb{K}\left(\Sigma', \mathcal{A}', \vec{\tau}'\right)$$
so that the same expression of these variables is zero in the field $\mathbb{K}_n = \mathbb{K}(\Sigma, \mathcal{A}, \vec{\tau})$, then it must also be zero in the field $\mathbb{K}'_n$. An important case (when (2) has special significance) is when $N \notin \sigma_k \cup \sigma_l$ ($1 \leq k < l \leq n$), $N \in \bigcap_{j=k+1}^{l-1} \sigma_j$, and $\mathbb{K}_{l-1} = \mathbb{K}_l$. Then the time slot $\tau_l' = \sum_{j=k+1}^{l} \tau_j$ (here the indexing of the time slots $\tau_l'$ of $(\Sigma', \mathcal{A}', \vec{\tau}')$ is lacunary.) should be expressed as a *rational function* of the earlier defined $(\Sigma, \mathcal{A}, \vec{\tau})$-variables (the variables constituting the field $\mathbb{K}_{l-1}$), just as it is indicated by the rational dependence of $\tau_l$ on the variables constituting the field $\mathbb{K}_{l-1}$. Therefore, $\tau_l' = \sum_{j=k+1}^{l} \tau_j$ is a rational expression

$$(4.14) \qquad \tau_l' \equiv g\left(\vec{x}; \tilde{q}_N^0, v_N^0, m_N; \tau_1, \ldots, \tau_{l-1}\right)$$

of the initial variables $\left(\vec{x}; \tilde{q}_N^0, v_N^0, m_N\right) \in \mathbb{C}^{(2\nu+1)N+1}$ and the time slots $\tau_1, \ldots, \tau_{l-1}$. Note that the identity (4.14) is also true if we impose additionally the side condition $m_N = 0$:

$$(4.15) \qquad \tau_l'^0 \equiv g\left(\vec{x}; \tilde{q}_N^0, v_N^0, 0; \tau_1^0, \ldots, \tau_{l-1}^0\right).$$

However, the $\{1, \ldots, N-1\}$-dynamics evolves independently of the data of the $N^{\text{th}}$ ball with zero mass $m_N$. This means that

$$g\left(\vec{x}; \tilde{q}_N^0(1), v_N^0(1), 0; \tau_1^0(1), \ldots, \tau_{l-1}^0(1)\right)$$
$$= g\left(\vec{x}; \tilde{q}_N^0(2), v_N^0(2), 0; \tau_1^0(2), \ldots, \tau_{l-1}^0(2)\right)$$

whenever all cumulative sums $\sum_{j=k_1+1}^{k_2} \tau_j^0(1)$ are equal to the corresponding cumulative sums $\sum_{j=k_1+1}^{k_2} \tau_j^0(2)$ for $1 \leq k_1 < k_2 < l$, $N \notin \sigma_{k_1} \cup \sigma_{k_2}$, $N \in \bigcap_{j=k_1+1}^{k_2-1} \sigma_j$. This precisely means that the rational function $g(\vec{x}; \tilde{q}_N^0, v_N^0, 0; \tau_1^0, \ldots, \tau_{l-1}^0)$

(i) does not contain the variables $\left(\tilde{q}_N^0\right)_\beta$, $\left(v_N^0\right)_\beta$ ($1 \leq \beta \leq \nu$), and

(ii) only depends on the above cumulative sums $\sum_{j=k_1+1}^{k_2} \tau_j^0 = \tau_{k_2}'^0$.

(The indexing of $\tau_j'^0$ is still lacunary!) Thus, we get that $g(\vec{x}; \tilde{q}_N^0, v_N^0, 0; \tau_1^0, \ldots, \tau_{l-1}^0)$ is, in fact, a rational expression of $\vec{x}$ and the time slots $\tau_j'^0$ preceding $\tau_l'^0$. This shows that $\tau_l'^0 \in \mathbb{K}_{l-1}'$, i.e. $\mathbb{K}_{l-1}' = \mathbb{K}_l'$ follows from the similar reducibility $\mathbb{K}_{l-1} = \mathbb{K}_l$.

The main point is the following one: Once $\tau_l$ is the root of a reducible polynomial, it is selected "artificially", following the instructions encoded in the scheme $\vec{\tau}$. Now we want to avoid the selection of the other root when defining the derived scheme $(\Sigma', \mathcal{A}', \vec{\tau}')$. This is why (2) is a requirement, and not a corollary, and this is why, in Remark 4.12, we use the wording "the time slot $\tau_l'$ should be expressed."

The fact that the reducibility of the polynomial defining $\tau_l$ is inherited during a substitution $m_N = 0$ followed from a general principle: If some algebraic function in $\mathbb{K}_{l-1}$ (the discriminant of the quadratic polynomial defining $\tau_l$) is the square of another such algebraic expression, then this square relation remains even "more true" if one annihilates a variable, i.e. considers that square relation on a smaller set. It is now the nature of our model that $m_N = 0$ immediately kills all variables involving the motion of the $N^{\text{th}}$ particle, since the considered new discriminant of the quadratic equation defining $\tau'_l$ is obviously independent of those variables.

*Remark* 4.16. (a) Note that the construction of such schemes $(\Sigma', \mathcal{A}', \vec{\tau}')$, by just requiring (1)–(2), is far from unique: It is very possible that, originally, for some $k$ we have $|\mathbb{K}_k : \mathbb{K}_{k-1}| = 2$ (i.e. the polynomial (3.8) determining $\tau_k$ is irreducible), but after the substitution $m_N = 0$ the corresponding field extension collapses to a trivial, degree-one extension. (The paramount example for this phenomenon is just the set-up of the upcoming Main Lemma 4.21.) If this phenomenon of "reducible splitting" takes place, then we have to make our choice: Which one of the two (now distinguishable) roots of the reducible (3.8) to include in $\vec{\tau}'$?

The requirement (A) for the derived scheme $(\Sigma', \mathcal{A}', \vec{\tau}')$, however, ties our hands when selecting this root in the scenario of Main Lemma 4.21, and it actually implies the uniqueness of the derived scheme $(\Sigma', \mathcal{A}', \vec{\tau}')$ (see Corollary 4.36).

The phenomenon of this reducible splitting can also be well represented by the following, simple paradigm: Let $a \in \mathbb{C}$ be a complex parameter (playing the role of $m_N$). Consider the Riemann-leaf of the function $y = f(z) = \sqrt{z^2 + a}$, i.e. the algebraic variety defined by the equation $y^2 - z^2 - a = 0$. This is irreducible (undecomposable) for $a \neq 0$, and it splits into two components ($y - z = 0$ and $y + z = 0$) for $a = 0$. However, given a sequence $a_k \to 0$ ($a_k \neq 0$), every point $(z, y)$ on the curve $y^2 - z^2 = 0$ turns out to be a limit of a sequence $(z_k, y_k)$ with $y_k^2 - z_k^2 - a_k = 0$.

(b) The paradigm just given also makes it possible to illustrate why (2) in Definition 4.11 is a requirement and not a corollary. Indeed, consider the field $\mathbb{K}_0 = \mathbb{C}(z, m_N)$; i.e., $\mathbb{K}_0$ is the complex field extended by two transcendental elements $z$ and $m_N$. Let $\mathbb{K}_1 := \mathbb{K}_0(\tau_1) = \mathbb{K}_0(\sqrt{z^2 + 1 + m_N})$ be the extension of $\mathbb{K}_0$ by the solution $\tau_1$ of the irreducible equation $\tau_1^2 - z^2 - 1 - m_N = 0$. Let us now extend $\mathbb{K}_1$ by the root $\tau_2$ of the quadratic polynomial $x^2 - 5\tau_1 x + 6\tau_1^2$. This polynomial is reducible over the field $\mathbb{K}_1$: its two roots are $2\tau_1$ and $3\tau_1$. Now the extended field $\mathbb{K}_2$ will be just equal to $\mathbb{K}_1$, and the selection of a designated root (as $\tau_2$) of the above polynomial amounts to selecting either $2\tau_1$ or $3\tau_1$. Observe that the root of the considered reducible polynomial could have been an arbitrary rational expression of the elements $z$, and $\tau_1$. In a general reducible step, such a rational expression shows the unique way of computing the actual $\tau_n$ from the initial data and the $\tau_i$'s with $i = 1, \ldots, n - 1$. However, when assigning actual complex values to the field element $\tau_1 = \sqrt{z^2 + 1 + m_N}$, we have two possiblities. Summing everything up, our selection does not mean giving values to these $\tau$'s but, rather, telling if we take (formally) $2\tau_1$ or $3\tau_1$ as $\tau_2$. (This selection determines the discrete algebraic structure, which becomes transparent if we think about repeated $\sigma$'s with identical adjustment vectors.) When defining the derived scheme this means that once we have chosen $3\sqrt{z^2 + 1 + m_N}$ (instead of $2\sqrt{z^2 + 1 + m_N}$) we still want to choose $3\sqrt{z^2 + 1}$ after the substitution $m_N = 0$, and do not want to switch suddenly to the other root $2\sqrt{z^2 + 1}$.

## Transversalities of the degeneracies

*Definition-notation* 4.17. Denote by $\Phi_p$ ($1 \leq p \leq n$) the degeneracy $\left\| v_{i_p}^p - v_{j_p}^p \right\|^2 = 0$ and by $\Psi_p$ ($1 \leq p \leq n$) the degeneracy defined by $\left\langle v_{i_p}^p - v_{j_p}^p ; \tilde{q}_{i_p}^p - \tilde{q}_{j_p}^p - L \cdot a_p \right\rangle = 0$. (These degeneracies are represented here with post-collision velocities, but in view of Remark 1 to Definition 3.18 they also have an analogous representation in terms of pre-collision velocities.)

*Remark* 4.18. (a) In the following arguments we will use a slightly different working domain $D = D_{p+1,n-1}$, which is not exactly the one of Definition 3.18. Instead, in general, the working domain $D_{p+1,n-1}$ describes the orbit segments that are nondegenerate from the $(p+1)$-st collision up to the $(n-1)$-st one (inclusive). Here $p \leq n - 1$, and if $p = n - 1$, then, by definition, $D_{p+1,n-1} = \mathbb{C}^{(2\nu+1)N+1}$. The reference time (the upper index of the initial data) is now $p + 1$; i.e. the outgoing kinetic variables are taken at $\sigma_{p+1}$.

(b) Observe that if for some $p < p'$

(4.19) $\quad \sigma_p = \sigma_{p'}$ and, for $p < j < p'$, $\sigma_j \cap \sigma_p =$ either $\sigma_p$ or $\emptyset$,

then by Remark 1 to Definition 3.18 $\Phi_p = \Phi_{p'}$. Similarly, if for some $p < p'$
(4.20)
$(\sigma_p, a_p) = (\sigma_{p'}, a_{p'})$ and for $p < j < p'$ either $(\sigma_j, a_j) = (\sigma_p, a_p)$ or $\sigma_j \cap \sigma_p = \emptyset$,

then necessarily $\Psi_p = \Psi_{p'}$.

(c) In the upcoming Main Lemma transversality statements for $\Phi_p$ (or $\Psi_p$) and $\Phi_n$ (or $\Psi_n$) will be claimed on $D_{p+1,n-1}$. Because of our previous observations, these assertions can be nonvacuous only if the symbolic pair $(\Sigma, \mathcal{A})$ is minimal in the sense that

(1) if $\Phi_p$ is involved in the statement, then there is no $p' : p < p' < n$ such that (4.19) holds (this requirement is called $\Phi$-*minimality from the left*);
(2) if $\Psi_p$ is involved in the statement, then there is no $p' : p < p' < n$ such that (4.20) holds (this requirement is called $\Psi$-*minimality from the left*);

Analogous requirements can be formulated when $\Phi_n$ ($\Psi_n$) is involved, and will be called $\Phi$ (or $\Psi$)-*minimality from the right*.

In order to show that there really exist derived schemes (that is, fulfilling not only (1)–(2), but also (3) of Definition 4.11), we need to prove the following:

MAIN LEMMA 4.21. *Assume that $(\Sigma, \mathcal{A}, \vec{\tau})$ enjoys Property* (A). *We claim that*

(I) *The codimension-one submanifolds defined either by $m_i = 0$ ($1 \leq i \leq N$) or by $L = 0$ can not locally coincide with any degeneracy $\Phi_p$ or $\Psi_p$, $1 \leq p \leq n$;*

(II) *On $D_{p+1,n-1}$, the degeneracy $\Phi_p$ cannot locally coincide with the degeneracy $\Psi_n$ ($1 \leq p \leq n-1$) whenever $((\sigma_p, \ldots, \sigma_n), (a_p, \ldots, a_n))$ is $\Phi$-minimal from the left and $\Psi$-minimal from the right;*

(III) *On $D_{p+1,n-1}$, the degeneracy $\Psi_p$ cannot locally coincide with the degeneracy $\Psi_n$, whenever $((\sigma_p, \ldots, \sigma_n), (a_p, \ldots, a_n))$ is $\Psi$-minimal from both sides;*

(IV) *The degeneracies $m_i = 0$ ($1 \leq i \leq N$), $L = 0$ intersect $D_{p+1,n-1}$ ($1 \leq p \leq n-1$); moreover the degeneracies $\Phi_p$, $\Phi_n$, $\Psi_p$ and $\Psi_n$ intersect $D_{p+1,n-1}$ ($1 \leq p \leq n-1$) whenever $((\sigma_p, \ldots, \sigma_n), (a_p, \ldots, a_n))$ is $\Phi$-minimal from the left or from the right, or is $\Psi$-minimal from the left or from the right, respectively;*

(V) *If $(\Sigma, \mathcal{A}, \vec{\tau})$ is $\Psi$-minimal from the right, then there exists a (nondegenerate) orbit segment*

$$\omega = \{\tilde{q}_i^\alpha(\omega), v_i^\alpha(\omega), m_i(\omega), L(\omega) \,|\, i = 1, \ldots, N;\ \alpha = 0, 1, \ldots, n-1\}$$
$$\in \tilde{\Omega}(\sigma_1, \ldots, \sigma_{n-1}; a_1, \ldots, a_{n-1}; \tau_1, \ldots, \tau_{n-1}) = \tilde{\Omega}_{n-1}$$

*for which the discriminant*

$$\delta(\omega) := \left\langle v_{i_n}^{n-1}(\omega) - v_{j_n}^{n-1}(\omega); \tilde{q}_{i_n}^{n-1}(\omega) - \tilde{q}_{j_n}^{n-1}(\omega) - L \cdot a_n \right\rangle^2$$
$$- \left\| v_{i_n}^{n-1}(\omega) - v_{j_n}^{n-1}(\omega) \right\|^2 \cdot \left( \left\| \tilde{q}_{i_n}^{n-1}(\omega) - \tilde{q}_{j_n}^{n-1}(\omega) - L \cdot a_n \right\|^2 - 4r^2 \right)$$

*of (3.8) (with $k = n$) is equal to zero and $m_i(\omega) \neq 0$ ($i = 1, \ldots, N$), $L(\omega) \neq 0$;*

(VI) *If $1 \leq k < n$, $\sigma_k = \sigma_n$, $a_k = a_n$, $N \notin \sigma_n$, and*

(i) *for $k < j < n$ either $\sigma_j \cap \sigma_n = \emptyset$ or $N \in \sigma_j$,*

(ii) *there exists a $j_0$ with $k < j_0 < n$ and $\sigma_{j_0} \cap \sigma_n \neq \emptyset$,*

*then $|\mathbb{K}_n : \mathbb{K}_{n-1}| = 2$ or, equivalently, the polynomial (3.8), defining $\tau_n$, is irreducible over the field $\mathbb{K}_{n-1}$;*

(VII) *If $1 \leq k < n$, $\sigma_k = \sigma_n$, $a_k \neq a_n$, and for every $k < j < n$, $\sigma_j \cap \sigma_n = \emptyset$, then again $|\mathbb{K}_n : \mathbb{K}_{n-1}| = 2$.*

(VIII) *On $D_{p+1,n-1}$, the degeneracy $\Phi_p$ cannot locally coincide with the degeneracy $\Phi_n$ ($1 \leq p \leq n-1$) whenever $((\sigma_p,\ldots,\sigma_n),(a_p,\ldots,a_n))$ is $\Phi$-minimal from both sides.*

*Remark* 4.22. (a) Assertion (VI) of Main Lemma 4.21 is tailor-made to ensure that whenever the substitution $m_N = 0$ is carried out at every "critical" step with the possibility of selecting the wrong $\tau'_n$ (and thus violating Property (A)), we are not forced to do so, because at an irreducible $\tau_n$ either branch of the corresponding square root of the discriminant of (3.8) may be selected for continuing the dynamics. Quite similarly, assertion (VII) guarantees the same thing for the substitution $L = 0$.

(b) The careful reader has certainly observed that statement (VIII) is separated from the quite similar assertions (II)–(III) of the Main Lemma. The reason is the following: By a slight improvement of the inductive arguments verifying the transitions $N - 1 \to N$ when establishing (II)–(III) and (VIII), one can strengthen these claims by only requiring in each of them minimality from the right. The advantage of these stronger statements is that then (VIII) is not needed in the inductive proofs of (II)–(IV). As a result, (V)–(VII) can be proved completely without (VIII).

*Proof.* **1**. First we prove (I) and for the sake of brevity, only for the case $m_i = 0$. The handling of the submanifold $L = 0$ is completely analogous. Assume that there is a counterexample, and $p$ (the time of degeneracy) has the smallest possible value. Push forward the manifold $m_i = 0$ up to the collision $\sigma_{p-1}$ by the algebraic dynamics: the form of the equation $m_i = 0$ defining this manifold does not change. (This push-forward can be made without difficulty by the minimality of the index $p$.) In this way we can reduce the proof of (I) to the case $p = 1$. However, the statement is obviously true in that case. □

For later use we prove the following:

**Corollary** 4.23. *If for the combinatorial scheme $(\Sigma, \mathcal{A}, \vec{\tau})$ and for some $\omega_0 \in \tilde{\Omega}_{n-1}$, $\delta(\omega_0) = 0$, then there is an $\omega^* \in \tilde{\Omega}_{n-1}$ such that $\delta(\omega^*) = 0$, none of the masses $m_1(\omega^*),\ldots,m_N(\omega^*)$ is zero and, moreover, $L(\omega^*) \neq 0$.*

*Proof.* Indeed, the conditions imply that for a codimension-one algebraic subset $Z$ of $D\left(\tilde{\Omega}_{n-1}\right)$ it is true that for every $x \in Z$ there exists an $\omega \in \tilde{\Omega}_{n-1}$ such that $x(\omega) = x$ and $\delta(\omega) = 0$. By (I), $Z$ is not even locally identical with any submanifold $\{m_i = 0\}$ ($i = 1,\ldots,N$). Consequently, we can select a small perturbation $x^* \in Z$ of $x_0 = x(\omega_0)$ and an $\omega^* \in \tilde{\Omega}_{n-1}$ such that $x(\omega^*) = x^*$, $\delta(\omega^*) = 0$, and none of the masses $m_1(\omega^*),\ldots,m_N(\omega^*)$, nor $L(\omega^*)$, vanishes. □

**2**. The implication (II) and (III) and (IV) $\Longrightarrow$ (V) is obvious since, by the argument used to prove (I), it is no problem to avoid the additional degeneracies $m_i + m_j = 0$, $1 \leq i < j \leq N$.

**3**. We prove next that the existence of an $\omega \in \tilde{\Omega}_{n-1}$ with $\delta(\omega) = 0$ (guaranteed in (V)) implies $|\mathbb{K}_n : \mathbb{K}_{n-1}| = 2$. This will immediately prove the implications (V) $\Longrightarrow$ (VI) and (V) $\Longrightarrow$ (VII). (It is easy to check that the conditions in (VI) or (VII) actually imply the precondition of (V).)

Let, therefore, $\omega \in \tilde{\Omega}_{n-1}$ be an orbit segment with $\delta(\omega) = 0$. By relabeling the balls, we can assume that $\sigma_n = \{1, 2\}$. In a small neighborhood $U_0 \subset \tilde{\Omega}_{n-1}$ of $\omega$ the discriminant $\delta(\omega')$ from (V) ($\omega' \in U_0$) is an analytic function of $\omega'$ or, equivalently, $\delta(\omega')$ is an analytic, algebraic function $f(\vec{x})$ of the initial variables $\vec{x} = \vec{x}(\omega) \in \mathbb{C}^{(2\nu+1)N+1}$.

Consider now the "normal vector of collision"

$$w = w(\omega) = \tilde{q}_1^n(\omega) - \tilde{q}_2^n(\omega) - L \cdot a_n$$
$$= \tilde{q}_1^{n-1}(\omega) - \tilde{q}_2^{n-1}(\omega) - L \cdot a_n + \tau_n(\omega) \cdot \left[v_1^{n-1}(\omega) - v_2^{n-1}(\omega)\right],$$

where $\tau_n(\omega)$ is just the double root of (3.8) with $k = n$. Execute now the small perturbations $\omega(\varepsilon) \in U_0$, ($\varepsilon \in \mathbb{C}$, $|\varepsilon|$ is small) as follows:

$$\tilde{q}_1^{n-1}(\omega(\varepsilon)) = \tilde{q}_1^{n-1}(\omega) + \varepsilon \cdot w;$$
$$\tilde{q}_i^{n-1}(\omega(\varepsilon)) = \tilde{q}_i^{n-1}(\omega), \quad i = 2, \ldots, N;$$
$$v_i^{n-1}(\omega(\varepsilon)) = v_i^{n-1}(\omega), \quad i = 1, \ldots, N.$$

(Note: If $1 \in \sigma_{n-1}$, then the configuration perturbation defined above ought to be modified by an appropriate scalar multiple of the flow direction in order to ensure $\left\|\tilde{q}_{i_{n-1}}^{n-1} - \tilde{q}_{j_{n-1}}^{n-1} - L \cdot a_{n-1}\right\|^2 = 4r^2$.) Easy geometric consideration yields that

$$\frac{d}{d\varepsilon}\delta\left(\omega(\varepsilon)\right)\Big|_{\varepsilon=0} = \frac{d}{d\varepsilon}f\left(\vec{x}\left(\omega(\varepsilon)\right)\right)\Big|_{\varepsilon=0} = -8r^2 \left\|v_1^{n-1}(\omega) - v_2^{n-1}(\omega)\right\|^2 \neq 0.$$

That is, the holomorphic function $g(\varepsilon) := f\left(\vec{x}\left(\omega(\varepsilon)\right)\right)$ of $\varepsilon \in \mathbb{C}$ ($|\varepsilon|$ is small) has the properties $g(0) = 0$, $g'(0) \neq 0$.

Suppose now, indirectly, that the algebraic function $f(\vec{x}) \in \mathbb{K}_{n-1}$ is a square $(h(\vec{x}))^2$ in $\mathbb{K}_{n-1}$. Consider the function $z(\varepsilon) := h\left(\vec{x}(\omega(\varepsilon))\right)$ of the single complex variable $\varepsilon \in \mathbb{C}$, $|\varepsilon|$ small. It is clear that $z(\varepsilon)$ is a meromorphic function of $\varepsilon$, being a rational function of the variables $\{\tilde{q}_i^\alpha, v_i^\alpha, m_i, L \,|\, i = 1, \ldots, N; \alpha < n\}$, and thus being the ratio of two holomorphic functions. This is the point where we use the indirect assumption $h(\vec{x}) \in \mathbb{K}_{n-1}$. However, the square $[z(\varepsilon)]^2 = g(\varepsilon)$ of the meromorphic function $z(\varepsilon)$, obviously, cannot have the properties $g(0) = 0$, $g'(0) \neq 0$.

This proves that, indeed, $|\mathbb{K}_n : \mathbb{K}_{n-1}| = 2$. $\square$

So far we have proved that (II), (III) and (IV) together imply assertions (V)–(VII) of Main Lemma 4.21.

**4.** We prove now (II)–(IV) ((VIII) will be settled later in point **5**.) Let us first observe that (II)–(IV) (and similarly (VIII)) are consequences of somewhat weaker statements

(II*) On $D_{2,n-1}$, the degeneracy $\Phi_1$ cannot locally coincide with the degeneracy $\Psi_n$ whenever $((\sigma_1,\ldots,\sigma_n),(a_1,\ldots,a_n))$ is $\Phi$-minimal from the left and $\Psi$-minimal from the right;

(III*) On $D_{2,n-1}$, the degeneracy $\Psi_1$ cannot locally coincide with the degeneracy $\Psi_n$ whenever $((\sigma_1,\ldots,\sigma_n),(a_1,\ldots,a_n))$ is $\Psi$-minimal from both sides;

(IV*) The degeneracies $m_i = 0$ $(1 \le i \le N)$, $L = 0$ intersect $D_{2,n-1}$; moreover, the degeneracies $\Phi_1$, $\Phi_n$, $\Psi_1$ and $\Psi_n$ intersect $D_{2,n-1}$ whenever $((\sigma_1,\ldots,\sigma_n),(a_1,\ldots,a_n))$ is $\Phi$-minimal from the left or from the right, or is $\Psi$-minimal from the left or from the right, respectively;

(VIII*) On $D_{2,n-1}$, the degeneracy $\Phi_1$ cannot locally coincide with the degeneracy $\Phi_n$ whenever $((\sigma_1,\ldots,\sigma_n),(a_1,\ldots,a_n))$ is $\Phi$-minimal from both sides.

Indeed, we only have to shift the index (time) of the initial variables from 2 to $p+1$ and use the observation: If the triple $(\Sigma,\mathcal{A},\vec{\tau})$ has Property(A) and $(\Sigma',\mathcal{A}') = (\sigma_p,\ldots,\sigma_n; a_p,\ldots,a_n)$ $(1 < p \le n)$ is an end segment of $(\Sigma,\mathcal{A})$, then the $\Sigma'$-restriction of every $(\Sigma,\mathcal{A},\vec{\tau})$-orbit segment $\omega$ is automatically a $(\Sigma',\mathcal{A}',\tau')$-orbit segment with *some* $\tau'$, so that $(\Sigma',\mathcal{A}',\tau')$ *enjoys Property* (A). Observe, moreover, that if $\tau'_j$ $(p < j \le n)$ is reducible in the truncated dynamics $(\Sigma',\mathcal{A}',\vec{\tau}')$, then $\tau_j$ is certainly reducible in the original dynamics $(\Sigma,\mathcal{A},\vec{\tau})$. In other words, any possible additional irreducibility in the truncated dynamics permits a larger variety of orbit segments to study than just the restrictions of the orbit segments of $(\Sigma,\mathcal{A},\vec{\tau})$ to the collisions $\sigma_p,...,\sigma_n$. Therefore it is sufficient to establish the claimed transversalities for $(\Sigma',\mathcal{A}',\vec{\tau}')$, and this is exactly the content of (II*)–(III*) and (VIII*).

We will establish now (II*)–(IV*) (and (VIII*) later in point **5**) by using a double induction: we assume that the statements are true for any number $N' \le N$ of balls and arbitrary length $n' \le n$ such that the pair $(n',N')$ is not identical with $(n,N)$. More precisely, when establishing (IV*) for the pair $(n,N)$, we will only use (IV*) for pairs $(n',N)$ with $n' < n$. On the other hand, when demonstrating (II*)–(III*) (and (VIII*)), the validity of (II*)–(IV*) (and also of (VIII*)) for pairs $(n',N')$ with $n' < n$ and $N' < N$ will be exploited.

*The proof is formulated in such a way that we take first the outer induction on $n$, and then the statements for this fixed $n$ are shown by using an inner induction on $N$.*

Observe that, for an arbitrarily fixed number $N$ of balls, assertion (IV*) is certainly true for the value $n$ once (II*)–(IV*) and (VIII*) are known to hold for smaller values of $n$. Indeed, for the degeneracies $m_i = 0$ ($1 \leq i \leq N$) and $L = 0$, the statement is obvious. Our inductive hypothesis (i.e. the validity of (II*)–(III*) and (VIII*) for all smaller values of $n$) guarantees that the claimed intersections of the singularities with the working domain $D_{2,n-1}$ are nonempty, as long as the appropriate minimality conditions of (IV*) hold.

Let us turn to the inductive proofs of (II*)–(III*) (and to that of (VIII*) later in point **5**).

The base of the induction on $n$ is $n = 2$: This case is a straightforward inspection.

Assume now that $n \geq 3$, and that (II*)–(IV*) and (VIII*) have been proved for all smaller values of $n$. The base of the induction on $N$ is $N = 2$ or $3$ or $4$, depending on the form of $(\Sigma, \mathcal{A})$, when establishing (II*) and (III*). In both proofs the inductive steps will first be presented, and afterwards the base cases will be treated.

By re-labeling the balls, we can obtain that $\sigma_n = \{1, 2\}$, $\sigma_1 = \{i, j\}$, $1 \leq i < j \leq 4$, and if $j = 4$, then $i = 3$ ("densely packed" indices).

**4/a**: *Proof of* (II*) *by using the inductive assumption, i.e. by assuming that* (II*)–(IV*) *and* (VIII*) *are valid for pairs* $(n', N')$ *satisfying* $n' < n, N' < N$. We argue by contradiction. Suppose, therefore, that, on $D_{2,n-1}$, $\Phi_1$ locally coincides with the degeneracy $\Psi_n$ defined by the equation $\delta(\omega) = 0$, i.e. their intersection contains a nonempty, smooth, codimension-one submanifold of $D_{2,n-1}$. We apply our favorite:

*Globalization/substitution argument.* In the sense of Remark 4.18, our variables ($\in D_{2,n-1}$) will now be represented at reference time 2, and the (algebraic) dynamics exists (strictly) between the collisions 1 and $n$. By using two quadratic equations:

(1) the one expressing that $\sigma_1$ is a collision of two particular balls;
(2) the other defining our considered $\sigma_1$-degeneracy $\Phi_1$,

we can eliminate two variables, say $(\tilde{q}_i^2)_1$ and $(v_i^2)_1$, by expressing them as some algebraic functions $f_1(\vec{y})$, $f_2(\vec{y})$ of the other variables

$$\vec{y} = \left\{ (\tilde{q}_\alpha^2)_\beta, (v_\alpha^2)_\beta, m_\alpha, L \,\middle|\, \alpha = 1, \ldots, N;\ 1 \leq \beta \leq \nu \right\} \setminus \left\{ (\tilde{q}_i^2)_1, (v_i^2)_1 \right\}$$

measured at the collision $\sigma_2$. The kinetic variables of $\vec{y}$ measured at $\sigma_2$ (i.e. with superscript 2) are now serving as the initial, or reference variables in the upcoming globalization/substitution argument. As usual, the algebraic functions $f_1(\vec{y})$ and $f_2(\vec{y})$ only contain field operations and square roots. Then, by using the recursive feature of the $(\Sigma, \mathcal{A}, \vec{\tau})$-dynamics and the inductive assumption of 4.21, we can express the discriminant $\delta(\omega)$ as another (also essentially unique) algebraic function $g(\vec{y})$ of the variables $\vec{y}$.

We argue by contradiction. Assume, therefore, that the degeneracies $\Phi_1$ and $\Psi_n$ locally coincide in the working domain $D_{2,n-1}$. The local coincidence of these degeneracies precisely means:

(*) The algebraic expression $g(\vec{y})$ is identically zero on some nonempty, open set of vectors $\vec{y} \in \mathbb{C}^{(2\nu+1)N-1}$.

Then, according to Lemma 4.8, this provides, on $D_{2,n-1}$, the *global* containment: $\Phi_1 \subset \Psi_n$. This statement remains valid after the substitution $m_N = 0$, provided that it was valid in the original set-up. We remind the reader that the almost everywhere (nondegenerate) definability of the algebraic dynamics subjected to the side condition $m_N = 0$ easily follows from (I).

*Remark* 4.24. Since the substitutions $m_N = 0, L = 0$ play a basic role in the forthcoming arguments, it is very important to note here that, thanks to the inductive hypothesis of 4.21, we can assume that Property (A) is inherited by the algebraic billiard trajectory segment subjected to the side conditions $m_N = 0$ and $L = 0$. Indeed, as we mentioned in Remark 4.22, parts (VI)–(VII) are designed to ensure the inheritance of Property (A).

Let us return now to our proof of the inductive step of (II*). Assume, of course, the inductive hypothesis. We have, moreover, already assumed indirectly that $\Phi_1$ and $\Psi_n$ coincide locally in $D_{2,n-1}$. By the inductive assumption, $\Phi_1$ and $\Psi_n$ intersect $D_{2,n-1}$, whereas by the indirect assumption and by the globalization argument

$$(4.25) \qquad \Phi_1 \subset \Psi_n$$

holds on $D_{2,n-1}$. If, under the substitution $m_N = 0$ (say), the degeneracy $\Phi_1$ (and thus $\Psi_n$, too) intersects $D'$, the working domain (cf. Remark 4.20) for the derived scheme $(\Sigma', \mathcal{A}', \vec{\tau}')$ of $(\Sigma, \mathcal{A}, \vec{\tau})$ corresponds to $m_N = 0$, then the contradiction to the inductive hypothesis (II*) is evident. For every $k$ with $\sigma'_k \in \Sigma'$, denote by $\Phi'_k$ and $\Psi'_k$ the projections of $\Phi_k$ and $\Psi_k$ to the phase space of the $(N-1)$-ball system corresponding to taking $m_N = 0$ and cancelling the data of the $N^{\text{th}}$ ball (as many times before, the indexing in the derived schemes is again lacunary).

Consider now the case when, under the substitution $m_N = 0$, the degeneracy $\Phi'_1$ or $\Psi'_n$ does not intersect $D'$. The relation $\Phi'_1 \cap D' = \emptyset$ can only hold if the symbolic sequence $\Sigma'$ is not $\Phi$-minimal from the left and, similarly, $\Psi'_n \cap D' = \emptyset$ can only happen if the symbolic sequence $(\Sigma', \mathcal{A}')$ is not $\Psi$-minimal from the right.

In any case, there exists a subsegment $(\Sigma'', \mathcal{A}'', \vec{\tau}'')$ of $(\Sigma', \mathcal{A}', \vec{\tau}')$, which is minimal from both sides. Moreover, by denoting by $p''$ and $n''$ the indices of its initial and final elements (in the original indexing of $(\Sigma, \mathcal{A}, \vec{\tau})$), then, on one hand, $\sigma'_1 = \sigma'_{p''}$, and for every $\sigma'_j \in \Sigma'$ with $1 < j < p''$, $\sigma'_j \cap \sigma'_1 =$ either $\sigma'_1$ or $\emptyset$, and on the other hand, $(\sigma'_n, a'_n) = (\sigma'_{n''}, a'_{n''})$, and for every $\sigma'_j \in \Sigma'$ with $n'' < j < n$ either $(\sigma'_j, a'_j) = (\sigma'_n, a'_n)$ or $\sigma'_j \cap \sigma'_n = \emptyset$. Note that, by virtue of Remark 4.18/b, the just formulated conditions also imply $\Phi'_1 = \Phi'_{p''}$ and $\Psi'_{n''} = \Psi'_n$. Also, by inductive assumption (IV*), $\Phi'_{p''} \cap D_{p''+1,n''-1} \neq \emptyset$ and $\Psi'_{n''} \cap D_{p''+1,n''-1} \neq \emptyset$ and further, by (II*), they cannot coincide locally in $D_{p''+1,n''-1}$.

When completing the argument we will be exploiting the following lemma, which is a limiting variant of the observation of point (b) of Remark 4.18 and can be proved easily.

LEMMA 4.26. *Consider a symbolic collision sequence $(\Sigma, \mathcal{A}, \vec{\tau})$ which satisfies Property* (A) *and is $\Phi$-minimal from the left and $\Psi$-minimal from the right. Suppose we set $m_N = 0$ and denote by $(\Sigma', \mathcal{A}', \vec{\tau}')$ the corresponding derived scheme. Now, there are two statements (whose duals for $\Phi$-minimality from the right and $\Psi$-minimality from the left are also true):*

(1) *Assume that $\sigma'_1 = \sigma'_{p'}$ ($= \{1, 2\}$, say) $1 < p' \leq n$, and for every $\sigma'_j \in \Sigma'$ with $1 < j < p'$ one has $\sigma'_j \cap \sigma'_1 = \sigma'_1$ or $\emptyset$ (lacunary indices). Take a sequence $\omega_k \in \tilde{\Omega}(\Sigma; \mathcal{A}; \vec{\tau})$ such that*

  (i)   all kinetic data of $\omega_k$ measured right after $\sigma_1$ have a limit,

  (ii)  $\lim_{k \to \infty} m_N(\omega_k) = 0$,

  (iii) $\lim_{k \to \infty} m_i(\omega_k) \neq 0$, $\quad 1 \leq i \leq N-1$,

  (iv)  $\lim_{k \to \infty} m_i(\omega_k) + m_j(\omega_k) \neq 0$, $\quad 1 \leq i < j \leq N-1$,

  (v)   $\lim_{k \to \infty} \|v_1^1(\omega_k) - v_2^1(\omega_k)\|^2 = 0$.

  *Now $\lim_{k \to \infty} \|v_1^{p'}(\omega_k) - v_2^{p'}(\omega_k)\|^2 = 0$.*

(2) *Assume that $(\sigma'_{n'}, a'_{n'}) = (\sigma'_n, a'_n)$ ($\sigma'_n = \{1, 2\}$, say) $1 \leq n' < n$, and for every $\sigma'_j \in \Sigma'$ with $n' < j < n$ one either has that $(\sigma'_j, a'_j) = (\sigma'_n, a'_n)$ or $\sigma'_j \cap \sigma'_n = \emptyset$ is true.*

*Take a sequence $\omega_k \in \tilde{\Omega}(\Sigma; \mathcal{A}; \vec{\tau})$ such that* (i)–(iv) *hold (with $\sigma_1$ replaced by $\sigma_n$ in* (i)*), and further*

  (v)$'$  $\lim_{k \to \infty} \langle v_1^n(\omega_k) - v_2^n(\omega_k); \tilde{q}_1^n(\omega_k) - \tilde{q}_2^n(\omega_k) - L \cdot a_n \rangle = 0$

*is also true.* Now

$$\lim_{k\to\infty} \left\langle v_1^{n'}(\omega_k) - v_2^{n'}(\omega_k);\ \tilde{q}_1^{n'}(\omega_k) - \tilde{q}_2^{n'}(\omega_k) - L \cdot a_{n'} \right\rangle = 0.$$

Take now convergent sequences of points $\vec{x}_k \in \Phi_1 \cap D_{2,n-1} \subset \Psi_n \cap D_{2,n-1}$ such that $m_N(\vec{x}_k) \to 0$ ($k \to \infty$). Then, according to Lemma 4.26, the projections $\vec{x}'$ of the points $\vec{x}$ to the $(N-1)$-ball phase space of the system necessarily fill up a nonempty open piece $U'$ of $\Phi'_{p''} = \Phi'_1$ on one hand, which necessarily also belongs to $\Psi'_{n''} = \Psi'_n$ on the other hand. Consequently, $\Phi'_{p''}$ coincides locally with $\Psi'_{n''}$, a contradiction. The conclusion is that if (II*) holds for the $(N-1)$-ball system, then it does also for the $N$-ball system. This completes the proof of the inductive step $(N-1) \to N$ for statement (II*).

As mentioned before, the base of the induction on $N$ is $N = 2$ or 3 or 4 depending on the form of $(\Sigma, \mathcal{A})$, and we treat two cases separately. We remind the reader of our notational convention to use densely packed indices.

*Case 4/a/1.* $\sigma_1 \neq \sigma_n$. In this case the base of the induction is $N = 3$ or 4.

The globalization argument provides that the indirect assumption that $\Phi_1$ and $\Psi_n$ coincide locally leads to the global containment $\Phi_1 \subset \Psi_n$ on $D_{2,n-1}$. Then this ought to hold under the side conditions $m_3 = m_4 = 0$, too. But this is impossible for the following reason: In the case $m_3 = m_4 = 0$ the motion of the particles $\{3, 4\}$ has absolutely no effect on the time evolution of the $\{1, 2\}$ subsystem. However, the equation defining $\Phi_1$ is obviously sensitive to the kinetic data (position and velocity) of the third ball, while the equation $\delta(\omega) = 0$ is written fully in terms of the subsystem $\{1, 2\}$. The obtained contradiction proves the base case of the induction for (II*) in the case $\sigma_1 \neq \sigma_n$.

*Case 4/a/2.* $\sigma_1 = \sigma_n$. Now the base of the induction is $N = 2$.

Similarly to the previous case, the globalization argument implies the containment $\Phi_1 \subset \Psi_n$. This is an assertion completely in terms of the now independently evolving two-particle system $\{1, 2\}$. We will show that in such a complexified, two-particle algebraic system the tangential degeneracy $\delta(\omega) = 0$ of the collision $\sigma_n$ cannot completely contain the degeneracy $\|v_1^2 - v_2^2\|^2 = 0$. (Now the $\sigma$'s only denote the collisions of the two-ball system $\{1, 2\}$.) The proof will be a simple induction on the positive integer $n$. For $n = 2$ the statement is obviously true. Assume now that $n > 2$, and the assertion has been proven for all smaller values of $n$. According to our inductive hypothesis, for a typical element $\omega$ of the $\sigma_n$-degeneracy $\delta(\omega) = 0$ we have that $\left\|v_1^3(\omega) - v_2^3(\omega)\right\|^2 \neq 0$. On the other hand, the collision $\sigma_3$ is known to preserve the quantity $\|v_1 - v_2\|^2$ (the velocity reflection is complex orthogonal), and, therefore, we get that $\left\|v_1^2(\omega) - v_2^2(\omega)\right\|^2 \neq 0$, as well.

The obtained contradiction finishes the discussion of 4/a.

**4/b**. *Proof of* (III$^*$) *by using the inductive hypothesis.* Assume now, besides the inductive hypothesis, the precondition of (III$^*$).

The proof of the inductive step $(N-1) \to N$ is absolutely analogous to that given in part 4/a: one has just to substitute $\Phi_1$ with $\Psi_1$. The verification of the basis of the induction will be split into three parts depending on the form of $(\Sigma, \mathcal{A})$. In all cases we again assume indirectly that $\Psi_1$ coincides locally with $\Psi_n$ in $D_{2,n-1}$. By the globalization/substitution argument from the proof of (II$^*$) we then obtain $\Psi_1 \subset \Psi_n$ to hold globally in $D_{2,n-1}$ for the base cases, too. What we show is that this inclusion leads to a contradiction in all cases.

*Case 4/b/1.* $\sigma_1 \ne \sigma_n$. In this case the base of the induction is $N = 3$ or 4, and the impossibility of the containment $\Psi_1 \subset \Psi_n$ for $N = 3$ or 4 can be obtained in the same way as in 4/a/1.

*Case 4/b/2.* $\sigma_1 = \sigma_n$ and there is a $k$, $1 < k < n$, such that $\sigma_k \cap \sigma_n \ne \emptyset$ and $\sigma_k \ne \sigma_n$.

Now the base case will be $N = 3$ and, as explained above, we can assume the global containment $\Psi_1 \subset \Psi_n$ on $D_{2,n-1}$. By re-labeling the particles, we can assume that $\sigma_1 = \sigma_n = \{1,2\}$ and $\sigma_k = \{2,3\}$.

The global containment $\Psi_1 \subset \Psi_n$ will naturally be preserved under the substitution $L = 0$, too. This substitution raises the same questions as the substitution $m_N = 0$, and they can be handled analogously. For instance, it is easily possible that though $(\Sigma, \mathcal{A}, \vec{\tau})$ is minimal from both sides, the scheme $(\Sigma', \mathcal{A}', \vec{\tau}')$ derived from it via the substitution $L = 0$ is not minimal. Our way of taking care of this possibility in the case of the substitution $m_N = 0$ can also be repeated here and, therefore, by our combined indirect and inductive arguments, we can again assume that the containment $\Psi_1 \subset \Psi_n$ now holds even in the 3-particle system $\{1,2,3\}$ (having a time evolution with Property (A)). On the other hand, the fact that this 3-particle system $\{1,2,3\}$ also has the property $L = 0$ amounts to $a_l = 0$ for $l = 1, \ldots, n$. Obviously, we can assume that $k = 2$, i.e., $\sigma_2 = \{2,3\}$. Case 4/b/2 will be taken care of by the following:

LEMMA 4.27. *Suppose that the triple $(\Sigma = (\sigma_0, \ldots, \sigma_n), 0, \vec{\tau})$ enjoys Property* (A), $N = 3$, $\sigma_0 = \sigma_n = \{1,2\}$ *and, say,* $\sigma_1 = \{2,3\}$. *Then $\Psi_0$ cannot be (globally) a subset of $\Psi_n$ and, therefore, the tangential degeneracies at $\sigma_0$ and $\sigma_n$ (i.e. $\Psi_0$ and $\Psi_n$) cannot even locally coincide in $D_{1,n-1}$.*

*Proof.* In order to simplify the notation we assume that $r = 1/2$, and $\sigma_{n-1} \ne \{1,2\}$. Denote by $\Sigma^* = (\sigma_1^*, \ldots, \sigma_s^*)$ the symbolic sequence that can be obtained from our original sequence $(\sigma_0, \ldots, \sigma_n)$ by discarding all symbols $\{1,2\}$. The truncated sequence $\Sigma^*$ consists of a pure block of $n(1)$ $\{2,3\}$-collisions first, then a pure block of $n(2)$ $\{1,3\}$-collisions etc., up to the closing pure block of $n(k)$ collisions of type $\{(3+(-1)^{k+1})/2, 3\}$. Consider first the following, essentially two-dimensional and completely real, limiting set-up: $m_1 = m_2 = 1$, $m_3 = 0$, $\tilde{q}_1^0 = 0$, $\tilde{q}_2^0 = e_1$, $\tilde{q}_3^0 = -e_2$, $v_1^0 = v_2^0 = 0$, $v_3^0 = (1+\varepsilon^2)^{-1/2} \cdot (e_1 + \varepsilon e_2)$, where $e_1 = (1, 0, \ldots, 0) \in \mathbb{R}^\nu$, $e_2 = (0, 1, 0, \ldots, 0) \in \mathbb{R}^\nu$, and $\varepsilon > 0$, $\varepsilon \to 0$ is a positive parameter of perturbation. By using Property (A) of $(\Sigma, 0, \vec{\tau})$ and the entire inductive hypothesis of 4.21, we can now construct a $(\Sigma, 0, \vec{\tau})$ orbit segment

$$\omega = \{\tilde{q}_i^\alpha, v_i^\alpha \mid i = 1, 2, 3;\ \alpha = 0, \ldots, n\}$$

by extending the fixed initial values $m_1 = m_2 = 1$, $m_3 = 0$, $\tilde{q}_1^0 = 0$, $\tilde{q}_2^0 = e_1$, $\tilde{q}_3^0 = -e_2$, $v_1^0 = v_2^0 = 0$, $v_3^0 = (1+\varepsilon^2)^{-1/2} \cdot (e_1 + \varepsilon e_2)$ in such a way that

(i) At every collision $\sigma_j = \{2, 3\}$ $(0 < j < n)$ the smaller root $\tau_j$ is always selected out of the two distinct, real roots;

(ii) At every collision $\sigma_j = \{1, 3\}$ $(0 < j < n)$ the greater root $\tau_j$ is always selected out of the two distinct, real roots.

We note here that, obviously, the above constructed orbit segment $\omega$ is a highly degenerate one because of the degenerate $\{1, 2\}$-collisions. However, we will transform $\omega$ into a nondegenerate orbit segment (except one tangential collision at $\sigma_0$) by the upcoming, final perturbation.

For $j = 1, \ldots, k$ we denote by $\nu(j)$ the total number of collisions $\sigma_i$ (i.e. together with the $\{1, 2\}$-collisions) from $\sigma_1$ up to the closing collision $\sigma_h^*$ ($h = n(1) + \cdots + n(j)$) of the $j^{\text{th}}$ $\Sigma^*$ block. We note that, due to the selection of the initial kinetic data, the $\{1, 2\}$-collisions have no effect on the motion of the third disk.

An easy geometric calculation shows that the asymptotic formulae

$$\tilde{q}_3^{\nu(1)}(\varepsilon) = e_1 - e_2 - (2n(1) - 1) \cdot (2\varepsilon)^{1/2} \cdot e_1 + o\left(\varepsilon^{1/2}\right),$$

$$v_3^{\nu(1)}(\varepsilon) = e_1 - 2n(1)(2\varepsilon)^{1/2} \cdot e_2 + o\left(\varepsilon^{1/2}\right)$$

hold. By similar arguments we obtain the asymptotics

$$\tilde{q}_3^{\nu(2)}(\varepsilon) = -e_2 + (2n(2) - 1) \cdot \left[4n(1) \cdot (2\varepsilon)^{1/2}\right]^{1/2} \cdot e_1 + o\left(\varepsilon^{1/4}\right),$$

$$v_3^{\nu(2)}(\varepsilon) = e_1 + 2n(2) \cdot \left[4n(1) \cdot (2\varepsilon)^{1/2}\right]^{1/2} \cdot e_2 + o\left(\varepsilon^{1/4}\right).$$

Then a simple induction yields that, in general, the following asymptotic formulae hold for $j = 0, 1, \ldots, k$:

(4.28)
$$\tilde{q}_3^{\nu(j)}(\varepsilon) = \frac{1 + (-1)^{j+1}}{2} \cdot e_1 - e_2 + (-1)^j b_j \varepsilon^{2^{-j}} \cdot e_1 + o\left(\varepsilon^{2^{-j}}\right),$$

$$v_3^{\nu(j)}(\varepsilon) = e_1 + (-1)^j c_j \varepsilon^{2^{-j}} \cdot e_2 + o\left(\varepsilon^{2^{-j}}\right),$$

where the positive constants $b_j$ and $c_j$ can be computed by the recursion $b_0 = 0$, $c_0 = 1$, $b_j = (2n(j) - 1)\sqrt{2c_{j-1}}$, $c_j = 2n(j)\sqrt{2c_{j-1}}$, $j = 1, \ldots, k$.

Let us now focus our attention on the *almost vertical* velocity changes

$$\Delta V_j(\varepsilon) = v_3^{\nu(j)}(\varepsilon) - v_3^{\nu(j-1)}(\varepsilon)$$

and their horizontal components $\langle \Delta V_j(\varepsilon); e_1 \rangle$, $j = 1, \ldots, k$. The asymptotic formulae (4.28) for the velocities immediately yield the expressions

(4.29)
$$\Delta V_j(\varepsilon) = (-1)^j c_j \varepsilon^{2-j} \cdot e_2 + o\left(\varepsilon^{2-j}\right),$$

$$\langle \Delta V_j(\varepsilon); e_1 \rangle = -\frac{1}{2} c_j^2 \varepsilon^{2^{1-j}} + o\left(\varepsilon^{2^{1-j}}\right),$$

$j = 1, \ldots, k$. Here we used the fact that the norm $||v_3^i(\varepsilon)||$ is always 1.

From the asymptotic formulae (4.29) we conclude that the order of magnitude of the horizontal velocity change is the highest one in the last pure block of collisions $B_k = \left(\sigma^*_{s-n(k)+1}, \ldots, \sigma^*_s\right)$ (actually, $s - n(k) = \sum_{j=1}^{k-1} n(j)$). These collisions take place between the third disk and the disk with the label $i(k) = \left(3 + (-1)^{k+1}\right)/2$. However, after a perturbation (to be carried out a bit later) $m_3 \approx 0$, $m_i \in \mathbb{R}_+$ ($i = 1, 2, 3$), only a certain part of the horizontal momentum given by the third disk to the $i(k)$-disk (during the collisions of $B_k$) will be retained by the disk $i(k)$: the rest will be transported to the other heavy disk with label $3 - i(k)$ by the almost tangential collisions $\{1, 2\}$ (taking place in between the collisions of $B_k$) with almost horizontal normals of impact. This justifies the following definition:

$$H := \{q \in \mathbb{Z} |\ 1 \leq q \leq n(k) \text{ and the number of the}$$
$$\{1,2\}\text{-collisions between } \sigma^*_{s-n(k)+q} \text{ and } \sigma^*_s \text{ is even}\},$$
$$\overline{H} := \{1, 2, \ldots, n(k)\} \setminus H.$$

By using the method of the proof of the second line of (4.29) we see that the order of magnitude of the horizontal velocity change at the $q^{\text{th}}$ collision of $B_k$ ($1 \leq q \leq n(k)$) is

$$4(2q - 1)c_{k-1} \varepsilon^{2^{1-k}} + o\left(\varepsilon^{2^{1-k}}\right).$$

Note that the cumulative horizontal velocity change caused by the collisions

$$\sigma^*_{s-n(k)+1}, \ldots, \sigma^*_{s-n(k)+q}$$

has the asymptotic formula $4q^2 c_{k-1} \varepsilon^{2^{1-k}}$; note how one can deduce the formula $c_k = 2n(k)\sqrt{2c_{k-1}}$.

There are now two possibilities:

*Case* A. $\sum_{q \in H}(2q-1) \neq \frac{1}{2}[n(k)]^2$. (Note that $\sum_{q \in H \cup \overline{H}}(2q-1) = [n(k)]^2$.)
In this case the small perturbation $m_3(\varepsilon)$,

(4.30) $\quad 0 < m_3(\varepsilon) << \varepsilon$, is very small, depending on $\varepsilon$,
$$v_1^0(\varepsilon) = -v_2^0(\varepsilon) = \eta(\varepsilon) \cdot e_2, \ 0 < \eta(\varepsilon) = \eta(m_3(\varepsilon)) << m_3(\varepsilon)$$
and $\eta(\varepsilon)$ is very small compared to $m_3(\varepsilon)$

(still $m_1 = m_2 = 1$) serves our purpose: the horizontal component of the relative velocity right before the collision $\sigma_n = \{1, 2\}$ will be of the order of magnitude

$$(4.31) \quad 4m_3 c_{k-1} \left[ \sum_{q \in H}(2q-1) - \sum_{q \in \overline{H}}(2q-1) \right] \varepsilon^{2^{1-k}} \neq 0.$$

We note that the asymptotic formulae (4.28–4.29) will not be destroyed as long as the quantities $m_3(\varepsilon)$ and $\eta(\varepsilon) << m_3(\varepsilon)$ are small enough depending on $\varepsilon$. Also note that the collision $\sigma_n$ (as well as all other collisions $\{1, 2\}$ amongst $\sigma_0, \ldots, \sigma_n$) will be an asymptotically tangential collision with asymptotically horizontal normal of impact. These facts immediately follow from the asymptotic formulae (4.29). Hence, in Case A we managed to perturb out of the degeneracy $\Psi_n$ from within $\Psi_0$.

*Case* B. $\sum_{q \in H}(2q-1) = \sum_{q \in \overline{H}}(2q-1)$.
In this case the equation $m_1 = m_2$ does not suffice: we need to use the possibility of selecting different "heavy" masses. First we select and fix the positive values $m_1 \neq m_2$, then we carry out the perturbation (4.30) and compute the orders of magnitude of the horizontal components $V_1^h$, $V_2^h$ of the velocities of the first and second disks measured right after the collision $\sigma_{n-1} = \sigma_s^*$. By applying the collision laws for horizontal "head-on" collisions with masses $m_1$ and $m_2$ we easily obtain the asymptotics

(4.32)
$$V_1^h = 2c_{k-1}[n(k)]^2 \left[ \frac{2m_1 m_3}{m_2(m_1+m_2)} + \frac{m_3}{m_1} - \frac{m_3}{m_2} \right] \cdot \varepsilon^{2^{1-k}} + o\left(\varepsilon^{2^{1-k}}\right),$$
$$V_2^h = 4c_{k-1}[n(k)]^2 \frac{m_1 m_3}{m_2(m_1+m_2)} \cdot \varepsilon^{2^{1-k}} + o\left(\varepsilon^{2^{1-k}}\right).$$

These asymptotics are different, thanks to the inequality $m_1 \neq m_2$. Therefore, we again managed to perturb out of the degeneracy $\Psi_n$ from within $\Psi_0$ by using the fact that all collisions $\{1, 2\}$ amongst $\sigma_0, \ldots, \sigma_n$ are asymptotically tangential collisions with asymptotically horizontal normals of impact.

Finally, by using the method of obtaining different asymptotics (4.31)–(4.32) for the horizontal components of $v_1^{n-1}$ and $v_2^{n-1}$, we can also achieve different asymptotics for the horizontal components of $v_1$ and $v_2$ right before every collision $\sigma_i = \{1, 2\}$ $(i > 0)$ by appropriately selecting the perturbation. In that way all collisions $\sigma_1, \ldots, \sigma_n$ become nonsingular, and all kinetic data of the constructed orbit segment will be real.

This finishes the proof of Lemma 4.27. □

Our last outstanding task in proving (I)–(VII) is the handling of

*Case 4/b/3.* $\sigma_1 = \sigma_n$ $(= \{1, 2\})$, *and for every* $k \in \{1, \ldots, n\}$ *it is true that either* $\sigma_k = \sigma_n$ *or* $\sigma_k \cap \sigma_n = \emptyset$.

In this case not all adjustment vectors $a_k$ with $\sigma_k = \sigma_n$ $(p \le k \le n)$ are the same.

The base case now is $N = 2$, and according to the indirect assumption we have a global containment of the tangential degeneracies: $\Psi_1 \subset \Psi_n$ in the 2-ball system $\{1, 2\}$. The discussion of Case 4/b/3 will be finished as soon as we prove the following:

LEMMA 4.33. *Suppose that* $N = 2$, *the triple* $(\Sigma, \mathcal{A}, \vec{\tau})$ *enjoys Property* (A), *and not all adjustment vectors* $a_1, \ldots, a_n$ *of* $\mathcal{A}$ *are the same. Then* $\Psi_1$ *and* $\Psi_n$ *cannot locally coincide in* $D_{2,n-1}$.

*Proof.*
$1^\circ$. We will construct a one-parameter family
$$\omega(s) \in \tilde{\Omega}(\Sigma, \mathcal{A}, \vec{\tau}) = \tilde{\Omega}, \qquad (s \in I \subset \mathbb{R}),$$
($I$ is an open and bounded interval of $\mathbb{R}$) of $(\Sigma, \mathcal{A}, \vec{\tau})$-orbit segments in such a way that this bunch of trajectories will explicitly show that $\Psi_1$ and $\Psi_n$ are not the same. In fact, we will see that $\Psi_1$ and $\Psi_n$ intersect each other transversally at a large collection of points of intersection.

$2^\circ$. Throughout the construction of the family $\{\omega(s) | s \in I\}$ we shall follow these guiding principles:

(a) The dimension $\nu$ in this construction will be two: The transversality of $\Psi_1$ and $\Psi_n$ in a two-dimensional section $\tilde{q}_i^t(s) \in \mathbb{C}^2 \subset \mathbb{C}^\nu$, $v_i^t(s) \in \mathbb{C}^2 \subset \mathbb{C}^\nu$ clearly implies the transversality in the general set-up $\nu \ge 2$.

(b) All data $\tilde{q}_i^t(s) \in \mathbb{R}^2$, $v_i^t(s) \in \mathbb{R}^2$, $m_1(s) = m_2(s) = 1$, $r = 1/2$, $L(s) = L \gg 1$ will be real, the side length $L$ of the 2-torus being a very large, fixed positive constant $L$.

(c) Due to the translation invariance of the system, we will only deal with the relative motion
$$q^t(s) := \tilde{q}_1^t(s) - \tilde{q}_2^t(s) \in \mathbb{R}^2,$$
$$v^t(s) := v_1^t(s) - v_2^t(s) \in \mathbb{R}^2, \quad \|v^t(s)\|^2 = 1, \ (s \in I).$$

(d) In virtue of the condition $\sum_{i=1}^{n} ||a_i - a_1||^2 > 0$, by cutting off a suitable starting segment of $\Sigma$ (if necessary), we can assume that $0 = a_1 \neq a_2$.

(e) The construction of the family $\omega(s)$ will be such that for every $k = 2, 3, \ldots, n$ the polynomial (3.8) has two distinct, nonnegative real roots. We will always take the greater one as $\tau_k$. The inductive assumption on the irreducibility of these polynomials (3.8) over the field $\mathbb{K}_{k-1}$ (which are, in fact, irreducible, unless they have to be reducible, because of the repetitive adjustment vectors.) Part of the inductive hypothesis of 4.21 and Property (A) of $(\Sigma, \mathcal{A}, \vec{\tau})$ together ensure that the time evolution of $\omega(s)$ will follow the algebraic pattern encoded in $(\Sigma, \mathcal{A}, \vec{\tau})$.

(f)
$$q^0(s) := (\cos(\phi_0 + s), \sin(\phi_0 + s)) + (1 - s) \cdot (-\sin(\phi_0 + s), \cos(\phi_0 + s)),$$
$$v^0(s) := (-\sin(\phi_0 + s), \cos(\phi_0 + s)),$$

where the ray $\{q^0(0) + tv^0(0) \mid t \geq 0\}$ passes through the point $L \cdot a_2$.

The above constructed family $\{q^0(s), v^0(s) \mid |s| < \varepsilon_0, s \in \mathbb{R}\}$ is a so-called convex family of rays, or convex wave front; see, for example, [B(1979)] or [W(1986)]. This means that

$$\left\langle \frac{d}{ds} q^0(s); v^0(s) \right\rangle = 0, \quad \left\langle \frac{d}{ds} q^0(s); \frac{d}{ds} v^0(s) \right\rangle > 0.$$

The time evolution of such diverging families of trajectories is discussed in the literature; see [W(1986)]. It is shown (see Section 3, especially Theorem 3) that our two-dimensional billiard with *inner reflections* at the circles

$$C_a := \{q \in \mathbb{R}^2 \mid d(q, L \cdot a) = 1\} \quad (a \in \mathbb{Z}^2)$$

has the so-called convex scattering property; i.e. every such convex family of rays (with not too big curvature, see [W(1986)]) focuses after any inner reflection at $C_a$ not later than when it reaches the midpoint of the chord of $C_a$ along which it is traveling after the considered reflection.

$3^o$. Plainly, the trajectory $\omega(s)$ intersects the circle $C_0$ tangentially in negative time $\tau_1(s) = s - 1$; that is, the convex family of rays $\{q^0(s), v^0(s) \mid |s| < \varepsilon_0\}$ is the "synchronized" family of trajectories belonging to $\Psi_1$. We recall that $r = 1/2$ and $a_1 = 0 \in \mathbb{Z}^2$.

$4^o$. The *inner reflections* at the circles $C_{a_2}, \ldots, C_{a_n}$ (yet to be constructed) will be in good harmony with the principle $2^0(e)$ above.

$5^o$. If $L$ is large enough, then there is a small, open subinterval $I_2 \subset (-\varepsilon_0, \varepsilon_0)$ (even the closure of $I_2$ will be a subset of $(-\varepsilon_0, \varepsilon_0)$) such that the rays $\{q^0(s) + tv^0(s) \mid t \geq 0\}$ intersect the circle $C_{a_2}$ in two distinct points for $s \in I_2$, while these rays are tangent to $C_{a_2}$ for $s \in \partial I_2$. For $s \in I_2$ we elongate the above ray up to the second point of intersection with the circle $C_{a_2}$, and this way obtain the orbit segment up to the first collision with that circle.

6°. Let $a_2 = a_3 = \cdots = a_{k(2)}$ $(2 \leq k(2) \leq n)$, where either $a_{k(2)+1} \neq a_{k(2)}$ or $k(2) = n$. Suppose that $k(2)$ is still less than $n$. By inspection of the geometry of the sequence of inner reflections $\sigma_2, \ldots, \sigma_{k(2)}$ at the circle $C_{a_2}$ it is easy to see that the outgoing velocities $v^{k(2)}(s)$ $(s \in I_2)$ sweep out the entire space $\mathbb{S}^1$ of velocity directions. Therefore, there is an even smaller interval $I_3$, $\text{Cl}(I_3) \subset I_2$, such that the rays $\{q^{k(2)}(s) + tv^{k(2)}(s) \mid t \geq 0\}$ intersect the circle $C_{a_{k(2)+1}}$ in two distinct points for $s \in I_3$, while these rays are tangent to $C_{a_{k(2)+1}}$ for $s \in \partial I_3$. (Here $q^{k(2)}(s)$ and $v^{k(2)}(s)$ denote the position and the outgoing velocity at the reflection $\sigma_{k(2)}$ taking place at time $\sum_{i=1}^{k(2)} \tau_i$.)

For $s \in I_3$ we again elongate the above ray up to the second point of intersection with the circle $C_{a_{k(2)+1}}$, and this way obtain the orbit segment up to the first collision with that circle.

7°. Now we keep going with the iterative construction of the orbit segments $\omega(s) \in \tilde{\Omega}(\Sigma, \mathcal{A}, \vec{\tau})$ described above, by also constructing the nested sequence of open intervals $I_2 \supset I_3 \supset \cdots \supset I_l := I$, where $k(l) = n$. The two orbit segments $\{\omega(s) \mid s \in \partial I\}$ are now obviously in the transversal intersection of $\Psi_1$ and $\Psi_n$. Thus the lemma is proved. □

So far we have settled statements (I)–(VII) of the Main Lemma.

**5.** Our only remaining task in the proof of the Main Lemma is to verify statement (VIII). Similarly to (II)–(III) it is again sufficient to demonstrate the slightly weaker assertion (VIII$^*$) formulated in point **4.**

Assume now, besides the inductive hypothesis, the precondition of (VIII$^*$). We again use double induction in the same way it was explained before 4/a. Also, the proof of the inductive step $(N-1) \to N$ is absolutely analogous to that given in part 4/a: one has just to substitute $\Psi_n$ with $\Phi_n$. The verification of the basis of the induction will be split into two parts depending on the form of $\Sigma$. In all cases we again assume indirectly that $\Phi_1$ coincides locally with $\Phi_n$. By the globalization/substitution argument from the proof of (II$^*$) we then obtain $\Phi_1 \subset \Phi_n$ to hold globally for the base cases, too. What we show is that this inclusion leads to a contradiction in both cases.

*Case* 5/1. $\sigma_1 \neq \sigma_n$. In this case the base of the induction is $N = 3$ or $4$, and the impossibility of the containment $\Phi_1 \subset \Phi_n$ for $N = 3$ or $4$ can be obtained in the same way as in the proof of (II$^*$).

*Case* 5/2. $\sigma_1 = \sigma_n$ and there is a $k$, $1 < k < n$, such that $\sigma_k \cap \sigma_n \neq \emptyset$ and $\sigma_k \neq \sigma_n$.

Now the base case will be $N = 3$, and as explained above we can derive the global containment $\Phi_1 \subset \Phi_n$. By re-labeling the particles, we can assume that $\sigma_1 = \sigma_n = \{1, 2\}$ and $\sigma_k = \{2, 3\}$.

The global containment $\Phi_1 \subset \Phi_n$ will naturally be preserved under the substitution $L = 0$, too. The assumed containment $\Phi_1 \subset \Phi_n$ now implies that even in the 3-particle system $\{1, 2, 3\}$ (having a time evolution with Property (A)) one has a global containment of the $\Phi$-degeneracies at $\sigma_1$ and $\sigma_n$. The fact that this 3-particle system $\{1, 2, 3\}$ also has the property $L = 0$ amounts to $a_l = 0$ for $l = 1, \ldots, n$. Obviously, we can assume that $k = 2$, i.e. $\sigma_2 = \{2, 3\}$. Case 5/2 will be taken care of by the following analogue of Lemma 4.27:

LEMMA 4.34. *Suppose that the triple $(\Sigma = (\sigma_0, \ldots, \sigma_n), 0, \vec{\tau})$ enjoys Property (A), $N = 3$, $\sigma_0 = \sigma_n = \{1, 2\}$ and, say, $\sigma_1 = \{2, 3\}$. Then $\Phi_0$ cannot be (globally) a subset of $\Phi_n$ and, therefore, the $\Phi$-degeneracies at $\sigma_0$ and $\sigma_n$ (i.e. $\Phi_0$ and $\Phi_n$) cannot even locally coincide in $D_{1,n-1}$.*

*Proof.* The construction used in the proof of Lemma 4.27 will be used in this construction with $\eta(\varepsilon) = 0$. Then, because $v_1^0(\varepsilon) = v_2^0(\varepsilon) = 0$, the initial perturbation lies in $\Phi_1$. At the same time, in view of (4.31,4.32), for the horizontal components $V_1^h, V_2^h$ of the velocities of the first and second discs measured right after the collision $\sigma_{n-1} = \sigma_s^*$ we have $\|V_1^h - V_2^h\|^2 \neq 0$ showing that the perturbation necessarily leads out of $\Phi_n$. Thus the globalization/substitution argument provides the claimed transversality of $\Phi_1$ and of $\Phi_n$, the statement of Lemma 4.34. □

This finishes the proof of Main Lemma 4.21, as well. □

COROLLARY 4.35. *For every triple $(\Sigma, \mathcal{A}, \vec{\tau})$ with Property (A) and $N \geq 3$, there exists at least one derived scheme $(\Sigma', \mathcal{A}', \vec{\tau}')$.*

COROLLARY 4.36. *For combinatorial schemes satisfying Property (A), the quadratic polynomial $b_l x^2 + c_l x + d_l$ is reducible over the field $\mathbb{K}_{l-1}$ (i.e. $\mathbb{K}_{l-1} = \mathbb{K}_l$) only in the case described above in 3.31.*

*Proof.* Unless we are considering the case described in Remark 3.31, statement (V) of Main Lemma 4.21 gives the existence of an $\omega \in \tilde{\Omega}_{n-1}$ with $\delta(\omega) = 0$. But then part **3** of the proof of the main lemma provides $|\mathbb{K}_n : \mathbb{K}_{n-1}| = 2$. □

## Proof of the key lemma. Final induction

The next lemma establishes the link between the orbit segments of $(\Sigma, \mathcal{A}, \vec{\tau})$ and $(\Sigma', \mathcal{A}', \vec{\tau}')$. It is a direct consequence of Remark 4.13 (actually, of the fact that the reducibility $\mathbb{K}_{l-1} = \mathbb{K}_l$ is inherited by the $(\Sigma', \mathcal{A}', \vec{\tau}')$-dynamics during the substitution $m_N = 0$) and of part (2) of the definition of the derived schemes.

LEMMA 4.37. *Fix a combinatorial scheme $(\Sigma, \mathcal{A}, \vec{\tau})$ satisfying Property (A) and an arbitrary scheme $(\Sigma', \mathcal{A}', \vec{\tau}')$ derived from it. Then, for almost every choice of the vector of data $(\tilde{q}_N^0, v_N^0) \in \mathbb{C}^\nu \times \mathbb{C}^\nu$,*

(1) *The subset $\{(x', \tilde{q}_N^0, v_N^0, m_N = 0) | \, x' \in \mathbb{C}^h\} \cap D(\Sigma, \mathcal{A}, \vec{\tau})$ is open and dense in the submanifold $\{(x', \tilde{q}_N^0, v_N^0, m_N = 0) | \, x' \in \mathbb{C}^h\}$ (here $h = (2\nu+1)(N-1)+1$);*
(2) *For a nonempty, open set of orbit segments $\omega' \in \tilde{\Omega}(\Sigma', \mathcal{A}', \vec{\tau}')$ there exists a trajectory segment $\omega \in \tilde{\Omega}(\Sigma, \mathcal{A}, \vec{\tau})$ for which $m_N(\omega) = 0$, $\text{trunc}(\omega) = \omega'$, and $(\tilde{q}_N^0(\omega), v_N^0(\omega)) = (\tilde{q}_N^0, v_N^0)$;*
(3) *For an open and dense set of orbit segments $\omega' \in \tilde{\Omega}(\Sigma', \mathcal{A}', \vec{\tau}')$ there exists a trajectory segment $\omega \in \tilde{\Omega}(\Sigma, \mathcal{A}, \vec{\tau})$ for which $m_N(\omega) = 0$, $\text{trunc}(\omega) = \omega'$, and $(\tilde{q}_N^0(\omega), v_N^0(\omega)) = (\tilde{q}_N^0, v_N^0)$.*

(Here "trunc" denotes the cancellation of all kinetic data of the $N^{\text{th}}$ ball with zero mass.)

*Proof.* Statement (1) is an obvious consequence of Main Lemma 4.21/(I) and the Fubini theorem. The second assertion follows immediately from the definition of the derived schemes; see also Remark 4.12. To demonstrate (3), we select and fix a typical pair $(\tilde{q}_N^0, v_N^0)$ according to (1) of Lemma 4.37. Intersect the open and dense subset of $\mathbb{C}^h$ – featured in (1) – with the open and dense domain $D(\Sigma', \mathcal{A}', \vec{\tau}')$ and obtain the open and dense intersection $D' \subset \mathbb{C}^h$. Suppose now that $\omega' \in \tilde{\Omega}(\Sigma', \mathcal{A}', \vec{\tau}')$, and the vector $x' = x(\omega')$ of the initial variables of $\omega'$ belongs to $D'$. (This is an open and dense condition on $\omega'$.) Launch the $(\Sigma, \mathcal{A}, \vec{\tau})$-dynamics with the initial data $(x', \tilde{q}_N^0, v_N^0, m_N = 0) \in D(\Sigma, \mathcal{A}, \vec{\tau})$ and try to continue the $(\Sigma, \mathcal{A}, \vec{\tau})$-dynamics by also preserving the time evolution of the given $\omega' \in \tilde{\Omega}(\Sigma', \mathcal{A}', \vec{\tau}')$. Remark 4.13 (more precisely, the inheritance of the reducibility proved there) and the definition of the derived schemes (especially (2) thereof) ensure that the continuation of the $(\Sigma, \mathcal{A}, \vec{\tau})$-dynamics is, indeed, possible by the preservation of the given time evolution of $\omega'$. Hence Lemma 4.37 follows. □

The next lemma will use the following:

*Definition* 4.38. Suppose that two indices $1 \le p < q \le n$ and two labels of balls $i, j \in \{1, \ldots, N\}$ are given with the additional requirement that if $i = j$, then $i \in \bigcup_{l=p+1}^{q-1} \sigma_l$. Following the proof of Lemma 4.2, denote by $Q_1(\vec{x}), Q_2(\vec{x}), \ldots, Q_\nu(\vec{x})$ ($\vec{x} \in \mathbb{C}^{(2\nu+1)N+1}$) the polynomials with the property that for every vector of initial data $\vec{x} \in D(\Sigma, \mathcal{A}, \vec{\tau})$ and for every $k$, $k = 1, \ldots, \nu$, the following equivalence holds true:

$$\left(\exists \omega \in \tilde{\Omega} \text{ such that } \vec{x}(\omega) = \vec{x} \text{ and } (v_i^p(\omega))_k = \left(v_j^{q-1}(\omega)\right)_k\right)$$
$$\iff Q_k(\vec{x}) = 0.$$

LEMMA 4.39. *Assume that the combinatorial-algebraic scheme* $(\Sigma, \mathcal{A}, \vec{\tau})$ *has Property* (A), *see Definition* 3.31, *and use the assumptions and notation of the above definition.*

*At least one of the polynomials* $Q_1, \ldots, Q_\nu$ *is nonzero, i.e.* $v_i^p(\omega) \neq v_j^{q-1}(\omega)$ *for almost every* $\omega \in \tilde{\Omega}(\Sigma, \mathcal{A}, \vec{\tau})$, *see also the Corollary* 4.7.

*Proof.* Induction on the number $N \geq 2$.

1. Base of the induction, $N = 2$: First of all, by performing the substitution $L = 0$, we can annihilate all adjustment vectors; see the proofs of parts 4/a,4/b in Main Lemma 4.21, in accordance with statement (VI) of that lemma. Then an elementary inspection shows that for any selection of *positive* masses $(m_1, m_2)$, indeed, $v_i^p(\omega) \neq v_j^{q-1}(\omega)$ almost surely in the section $\tilde{\Omega}(\Sigma, \mathcal{A}, \vec{\tau}, \vec{m})$ of $\tilde{\Omega}(\Sigma, \mathcal{A}, \vec{\tau})$ corresponding to the selected masses, since any trajectory segment of a two-particle system with positive masses and $\mathcal{A} = 0$ has a very nice, totally real (and essentially two-dimensional) representation in the relative coordinates of the particles: the consecutive, elastic bounces of a point particle moving uniformly inside a ball of radius $2r$. Therefore, the statement of the lemma is true for $N = 2$.

We note here that there are merely four points in the entire proof of Key Lemma 4.1 where we need to use Property (A) explicitly: The case $N = 2$, $i = j$ in the proof of this lemma, the Lemmas 4.27, 4.33 and 4.34. Whenever a collision of type $\{1,2\}$ is repeated an even number of times with the same adjustment vector and zero time slots, the resulting velocity change is identically zero, because an even power of the same velocity reflection is applied. (See Example I at the end of this section.)

2. Assume that $N \geq 3$ and the lemma has been proved for all smaller numbers of balls. By re-labeling the particles, if necessary, we can obtain

(i) $N \neq i$, $N \neq j$ and

(ii) if $i = j$, then the ball $i$ has at least one collision between $\sigma_p$ and $\sigma_q$ with a particle different from $N$.

For the fixed combinatorial scheme $(\Sigma, \mathcal{A}, \vec{\tau})$, select a scheme $(\Sigma', \mathcal{A}', \vec{\tau}')$ derived from it. By applying the inductive hypothesis to $(\Sigma', \mathcal{A}', \vec{\tau}')$, we get that for almost every $\omega' \in \tilde{\Omega}(\Sigma', \mathcal{A}', \vec{\tau}')$ the inequality $v_i^p(\omega') \neq v_j^{q-1}(\omega')$ holds. Moreover, by virtue of Lemma 4.37/(1), for almost every $(\tilde{q}_N^0, v_N^0) \in \mathbb{C}^\nu \times \mathbb{C}^\nu$ one has $(x', \tilde{q}_N^0, v_N^0, m_N = 0) \in D(\Sigma, \mathcal{A}, \vec{\tau})$ for a dense, open subset of $x' \in \mathbb{C}^h$.

Let us fix one such pair $(\tilde{q}_N^0, v_N^0)$. By Lemma 4.37/(2), for a nonempty, open set of $\omega' \in \tilde{\Omega}(\Sigma', \mathcal{A}', \vec{\tau}')$ it is true that $(x(\omega'), \tilde{q}_N^0, v_N^0, m_N = 0) \in D(\Sigma, \mathcal{A}, \vec{\tau})$, $v_i^p(\omega') \neq v_j^{q-1}(\omega')$, and there exists an $\omega \in \tilde{\Omega}(\Sigma, \mathcal{A}, \vec{\tau})$ with $m_N(\omega) = 0$ and $\text{trunc}(\omega) = \omega'$.

The assertion of the lemma now follows for the $N$-ball systems since the bare existence of a single orbit segment with the above nonequality implies that not all polynomials $Q_l(\vec{x})$ are zero. $\square$

LEMMA 4.40. *Assume that the combinatorial scheme $(\Sigma, \mathcal{A}, \vec{\tau})$ satisfying Property (A) is given for the $N$ ($N \geq 3$) ball system $\{1, 2, \ldots, N\}$ and, furthermore, there are two integers $1 \leq p < q \leq n$ such that $N \in \sigma_p \cap \sigma_q$, $N \notin \sigma_j$ for $p < j < q$, and if $\sigma_p = \sigma_q$ then we require additionally that there exists an index $j$, $p < j < q$, such that $\sigma_p \cap \sigma_j \neq \emptyset$. Denote by $(\Sigma', \mathcal{A}', \vec{\tau}')$ an arbitrary derived scheme of $(\Sigma, \mathcal{A}, \vec{\tau})$.*

*Denote by $P_1(\vec{x}), \ldots, P_s(\vec{x})$ the canonically determined complex polynomials, the simultaneous vanishing of which is the consequence of nonsufficiency of the $(\Sigma, \mathcal{A}, \vec{\tau})$ orbit segments in the sense of Lemma 4.2. If all these polynomials are zero, then all the analogous polynomials associated with the derived scheme $(\Sigma', \mathcal{A}', \vec{\tau}')$ are also zero.*

*In other words, if all $(\Sigma, \mathcal{A}, \vec{\tau})$-orbit segments are nonsufficient, then the same is true for all $(\Sigma', \mathcal{A}', \vec{\tau}')$-orbit segments; see also the corollary at the end of the proof of Lemma 4.2.*

*Proof.* Assume that the assertion is false for the $(\Sigma', \mathcal{A}', \vec{\tau}')$-orbit segments, that is, there exists at least one sufficient $(\Sigma', \mathcal{A}', \vec{\tau}')$-orbit segment or, equivalently, by Lemma 4.8, almost all $(\Sigma', \mathcal{A}', \vec{\tau}')$-orbit segments are sufficient.

By putting together Lemmas 4.37/(1,2), 4.39 and 4.8, we obtain a sufficient $(\Sigma, \mathcal{A}, \vec{\tau})$-orbit. Indeed, for typical elements $x' \in D(\Sigma', \mathcal{A}', \vec{\tau}')$ we have that

$$\text{(4.41)} \qquad \text{for all } \omega' \in \tilde{\Omega}(\Sigma', \mathcal{A}', \vec{\tau}') \text{ with } x(\omega') = x'$$
$$\omega' \text{ is sufficient and } v_{i_p}^p(\omega') \neq v_{i_q}^{q-1}(\omega').$$

By Lemma 4.37/(1) for almost every $(\tilde{q}_N^0, v_N^0) \in \mathbb{C}^\nu \times \mathbb{C}^\nu$ it is true that for a dense open subset of $x' \in \mathbb{C}^h$

$$\text{(4.42)} \qquad (x', \tilde{q}_N^0, v_N^0, m_N = 0) \in D(\Sigma, \mathcal{A}, \vec{\tau}).$$

Thus, for almost every $(\tilde{q}_N^0, v_N^0) \in \mathbb{C}^\nu \times \mathbb{C}^\nu$ we have the following: For a dense, open subset of $x' \in \mathbb{C}^h$ (4.41) and (4.42) hold. Further, the velocity $v_N^{q-1}(\omega) = v_N^p(\omega)$ can be chosen in such a *typical* manner that it is not contained by the complex line

$$v_{i_p}^p(\omega') + \lambda \cdot \left(v_{i_q}^{q-1}(\omega') - v_{i_p}^p(\omega')\right) \quad (\lambda \in \mathbb{C}),$$

so that condition (2) of Lemma 4.8 is satisfied. (The possibility of the desired choice for $v_N^p(\omega) = v_N^{q-1}(\omega)$ is easily seen if we represent our variables at the $p^{\text{th}}$ collison; cf. the proof of Lemma 3.19.)

By Lemma 4.37/(1, 2), all this implies that for a typical $(\tilde{q}_N^0, v_N^0) \in \mathbb{C}^\nu \times \mathbb{C}^\nu$ there exists a nonempty, open set $G$ of orbit segments $\omega' \in \tilde{\Omega}(\Sigma', \mathcal{A}', \vec{\tau}')$, with $(x(\omega'), \tilde{q}_N^0, v_N^0, m_N = 0) \in D(\Sigma, \mathcal{A}, \vec{\tau})$, such that for every $\omega' \in G$ there exists an $\omega \in \tilde{\Omega}(\Sigma, \mathcal{A}, \vec{\tau})$ satisfying $m_N(\omega) = 0$, $\text{trunc}(\omega) = \omega'$, $(\tilde{q}_N^0(\omega), v_N^0(\omega)) = (\tilde{q}_N^0, v_N^0)$, and, according to Lemma 4.8, $\omega$ is sufficient. Hence Lemma 4.40 follows. □

The next, purely combinatorial lemma is the last ingredient of the inductive proof of the key lemma.

LEMMA 4.43. *Define the sequence of positive numbers $C(N)$ recursively by taking $C(2) = 1$ and $C(N) = \frac{N}{2} \cdot \max\{C(N-1); 3\}$ for $N \geq 3$. Let $N \geq 3$, and suppose that the symbolic collision sequence $\Sigma = (\sigma_1, \ldots, \sigma_n)$ for $N$ particles is $C(N)$-rich. Then there is a particle, say the one with label $N$, and two indices $1 \leq p < q \leq n$ such that*

(i) $N \in \sigma_p \cap \sigma_q$,
(ii) $N \notin \bigcup_{j=p+1}^{q-1} \sigma_j$,
(iii) $\sigma_p = \sigma_q \implies \exists j\ p < j < q$ and $\sigma_p \cap \sigma_j \neq \emptyset$, and
(iv) $\Sigma'$ is $C(N-1)$-rich on the vertex set $\{1, \ldots, N-1\}$.

Here, just as in the case of derived schemes, $\Sigma'$ is the symbolic sequence obtained from $\Sigma$ by discarding all edges containing $N$.

*Proof.* The hypothesis on $\Sigma$ implies that there exist subsequences $\Sigma_1, \ldots, \Sigma_r$ of $\Sigma$ with the following properties:

(1) For $1 \leq i < j \leq r$ every collision of $\Sigma_i$ precedes every collision of $\Sigma_j$.

(2) The graph of $\Sigma_i$ $(1 \leq i \leq r)$ is a tree (a connected graph without loop) on the vertex set $\{1, \ldots, N\}$, and

(3) $r \geq C(N)$.

Since every tree contains at least two vertices with degree one and $C(N) = \frac{N}{2} \cdot \max\{C(N-1); 3\} = \frac{N}{2} \cdot b$ ($b$ is shorthand for $\max\{C(N-1); 3\}$), there is a vertex, say the one labelled by $N$, such that $N$ is a degree-one vertex of $\Sigma_{i(1)}, \ldots, \Sigma_{i(t)}$, where $1 \leq i(1) < \cdots < i(t) \leq r$ and $t \geq b$.

Now $t \geq C(N-1)$ and, therefore, (iv) obviously holds.

Let $\sigma_{p'}$ be the edge of $\Sigma_{i(1)}$ that contains $N$ and, similarly, let $\sigma_{q'}$ be the edge of $\Sigma_{i(t)}$ containing the vertex $N$. Then the fact that $t \geq 3$ ensures that the following properties hold:

(i)′ $N \in \sigma_{p'} \cap \sigma_{q'}$,
(iii)′ $\sigma_{p'} = \sigma_{q'} \implies \exists j,\ p' < j < q'$ and $\sigma_{p'} \cap \sigma_j \neq \emptyset$, $\sigma_j \neq \sigma_{p'}$.

Let $\sigma_p$, $\sigma_q$ ($1 \leq p < q \leq n$) be a pair of edges $\sigma_{p'}$, $\sigma_{q'}$ ($1 \leq p' < q' \leq n$) fulfilling (i)' and (iii)' and having the minimum possible value of $q' - p'$. Elementary inspection shows that then (ii) must also hold for $\sigma_p$, $\sigma_q$. Lemma 4.43 is now proved. $\square$

We are now able to prove the key lemma by induction on the number $N \geq 2$. Indeed, for $N = 2$ there are no nonsufficient trajectories with initial data in $D(\Sigma, \mathcal{A}, \vec{\tau})$ and, therefore, the assertion of the key lemma is obviously true.

The inductive step $(N-1) \to N$ can be obtained by simply putting together Lemmas 4.2, 4.8, 4.40 and 4.43. The proof of the key lemma is now complete in the complex case. $\square$

Finally, since a nonzero polynomial $P_i(\vec{x})$ ($1 \leq i \leq s$) takes nonzero values almost everywhere on the real space $\mathbb{R}^{(2\nu+1)N+1}$, we immediately obtain the validity of the real version of Key Lemma 4.1. $\square$

We conclude this paragraph by presenting to the reader two interesting examples. The first of them shows the absolute necessity of imposing Property (A) on the discrete structure $(\Sigma, \mathcal{A}, \vec{\tau})$, see Definition 3.31, while the second one sheds light on the limitations of our algebraic approach developed in Sections 3–4.

*Example* I. Suppose that the collision graph of $\Sigma$ is a tree (a connected graph without loop) with the possiblity of repeated edges. (This allows an arbitrarily large number of consecutive, connected collision subgraphs.) Assume further that the combinatorial-algebraic scheme $(\Sigma, \mathcal{A}, \vec{\tau})$ is given in such a way that it has the following property: There exists a sequence $1 = k_1 < k_2 < \cdots < k_r < k_{r+1} = n + 1$ such that

(1) $k_{j+1} - k_j$ is an even number for $j = 1, \ldots, r$, and
(2) $\sigma_p = \sigma_{k_j}$, $a_p = a_{k_j}$, and $\tau_p = 0$ for $k_j < p < k_{j+1}$, $j = 1, \ldots, r$.

Plainly, the orthogonal velocity reflections corresponding to the collisions

$$\sigma_{k_j}, \sigma_{k_j+1}, \ldots, \sigma_{k_{j+1}-1}$$

are the same ($j = 1, \ldots, r$) and, therefore, their product is the identity operator. Thus $v_i^0(\omega) = v_i^{k_j-1}(\omega)$ for all $i = 1, \ldots, N$, $j = 1, \ldots, r+1$, $\omega \in \tilde{\Omega}(\Sigma, \mathcal{A}, \vec{\tau})$. It is an easy exercise then to show that these orbit segments $\omega$ are very far from being sufficient. Actually, $\dim_{\mathbb{C}}\{\alpha_1, \ldots, \alpha_n\} = N - 1$ and $\dim_{\mathbb{C}} \mathcal{N}(\omega) = \nu + N - 1$; see also (4.4).

*Example* II. Suppose now that the $C(N)$-rich symbolic sequence $\Sigma$ only contains the collisions of type $\{i, i+1\}$ for $i = 1, \ldots, N-1$; i.e. the collision graph is the simple path of length $N$ with the allowance of repeated edges. Consider the set $X$ of orbit segments $\omega \in \tilde{\Omega}(\Sigma, \mathcal{A}, \vec{\tau})$ with the side conditions $m_1 = m_3 = 0$. (Of course, $m_2 \neq 0$.) It is easy to see that the orbit segments from $X$ are *never* sufficient; the vectors

$$(\delta \tilde{q}_1^0, \ldots, \delta \tilde{q}_N^0) = \left(\lambda \cdot (v_1^0 - v_2^0), 0, \ldots, 0\right)$$

$(0 \neq \lambda \in \mathbb{C})$ are nontrivial neutral vectors. The interesting feature of this example is that the set $X$ has two codimensions in $\tilde{\Omega}(\Sigma, \mathcal{A}, \vec{\tau})$, and the codimension of $X$ decreases to one at the crucial step of the induction when we perform the substitution $m_1 = 0$. (In earlier notation $m_N = 0$.)

## 5. Proof of the Main Theorem

Richness will be understood throughout in the sense of Definition 2.5 by choosing $C := C(N)$ as prescribed by Key Lemma 4.1 (see Remark 4.1/b after the key lemma).

According to the results of [V(1979)], [G(1981)] and [B-F-K(1998)], in a semi-dispersing billiard system satisfying the nondegeneracy condition mentioned before our Main Theorem, there are no orbits with a finite accumulation point of collision moments. Thus, it is sufficient to concentrate our attention on the set $M^0$ of orbits containing no singular collision, since it is well-known and easy to see that phase points with a singular orbit form a countable union of proper submanifolds.

THEOREM 5.1. *Consider a system of $N$ ($\geq 3$) particles on the $L$-torus $\mathbb{T}_L^\nu$ ($\nu \geq 2$) satisfying $r \in R_0$. Let $P = \{P_1, P_2\}$ be a given, two–class partition of the $N$ particles, where, for simplicity, $P_1 = \{1, \ldots, n\}$ and $P_2 = \{n+1, \ldots, N\}$ ($n < N - 1$). Then the closed set*

$$F_+ = \left\{x \in \mathbf{M} : S^{[0,\infty)}x \text{ is partitioned by } P\right\}$$

*has measure zero.*

This statement would, in fact, allow us to fix $C(N)$ (the number of consecutive, connected collision graphs required for sufficiency) arbitrarily large.

*Proof.* Our reasoning is reminiscent of the proof of Theorem 5.1 of 2 [Sim(1992)-I] irrespective of the fact that the notion of center of mass has lost its sense now. The two cases $\min\{n, N - n\} \geq 2$ and $\min\{n, N - n\} = 1$ can be treated similarly, and thus we only consider the first one.

Every point $x \in \mathbf{M}$ can be characterized by the following coordinates in an essentially unique way:

(1)  $\pi_{P_1}(x) = x_1 \in \mathbf{M}_1$,           (2)  $\pi_{P_2}(x) = x_2 \in \mathbf{M}_2$,

(3)  $C_1(x) = \tilde{q}_{n+1} - \tilde{q}_1 \in \mathbb{R}^\nu$,   (4)  $\dfrac{I_1(x)}{\|I_1(x)\|} \in \mathbb{S}^{\nu-1}$,

(5)  $\|I_1(x)\| \in \mathbb{R}_+$,         (6)  $E_1(x) = \dfrac{1}{2} \sum_{i=1}^n m_i v_i^2(x) \in \mathbb{R}_+$,

where $I_1(x) = \dfrac{\sum_{i=1}^n m_i v_i(x)}{\sum_{i=1}^n m_i}$. (Non-uniqueness only arises in choosing $C_1(x)$ as an arbitrary representative of $q_{n+1} - q_1 \in \mathbb{T}_L^\nu$.) In what follows the six coordinates corresponding to the characterization given above will, in general, be denoted by $b_1, \ldots, b_6$; thus it will always be assumed that $b_1 \in \mathbf{M}_1, b_2 \in \mathbf{M}_2, b_3 \in \mathbb{R}^\nu, b_4 \in \mathbb{S}^{\nu-1}, b_5, b_6 \in \mathbb{R}_+$. The relation $\mu(F_+) = 0$ will certainly follow if we show that for almost every such choice of the $b_i$'s

(5.2) $$\mu_{b_1,b_2,b_3,b_5,b_6}(F_+(b_1, b_2, b_3, b_5, b_6)) = 0,$$

where

$$F_+(b_1, b_2, b_3, b_5, b_6) = \{y_4 \in \mathbb{S}^{\nu-1}: \ (b_1, b_2, b_3, y_4, b_5, b_6) \in F_+\},$$

and $\mu_{b_1,b_2,b_3,b_5,b_6}$ denotes the conditional measure of $\mu$ under the conditions corresponding to fixing the values of $b_1, b_2, b_3, b_5, b_6$. (This conditional measure is equivalent to the Lebesgue measure on its support.)

The relation $x \in F_+$ is equivalent to saying that for every pair $i \in P_1, j \in P_2$ and every $t \geq 0$

(5.3) $$\varrho\left(q_i^t(x) - q_j^t(x), L \cdot \mathbb{Z}^\nu\right) \geq 2r,$$

where $\varrho(.,.)$ denotes the euclidean distance. For simplicity, fix $i = 1$ and $j = n+1$.

Now we will consider the subdynamics corresponding to our two-class partition (cf. the "Subsystems, decompositions" part of Section 2 in the paper [Sim(1992)-I]), and will denote them, for simplicity, by $S_1$ and $S_2$, respectively (their phase spaces are $\mathbf{M}_1$ and $\mathbf{M}_2$, of course).

It is worthwhile to note here that — according to the construction, given in the introduction, of the configuration space of our hard ball system with an arbitrary mass vector — we cannot talk about the absolute configuration vector of a particle but rather about the relative position $q_i^t(x) - q_j^t(x)$ for any pair of particles, only. Also, it makes sense to talk about the time-displacement $q_i^t(x) - q_i^0(x)$ for any $1 \leq i \leq N$, $t \in \mathbb{R}$ and $x \in \mathbf{M}$ since

$$q_i^t(x) - q_i^0(x) = \int_0^t v_i^s(x)ds \qquad (\mathrm{mod}\ L \cdot \mathbb{Z}^\nu).$$

Our remarks naturally apply to the subdynamics $S_1$ and $S_2$ as well.

Now, assuming for simplicity $1 < n < N-1$, a standard calculation, which is given for instance in the proof of Lemma 5.8 of [Sim(1992)-I], yields that

$$q_{n+1}^t(x) - q_1^t(x) = \int_0^{\beta t} \tilde{v}_{n+1}^s(x_2)ds - \int_0^{\alpha t} \tilde{v}_1^s(x_1)ds + C_1(x) + \tilde{I}t \quad (\text{mod } L \cdot \mathbb{Z}^\nu) \tag{5.4}$$

where $\tilde{v}_1^s(x_1)$ and $\tilde{v}_{n+1}^s(x_2)$ denote the time evolution of the corresponding velocity vectors under the subdynamics $S_1$ and $S_2$, respectively, and $\tilde{I}$ is the relative velocity of the "baricenters" of the second and first subsystems. This term appears since in $\mathbf{M}_1$ and $\mathbf{M}_2$ the moments of the subsytems are scaled to be equal to 0; actually,

$$\tilde{I} = \frac{-\mathcal{M}_1 - \mathcal{M}_2}{\mathcal{M}_2} I_1 = \frac{-\mathcal{M}_1 - \mathcal{M}_2}{\mathcal{M}_2} \|I_1\| \frac{I_1}{\|I_1\|} \tag{5.5}$$

where

$$\mathcal{M}_1 = \sum_{i=1}^n m_i \quad \text{and} \quad \mathcal{M}_2 = \sum_{i=n+1}^N m_i$$

are the masses of the subsystems. Finally,

$$\alpha = \sqrt{2E_1(x) - \mathcal{M}_1 \|I_1(x)\|^2}, \quad \text{and}$$

$$\beta = \sqrt{1 - 2E_1(x) - \frac{\mathcal{M}_1^2}{\mathcal{M}_2} \|I_1(x)\|^2}$$

are the corresponding time scalings.

Our task is to show that the event

$$\text{for all } t \geq 0, \ \varrho\left(\lambda t \frac{I_1}{\|I_1\|} + f(t), L \cdot \mathbb{Z}^\nu\right) \geq 2r \tag{5.6}$$

has measure zero for every fixed $b_1, b_2, b_3, b_5, b_6$ where

$$\lambda = (\mathcal{M}_2)^{-1}(-\mathcal{M}_1 - \mathcal{M}_2)\|I_1\|$$

and

$$f(t) : \mathbb{R} \to \mathbb{R}^\nu$$

is an arbitrary fixed function such that $f(t) = f(t, x_1, x_2, C_1, \|I_1\|, E_1)$. Actually, by (5.4),

$$f(t) = \int_0^{\beta t} \tilde{v}_{n+1}^s(x_2)ds - \int_0^{\alpha t} \tilde{v}_1^s(x_1)ds + C_1(x).$$

In (5.6), the canonical meaning of zero measure is that in $\frac{I_1}{\|I_1\|}$.

Denote by $\mathcal{L}_{2r,L}$ the lattice of balls of radius $2r$ centered at points of $L \cdot \mathbb{Z}^\nu$. Our proof of Theorem 5.1 will be based on the following well-known elementary lemma.

LEMMA 5.7. *Fix a vector $\vec{n} \in \mathbb{S}^{\nu-1}$ for which at least one ratio of coordinates is irrational. Consider arbitrary hyperplanes $H$ perpendicular to $\vec{n}$, and denote by $B_R(z)$ the $(\nu-1)$-dimensional ball of radius $R$ in $H$ centered at $z \in H$. Then, for a suitable $\gamma(\vec{n}) > 0$*

$$(5.8) \qquad \liminf_{R \to \infty} \inf_H \inf_{z \in H} \frac{\operatorname{meas}(B_R(z) \cap \mathcal{L}_{2r,L})}{\operatorname{meas}(B_R(z))} \geq \gamma(\vec{n}).$$

*Proof.* The lemma easily follows from the fact that the set $B_R(z)$ modulo $L \cdot \mathbb{Z}^\nu$ is "uniformly dense" in the torus $\mathbb{T}_L^\nu$; i.e. it is an $\varepsilon(R)$-dense set in $\mathbb{T}_L^\nu$, where $\varepsilon(R) \to 0$ as $R \to \infty$. This statement is, however, an easy consequence of the well-known density of the hyperplane $H$ (mod $L \cdot \mathbb{Z}^\nu$) in the torus $\mathbb{T}_L^\nu = \mathbb{R}^\nu / L \cdot \mathbb{Z}^\nu$ and the compactness of $\mathbb{T}_L^\nu$. Indeed, since being an $\varepsilon$-dense set is a translation invariant property, the positioning of the foot point $z$ is irrelevant. For a given positive $\varepsilon$, we choose first a finite $\varepsilon/2$-dense set $F \subset \mathbb{T}^\nu$, and then for every point $y \in F$ we select a radius $R_y > 0$ with the property that $B_{R_y}(z)$ is closer to $y$ than $\varepsilon/2$. The maximum value $R(\varepsilon)$ of the finitely many radii $R_y$ ($y \in F$) entering the problem provides now an $\varepsilon$-dense set $B_{R(\varepsilon)}(z)$. Now select and fix $\varepsilon = r$. Observe that at least a fixed, positive percentage of the ball $B_{R(r)}(z)$ belongs to $\mathcal{L}_{2r,L}$. However, for every large enough $R$, also a fixed, positive percentage of the ball $B_R(z)$ can be covered by disjoint balls $B_{R(r)}(z_i)$. Hence Lemma 5.7 follows. □

Assume that statement (5.6) is not true; i.e. the measure of the subset $K$ of $\mathbb{S}^{\nu-1}$ described by (5.6) is positive. Select then and fix a Lebesgue density point $\vec{n}$ of $K$ with the property that at least one ratio of the components of $\vec{n}$ is irrational. Denote by $G_\varepsilon \subset \mathbb{S}^{\nu-1}$ the ball of radius $\varepsilon$ around $\vec{n}$. By (5.8) we can choose $R_0$ so large that for $R \geq R_0$

$$\inf_H \inf_{z \in H} \frac{\operatorname{meas}(B_R(z) \cap \mathcal{L}_{2r,L})}{\operatorname{meas}(B_R(z))} \geq \frac{\gamma(\vec{n})}{2}.$$

The set $\lambda t G_\varepsilon$ can be arbitrarily well approximated by a ball of radius $\lambda t \varepsilon = R$ ($R$ is fixed, $R \geq R_0$) in the hyperplane orthogonal to $\vec{n}$ through the point $\lambda t \vec{n} + f(t)$ if only $t$ is sufficiently large. Consequently, if $R \geq R_0$, then by choosing $t$ sufficiently large and at the same time putting $\varepsilon = (\lambda t)^{-1} R$, we have

$$\frac{\operatorname{meas}((\lambda t G_\varepsilon + f(t)) \cap \mathcal{L}_{2r,L})}{\operatorname{meas}((\lambda t G_\varepsilon + f(t)))} \geq \frac{\gamma(\vec{n})}{4}.$$

But this inequality contradicts to the fact that $\vec{n}$ was chosen as a Lebesgue density point of the subset $K \subset \mathbb{S}^{\nu-1}$. Hence Theorem 5.1 follows. □

Return now to the demonstration of the Main Theorem, restricting attention to the subset $M^0$. From Theorem 5.1 it follows that almost every point is $C(N)$-rich whatever the value of $C(N)$ may be. Let us fix $C(N)$ according to the requirement of Key Lemma 4.1. Denote by $\hat{R} \subset M^0$ the subset of $C(N)$-rich phase points. By the real version of Key Lemma 4.1, apart from a countable union of proper analytic submanifolds of the outer geometric parameters $(\vec{m}, L)$, for every $(\Sigma, \mathcal{A}, \vec{\tau})$ satisfying Property (A), almost every nonsingular orbit segment with the combinatorial scheme $(\Sigma, \mathcal{A}, \vec{\tau})$ is sufficient. Since the number of combinatorial schemes is countable, and these schemes of actual (real) hard ball trajectories necessarily satisfy Property (A), we can conclude that, indeed, almost every orbit in $\hat{R}$ is sufficient.

At this point we will use the following simple consequence of Poincaré recurrence and the ergodic theorem.

LEMMA 5.9 (see [S-Ch(1987)]).   *If $S^{[a,b]}x$ is a nonsingular orbit segment and it is sufficient, then, in a suitable open neighbourhood of $x$, all relevant Lyapunov exponents are nonzero almost everywhere.*

Adding Lemma 5.9 to what has been said before, one also obtains that for almost every phase point in $\hat{R}$, the relevant Lyapunov exponents of the flow are nonzero. Thus our Main Theorem is proved.  □

*Proof of the Corollary to the Main Theorem (Section 1).* The positivity of the ergodic components is a consequence of the Katok-Strelcyn theory (cf. [K-S(1986)]), where the underlying assumption is that the relevant Lyapunov-exponents of the system do not vanish. (We note that, in fact, the Katok-Strelcyn theory is a generalization of Pesin's theory (see [P(1977)] for hyperbolic systems with singularities, also containing billiards isomorphic to hard ball systems).

Finally, it is worthwhile to note here that, in fact, on any positive ergodic component (i.e. on any ergodic component of positive measure) the Bernoulli property of the standard billiard flow $(\mathbf{M}, \{S^t\}, \mu)_{\vec{m}, L}$ follows from its hyperbolicity and weak mixing, as has been recently shown by Chernov-Haskell and Ornstein-Weiss; cf. [C-H(1996)] and [O-W(1998)]. For a brief summary of how the K-property follows from ergodicity and hyperbolicity, see, for instance, Section 7 of [Sim(1992)-I], or [K-S(1986)].  □

## 6. Concluding remarks

1. Had we been able to strengthen the statement of Key Lemma 4.1 from "*for almost every $\omega \in \tilde{\Omega}(\Sigma, \mathcal{A}, \vec{\tau})$*", i.e. from "for every $\omega \in \tilde{\Omega}(\Sigma, \mathcal{A}, \vec{\tau})$ apart from a closed algebraic subset of codimension one", to a similar statement with an *at least two-codimensional algebraic subset*, we would have obtained the global ergodicity of the system. A crucial step in getting to codimension two in the Key Lemma would be the strengthening of the last claim of Lemma 4.39 to stating $v_i^p(\omega) \neq v_j^{q-1}(\omega)$ apart from a subset of codimension two. Both lemmas would be implied by the following analogous algebraic properties: All common divisors of the polynomials

(1) $P_1, \ldots, P_s$    (in case of Key Lemma 4.1);
(2) $Q_1, \ldots, Q_\nu$    (in case of Lemma 4.39)

are purely homogeneous polynomials of the masses (example II at the end of Section 4 illustrates the essential difficulty arising in the inductive proof of this claim). In this way one can formulate a *hypothetical statement about the Boltzmann-Sinai ergodic hypothesis*: *If this stronger form of Key Lemma* 4.1 *is true, then, in the Main Theorem, we can also claim the* K-*property instead of just full hyperbolicity of the system.*

2. The system of $N = 2$ hard discs (i.e. $\nu = 2$) with arbitrary masses $m_1$, $m_2$ was studied by Simányi and Wojtkowski [S-W(1989)]. They considered the same dynamics as here, but in the larger phase space $\tilde{\mathbf{Q}} \times S^1$. By using the fact that this system is an isometric $\mathbb{T}^2$-extension of the standard billiard ball flow $(\mathbf{M}, S^{\mathbb{R}}, \mu)_{(m_1, m_2)}$, which is a planar dispersing billiard, and was known to possess the Bernoulli property since [G-O(1974)], they showed that this extension is also a Bernoulli flow as long as the ratio $m_1/m_2$ is irrational. If one thoroughly studies the steps of that proof, he realizes that, once the difficulties outlined in the above remark have been resolved, the methods of the present paper and [S-W(1989)] together with the results of [C-H(1996)] and [O-W(1998)] prove even the Bernoulli property of the flow on the extended phase space $\tilde{\mathbf{Q}} \times \mathcal{E}$ (without the factorization $\Psi$; see the introduction), provided that *not all ratios $m_i/m_j$ are rational.* (The latter condition means precisely that the "center of mass" is defined modulo a dense subgroup of $\mathbb{T}^\nu$.)

3. In [Sz(1994)], the second author obtained necessary and sufficient conditions for the K-property of the so-called "orthogonal cylindric billiards." The subclass of cylindric billiards, introduced in [Sz(1993)], within the family of semi-dispersing billiards, is of particular interest since it is this class for which the formulation of general — and at the same time in some sense constructively verifiable — necessary and sufficient conditions of the K-property seems possible. (We have seen in the proof of Lemma 2.1 that the system of $N$ balls with an arbitrary mass vector is also isomorphic to a cylindric billiard.) Indeed, such a condition is found and shown to be necessary for the systems of $N$ hard balls of identical masses in a forthcoming paper of the present authors [S-Sz(1997)]. One apparent difficulty in establishing the sufficiency of this condition is the fact that so far no characterization of the neutral subspace of a trajectory segment of a general cylindric billiard is known (in the case of hard ball systems such a characterization was obtained by using the CPF, whose derivation, on the other hand, used the conservation of momentum, a property present in hard ball systems but absent in general cylindric billiards).

*Acknowledgement.* The authors are deeply indebted to N. Chernov for the large number of very valuable remarks and corrections he made during a careful reading of the manuscript. Special thanks are due to János Kollár who called our attention to the advantage of the complexification of the dynamics which turned out to be very important in Sections 3 and 4. Also, he made useful remarks concerning some elementary facts from algebraic geometry. The authors express their sincere gratitude to A. Katok, C. Liverani, Ya. Pesin and L. Vaserstein for their helpful remarks on earlier versions of the paper.

Bolyai Institute, University of Szeged, Szeged, Hungary
*E-mail address*: simanyi@sol.cc.u-szeged.hu
Mathematical Institute of the Hungarian Academy of Sciences, Budapest, Hungary
*E-mail address*: szasz@math-inst.hu